\documentstyle[
11pt,twoside,amsfonts,amssymb,amsmath,amsthm,latexsym]{article}
{\tiny }
\newcommand {\rel} {{\mathbb R}}
\newcommand {\com} {{\mathbb C}}

\newcommand {\nat} {{\mathbb N}}

\newcommand {\sphere} {{\mathbb S}}
\newcommand {\stereo} {{\mathcal{S} }}
\newcommand {\Will} {{\mathcal{W} }}
\newcommand {\Wil} {{\mathcal{E} }}

\newcommand {\Lift} {{\mathcal{L} }}
\newcommand {\UU} {{\mathcal{U} }}

\newcommand {\dom} {{\mathcal{D} }}
\newcommand {\OO} {{\mathcal{O} }}
\newcommand {\PP} {{\mathcal{P} }}
\newcommand {\Tan} {{\mathcal{T} }}
\newcommand {\FF} {{\mathcal{F} }}

\newcommand {\CC} {{\mathcal{C} }}
\newcommand {\VV} {{\mathcal{V} }}

\newcommand {\HH} {{\mathcal{H} }}
\newcommand {\Mill} {{\mathcal{M} }} 
\newcommand {\NN} {{\mathcal{N} }}

\begin{document}
\thanks{Institute: Mathematics Department, 
		Technion: Israel Institute of Technology,
		3200003 Haifa, Israel.\\
		Phone: 00972 549298486.
		E-mail: rubenj at technion.ac.il\\\\
	The author was partially funded by the Ministry of Absorption of the State of Israel in the academic years 2019/2020 and 2020/2021.}
	
	\newtheorem{theorem}{Theorem}[section]
	\newtheorem{definition}{Definition}[section]
	\newtheorem{proposition}{Proposition}[section]
	\newtheorem{lemma}{Lemma}[section]
	\newtheorem{corollary}{Corollary}[section]
	\newtheorem{remark}{Remark}[section]
	\newtheorem{example}{Example}
	
	\author{Ruben Jakob}
	
	\title{Global existence and full convergence of the M\"obius-invariant Willmore flow in the $3$-sphere}
	
	\maketitle
	\begin{abstract}
		In this article, we prove two ``global existence and full  convergence theorems'' for flow lines of the M\"obius-invariant Willmore flow, and we use these results, in order to prove
		that fully and smoothly convergent flow lines of the M\"obius-invariant Willmore flow are stable w.r.t. small perturbations of their initial immersions in 
		any $C^{4,\gamma}$-norm, provided they converge either to a 
		smooth parametrization of ``a Clifford-torus'' in $\sphere^3$ or to a umbilic-free $C^4$-local minimizer of the Willmore functional among $C^4$-immersions of a compact torus 
		$\Sigma$ into either $\rel^3$ or $\sphere^3$. The proofs of our four main theorems rely 
		on the author's recent achievements about the M\"obius-invariant Willmore flow, on Escher's, Mayer's and Simonett's work from ``the 90s'' on ``invariant center manifolds'' for uniformly parabolic quasilinear evolution equations and their special applications to the ``Willmore flow'' and ``Surface diffusion flow'' near round $2$-spheres in $\rel^3$ and on Rivi\`ere's and Bernard's fundamental investigation of the Willmore functional on the basis of its conformal invariance and Noether's Theorem.     
	\end{abstract}
\noindent
	MSC-codes: 34C45, 35K46, 35R01, 53C42, 58J35
	
\section{Introduction and main results}  \label{Introduction}
	
In this article, the author continues his investigation
of the M\"obius-invariant Willmore flow (MIWF), which is 
the evolution equation 
\begin{equation}  \label{Moebius.flow}
	\partial_t f_t = -\frac{1}{2}\, \frac{1}{|A^0_{f_t}|^4} \,
	\Big{(} \triangle_{f_t}^{\perp} \vec H_{f_t} + Q(A^{0}_{f_t})(\vec H_{f_t}) \Big{)}
	\equiv -\frac{1}{|A^0_{f_t}|^4} \,\nabla_{L^2} \Will(f_t),
\end{equation}
extending the short time existence result of his article \cite{Jakob_Moebius_2016}.
As explained in \cite{Jakob_Moebius_2016} and \cite{Ruben.MIWF.III}, the flow (\ref{Moebius.flow}) is well-defined for differentiable families of $W^{4,2}$-immersions $f_t$ mapping some arbitrarily fixed smooth compact torus $\Sigma$ into $\sphere^3$, without any umbilic points, and $\Will$ denotes the Willmore-functional
\begin{equation} \label{Willmore.functional}
	\Will(f):= \int_{\Sigma} K^M_f + 
	\frac{1}{4} \, \mid \vec H_f \mid^2 \, d\mu_f,
\end{equation}
which can more generally be considered on $C^{2}$-immersions
$f:\Sigma \longrightarrow M$, mapping any closed smooth Riemannian orientable surface $\Sigma$ into an arbitrary 
smooth Riemannian manifold $M$, where $K^M_f(x)$ denotes the
sectional curvature of $M$ w.r.t. the ``immersed tangent plane''
$Df_x(T_x\Sigma)$ in $T_{f(x)}M$. In the cases relevant in
this article, we have $K_f\equiv 0$ for $M=\rel^n$ or 
$K_f\equiv 1$ for $M=\sphere^{n}$. Regarding the aims of this 
article, we will only have to consider the cases
$M=\rel^3$ or $M=\sphere^3$. Given some immersion 
$f:\Sigma \longrightarrow M$, we endow the torus 
$\Sigma$ with the pullback
$f^*g_{\textnormal{euc}}$ under $f$ of the Euclidean 
metric of either $\rel^3$ or $\rel^{4}$, i.e. with coefficients 
$g_{ij}:=\langle \partial_i f, \partial_j f \rangle$,
and we let $(A_f)_{\rel^{3}}$ and $(A_f)_{\sphere^{3}}$ denote the second fundamental forms of the immersion $f$, either mapping
into $\rel^{3}$ or into $\sphere^{3}$, defined on pairs of 
tangent vector fields $X,Y$ on $\Sigma$ by:
\begin{eqnarray*}  \label{second.fundam.form}
	(A_{f})_{\rel^3}(X,Y) := D_X(D_Y(f)) - P^{\textnormal{Tan}(f),\rel^3}(D_X(D_Y(f)))
	\equiv (D_X(D_Y(f)))^{\perp_{f,\rel^3}}                                            \\
	(A_{f})_{\sphere^3}(X,Y) := D_X(D_Y(f)) - P^{\textnormal{Tan}(f),\sphere^3}(D_X(D_Y(f)))
	\equiv (D_X(D_Y(f)))^{\perp_{f,\sphere^3}}                                    \nonumber
\end{eqnarray*}
where $D_X(V)\lfloor_{x}$ denotes the projection of the classical derivative of a vector field $V:\Sigma \longrightarrow \rel^{3}$ respectively $V:\Sigma \longrightarrow \rel^{4}$
in direction of the tangent vector field $X \in \Gamma(T\Sigma)$ into the respective fiber $T_{f(x)}\rel^3=\rel^3$ of 
$T\rel^3$ respectively $T_{f(x)}\sphere^3$ of $T\sphere^3$, 
and where
\begin{eqnarray*}  \label{Projections}
	P^{\textnormal{Tan}(f),\rel^{3}}:\bigcup_{x \in \Sigma} \{x\} \times \rel^{3}
	\longrightarrow \bigcup_{x \in \Sigma} \{x\} \times \textnormal{T}_{f(x)}(f(\Sigma))
	=:\textnormal{Tan}(f)               \\
	P^{\textnormal{Tan}(f),\sphere^3}: \bigcup_{x \in \Sigma} \{x\}
	\times \textnormal{T}_{f(x)}\sphere^3
	\longrightarrow \bigcup_{x \in \Sigma} \{x\} \times \textnormal{T}_{f(x)}(f(\Sigma))
	= \textnormal{Tan}(f)
\end{eqnarray*}
denote the bundle morphisms which project
the entire tangent spaces $\rel^3$ respectively $\textnormal{T}_{f(x)}\sphere^3$ orthogonally into their subspaces $\textnormal{T}_{f(x)}(f(\Sigma))$ --
the tangent spaces of the immersion $f$ in the points $f(x)$ for every $x \in \Sigma$. Furthermore, $(A_{f}^0)_{\rel^3}$ and $(A_{f}^0)_{\sphere^3}$ denote the tracefree parts of $(A_{f})_{\rel^3}$ and $(A_{f})_{\sphere^3}$ respectively, e.g.
$$
(A_{f}^0)_{\sphere^3}(X,Y) := (A_{f})_{\sphere^3}(X,Y) - \frac{1}{2} \,g_f(X,Y)\, \vec H_{f,\sphere^3}
$$
and 
$$
\vec H_{f,\sphere^3}:=\textnormal{trace}((A_{f})_{\sphere^3}) \equiv (A_{f})_{\sphere^3}(e_i,e_i)
$$
(``Einstein's summation convention'') denotes the mean curvature vector of the immersion $f:\Sigma \longrightarrow \sphere^3$, where $\{e_i\}$ denotes a local orthonormal frame along the tangent bundle $T\Sigma$. Finally, in both settings $Q(A_f)$ respectively $Q(A^0_f)$ operate on vector 
fields $\phi$ which are sections of the normal bundle of $f$, i.e. which are normal along $f$, by assigning 
$$
Q(A_f)(\phi):= A_f(e_i,e_j) \langle A_f(e_i,e_j),\phi \rangle,
$$ 
which is by definition again a section of the normal bundle of $f$. Moreover, in equation (\ref{Moebius.flow}) we consider the ``normal Beltrami-Laplace operator'' $\triangle_{f}^{\perp}$ 
for an arbitrary $C^2$-immersion $f:\Sigma \longrightarrow \sphere^3$. As introduced in Section 1 of \cite{Simons.1968} or also in Section 1 of \cite{Weiner}, this is a differential operator of $2$nd order acting 
on those sections of the pullback-bundle $f^*(T\sphere^3)$ which are normal along $f$ within $T\sphere^3$ and again outputting such sections, i.e. sections of the ``normal 
subbundle $Nf$'' of $f^*(T\sphere^3)$. 
It is constructed by means of the composition of the unique Riemannian connection $\nabla^{\sphere^3}$ on $f^*(T\sphere^3)$ with pointwise orthogonal projection of each fiber of $f^*(T\sphere^3)$ into the corresponding fiber of its normal subbundle $Nf$.  
Alternatively, this notion can be defined via concrete local coordinate systems on $\Sigma$, as introduced e.g. in 
Definition 3.1 in Section 3 of \cite{Ruben.MIWF.II}.\\
In the author's article \cite{Jakob_Moebius_2016}, 
the author has proved short-time existence and uniqueness of the MIWF (\ref{Moebius.flow}), starting to move in smooth and umbilic-free immersions. 
Moreover, we recall here that the author has proved in the third part of Theorem 1.2 in \cite{Ruben.MIWF.II}, that 
- up to smooth reparametrizations -
the unique flow lines of the classical Willmore-flow 
\begin{equation}  \label{Willmore.flow}
	\partial_t f_t = -\frac{1}{2} \,
	\Big{(} \triangle_{f_t}^{\perp} \vec H_{f_t} + Q(A^{0}_{f_t})(\vec H_{f_t}) \Big{)}
	\equiv -\,\nabla_{L^2} \Will(f_t)
\end{equation}
in $\sphere^3$ converge smoothly and fully into immersions $F^*$, which parametrize conformally transformed Clifford-tori in $\sphere^3$, provided those flow lines start moving in 
a smooth parametrization $F_0:\Sigma \longrightarrow \sphere^3$
of a ``Hopf-torus'' with initial Willmore energy  
``$\Will(F_0) < 16\, \sqrt{\frac{\pi^3}{3}}$''.
Interestingly, this result does not seem to hold for 
flow lines $\{\PP(\,\cdot\,,0,F_0)\}$ of the 
MIWF \eqref{Moebius.flow}, meeting the same start conditions, 
since even this narrow class of flow lines of the MIWF 
might develop curvature singularities, even under the condition that their initial Willmore energies are 
smaller than the prominent threshold $8\pi$; 
see here Theorem 1.3 in \cite{Ruben.MIWF.V} for a precise criterion for ``full convergence'' of global flow lines of the MIWF which start moving in smooth Hopf-tori in $\sphere^3$ 
with Willmore energy smaller than $8\pi$, 
whose proof partially relies on the result of Theorem \ref{Center.manifold} below. \\
On the other hand, comparing the MIWF (\ref{Moebius.flow}) 
to the classical Willmore flow (\ref{Willmore.flow}), 
the MIWF enjoys the big advantage, that its flow lines can be 
conformally mapped - by means of stereographic projection -  
either from $\rel^3$ into $\sphere^3$ or from $\sphere^3$ 
into $\rel^3$, depending on the concrete objectives of 
the investigation. 
The ``$3$-sphere'' is a simply-connected compact Lie-group -
diffeomorphic to SU$(2)$ - can be interpreted as the set of quaternions of length $1$, is fibered by the Hopf-fibration 
$\pi:\sphere^3 \longrightarrow \sphere^2$ and contains 
the Clifford-torus - an embedded minimal surface in $\sphere^3$, 
which is the global minimizer of the 
Willmore functional \eqref{Willmore.functional} 
among all compact surfaces of genus $\geq 1$ 
immersed into $\sphere^3$, on 
account of \cite{Marques.Neves.2014}, Theorem A.
This particular mathematical situation actually 
plays a key-role in the proof of the first main theorem 
of this paper and also in the proof of Theorem 1.3 in \cite{Ruben.MIWF.V}, which we had already quoted above.  
On the other hand, ``$\rel^3$'' is a linear ambient
space for flow lines of the MIWF (\ref{Moebius.flow}), 
allowing the application of heavy machinery 
from both linear and non-linear ``Functional Analysis'', 
in order to investigate 
the evolution operator of the MIWF - as an operator between 
appropriately chosen Banach spaces - 
via its linearization, as already performed in 
the author's papers \cite{Jakob_Moebius_2016} and 
\cite{Ruben.MIWF.III} and also in Section 2 of this paper.  \\
We will exploit this flexibility of the 
MIWF successfully in this article, showing the first two  
``global existence and full-convergence results'' for the 
MIWF and also two ``stability results'' for fully convergent 
flow lines of the MIWF to Willmore-surfaces in $\sphere^3$
respectively $\rel^3$. The proof of Theorem
\ref{Center.manifold} relies on a 
combination of basic techniques of the author's first article \cite{Jakob_Moebius_2016} about the MIWF - which are 
based on Schauder a-priori estimates 
and on the continuity method - with a variant of the first 
part of Theorem 3.1 in \cite{Ruben.MIWF.III}, with 
particular computations due to Weiner \cite{Weiner} and with 
Escher's, Mayer's and Simonett's technique in \cite{Escher.Mayer.Simonett.1998}, \cite{Escher.Simonett.1998}, \cite{Simonett.1994}, \cite{Simonett.1995} and \cite{Simonett.2001} of ``invariant center manifolds'' for uniformly parabolic quasilinear evolution equations and their special application to the Willmore flow near round $2$-spheres in $\rel^3$. Therefore, the first main result of this article 
obviously shares many similarities with Theorem 1.2 in \cite{Simonett.2001}, treating the classical Willmore flow moving spherical immersions into $\rel^3$.    
\begin{theorem} [Full Convergence Theorem I]  \label{Center.manifold} 
	Let $\Sigma$ be a smooth compact torus, and let 
	$F^*:\Sigma \stackrel{\cong}\longrightarrow 
	M\big{(}\frac{1}{\sqrt 2}(\sphere^1 \times \sphere^1)\big{)}$
	be a smooth diffeomorphic parametrization of a compact torus in $\sphere^3$, which is conformally equivalent to the standard Clifford-torus $\CC$ via some conformal transformation $M\in \textnormal{M\"ob}(\sphere^3)$, 
	and let some $\beta \in (0,1)$ and $k \in \nat$ be fixed. Then, there is some small neighborhood $W=W(\Sigma,F^*,k)$ about $F^*$ in $h^{2+\beta}(\Sigma,\rel^4)$, such that for every $C^{\infty}$-smooth initial immersion 
	$F_1:\Sigma \longrightarrow \sphere^3$, which is contained in $W$, the unique flow line $\{\PP(t,0,F_1)\}_{t \geq 0}$ of the MIWF exists globally and converges - up to smooth reparametrization - fully to a smooth and diffeomorphic parametrization of a torus in $\sphere^3$, which is again conformally equivalent to the standard Clifford-torus $\CC$. This full convergence takes place w.r.t. the $C^k(\Sigma,\rel^4)$-norm and at an exponential rate. 
	\qed
\end{theorem}
\noindent 
In Theorem \ref{Center.manifold}, ``$h^{2+\beta}(\Sigma,\rel)$'' 
denotes the ``little H\"older space'', modelled on $\Sigma$, 
of differentiation order $2+\beta$, see here e.g. 
\cite{Escher.Mayer.Simonett.1998}, p. 1419, or 
\cite{Shao.Simonett.2014}, p. 219,
for a precise definition. 
Combining this result with Theorem \ref{Frechensbergo} below, we obtain the following ``stability theorem'' for fully convergent flow lines of the MIWF into ``a Clifford-torus''. 
\begin{theorem} [Stability Theorem I] \label{main.result.1}
	Suppose that $\Sigma$ is a smooth compact torus and that 
	$F_0:\Sigma \longrightarrow \sphere^3$ is a $C^{\infty}$-smooth and umbilic-free immersion, whose corresponding unique flow line $\{\PP(t,0,F_0)\}_{t\geq 0}$ of the MIWF exists 
	globally and converges fully and smoothly - up to smooth 
	reparametrization - to a parametrization $F^*$ of a conformal image of the Clifford-torus in $\sphere^3$.   
	Then, for any given $k\in \nat$ and $\gamma \in (0,1)$ there is an open ball $B_r(F_0)$
	about $F_0$ in $C^{4,\gamma}(\Sigma,\rel^4)$, with $r=r(\Sigma,F_0,F^*,k,\gamma)>0$, such that: \\
	For every $C^{\infty}$-smooth immersion 
	$F:\Sigma \longrightarrow \sphere^3$ being contained in 
	the open ball $B_r(F_0)$, there is a smooth 
	family of smooth diffeomorphisms 
	$\Psi^F_t:\Sigma \stackrel{\cong}\longrightarrow \Sigma$, for $t \geq 0$, a M\"obius-transformation 
	$M^F\in \textnormal{M\"ob}(\sphere^3)$ 
	and a $C^{\infty}$-diffeomorphism 
	$E^*_F:\Sigma \stackrel{\cong}\longrightarrow  
	M^F\big{(}\frac{1}{\sqrt 2}(\sphere^1 \times \sphere^1)\big{)}$, such that the reparametrized flow line $\{\PP(t,0,F)\circ \Psi^F_t\}_{t\geq 0}$ of the MIWF (\ref{Moebius.flow}) 
	converges fully and exponentially fast to $E^*_F$ in $C^k(\Sigma,\rel^4)$; thus we have:
	\begin{equation} \label{full.convergence.MIWF.2}
		\PP(t,0,F) \circ \Psi^F_t \longrightarrow E^*_F \qquad \textnormal{in}\,\,\, C^k(\Sigma,\rel^4), \quad \textnormal{as} \quad t \to \infty.
	\end{equation} 
	\qed
\end{theorem}
\noindent
It should also be noted here, that 
the proofs of Theorems \ref{Center.manifold} and \ref{main.result.1} \underline{do not require} any type of the ``Lojasiewicz-Simon-gradient-inequality'' 
for the Willmore-functional. But we shall prove the following ``full convergence theorem'' by means of the well-known trick 
of using the ``Lojasiewicz-Simon-gradient-inequality'' for a certain real-analytic functional $\FF$ 
- see e.g. \cite{Chill.Schatz.2009} or \cite{Simon.1983} -  
in order to prove simultaneously global existence of a 
flow line of the corresponding $L^2$-gradient flow and also its full convergence, combined with Rivi\`ere's \cite{Riviere.2008} and Bernard's \cite{Bernard.2016} reformulation of the $L^2$-gradient $\nabla_{L^2}\Will$ of the Willmore functional $\Will$ as a differential operator \underline{in divergence form}.\\
\begin{theorem} [Full Convergence Theorem II]  \label{Convergence.to.local.minimizer} 
	Let $\Sigma$ be a smooth compact torus, $k\in \nat$, with $k\geq 4$, and $\alpha \in (0,1)$ be given, and let $F^*:\Sigma \longrightarrow \rel^3$ be a umbilic-free and $C^{\infty}$-smooth Willmore immersion, 
	which locally minimizes the Willmore functional in the $C^k$-norm, in the sense that there exists some 
	$\delta>0$, such that for any immersion 
	$f:\Sigma \longrightarrow  \rel^3$ with 
	$\parallel f - F^* \parallel_{C^k(\Sigma)}<\delta$ 
	there holds 
	\begin{equation}  \label{local.minimum}
		\Will(f) \geq \Will(F^*).
	\end{equation}  
	Then there exists some $\varepsilon=\varepsilon(\Sigma,F^*,k,\alpha)
	\in (0,\delta)$, such that for any $C^{\infty}$-smooth immersion $f_0:\Sigma \longrightarrow \rel^3$ with 
	$\parallel f_0 - F^* \parallel_{C^{k,\alpha}(\Sigma)}<\varepsilon$
	the corresponding flow line $\{\PP(t,0,f_0)\}_{t\geq 0}$ of the MIWF exists globally and converges fully in the 
	$C^{k-1,\alpha}$-norm - up to smooth reparametrization -
	to a umbilic-free Willmore immersion $F_{\infty}$ of 
	class $C^k$, as $t\to \infty$, and this limit immersion is a $C^k$-local minimizer of the Willmore energy as well, satisfying $\Will(F_{\infty})=\Will(F^*)$.  
\qed 
\end{theorem}       
\noindent 
Hence, the MIWF can be used in order to detect  
$C^k$-local minimizers of the Willmore energy among 
$C^k$-immersions of the torus $\Sigma$ into $\rel^3$.
Combining Theorem \ref{Convergence.to.local.minimizer} 
with Theorem \ref{Frechensbergo} below, we obtain the second ``stability theorem'', namely for fully convergent flow lines of the MIWF into $C^4$-local minimizers of the Willmore functional. 
\begin{theorem}   [Stability Theorem II] \label{main.result.2}
	Suppose that $\Sigma$ is a smooth compact torus and that $F_0:\Sigma \longrightarrow \rel^3$ is a $C^{\infty}$-smooth and umbilic-free immersion, whose corresponding flow line $\{\PP(t,0,F_0)\}_{t\geq 0}$ of the MIWF exists globally and converges fully and smoothly 
	- up to smooth reparametrization - to a $C^{\infty}$-smooth parametrization $F^*$ of a umbilic-free $C^4$-local minimizer of the Willmore functional $\Will$, in the sense of formula (\ref{local.minimum}) with $k=4$. Then, for any given $\gamma \in (0,1)$ there is an open ball $B_r(F_0)$
	about $F_0$ in $C^{4,\gamma}(\Sigma,\rel^3)$, with $r=r(\Sigma,F_0,F^*,\gamma)>0$, such that: \\
	For every $C^{\infty}$-smooth immersion 
	$F:\Sigma \longrightarrow \rel^3$ being contained in the open ball $B_r(F_0)$, there is a smooth family of smooth diffeomorphisms 
	$\Psi^F_t:\Sigma \stackrel{\cong}\longrightarrow \Sigma$, 
	for $t \geq 0$, such that the reparametrized flow 
	line $\{\PP(t,0,F)\circ \Psi^F_t\}_{t\geq 0}$ 
	of the MIWF (\ref{Moebius.flow}) converges fully in the
	$C^{3,\gamma}$-norm to a umbilic-free Willmore immersion $F_{\infty}$ of class $C^4$, as $t\to \infty$, and this limit immersion is a $C^4$-local minimizer of the Willmore functional $\Will$, satisfying $\Will(F_{\infty})=\Will(F^*)$. 
	\qed
\end{theorem}
\noindent	
\begin{remark}\,
\begin{itemize} 
	\item [1)]	
	Employing the conformal invariance of the Willmore functional, of its $L^2$-gradient and of the MIWF-equation \eqref{Moebius.flow}, the statements of Theorems \ref{Convergence.to.local.minimizer} 
	and \ref{main.result.2} also hold for flow lines of 
	the MIWF in $\sphere^3$.
	\item [2)] We should also mention here, that the limit Willmore immersion $F_{\infty}$ in Theorems 
	\ref{Convergence.to.local.minimizer} and \ref{main.result.2} 
	cannot have any branch points on account of estimate 
	\eqref{enclosure.infty} in the proof of Theorem 
	\ref{Convergence.to.local.minimizer}. 
	Hence, on account of Corollary 4.4 in 
	\cite{Riviere.Park.City.2013} the classical 
	Willmore immersion $F_{\infty}$ 
	can be reparametrized by a Lipschitz homeomorphism 
	$\Psi_{\infty}:\Sigma \longrightarrow \Sigma$, such 
	that $F_{\infty}\circ \Psi_{\infty}$ is conformal
	w.r.t. an appropriate complex structure $c$ on $\Sigma$, 
	and moreover $F_{\infty}\circ \Psi_{\infty}$ has to be 
	$C^{\infty}$-smooth on $\Sigma$ w.r.t. the complex 
	structure $c$ - defining a smooth structure 
	on $\Sigma$ in particular - on account of Theorem 7.11 in 
	\cite{Riviere.Park.City.2013}. Corollary 4.4 and Theorem 
	7.11 in \cite{Riviere.Park.City.2013} originate from 
	H\'elein's and Rivi\`ere's seminal work \cite{Helein.2002},   \cite{Riviere.2008}; see Theorem 5.1.1 in \cite{Helein.2002} 
	and Theorems I.2 and I.5 in \cite{Riviere.2008}.
	\qed 
\end{itemize}  
\end{remark} 

\section{Preparations for the proofs of Theorems 
	\ref{Center.manifold} -- \ref{main.result.2}} \label{Prep.proofs}

\begin{definition}  \label{Umbilic.free}
	Let $\Sigma$ be a smooth compact torus, and let $M=\rel^3$ or $M=\sphere^3$. We denote by 
	$C^{\infty}_{\textnormal{Imm}}(\Sigma,M)$
	the set of $C^{\infty}$-smooth immersions 
	$F:\Sigma \longrightarrow M$ of the torus $\Sigma$ into $M$.
	\qed
\end{definition}
\noindent
In order to transfer the method in \cite{Escher.Mayer.Simonett.1998},
\cite{Escher.Simonett.1998}, \cite{Simonett.1995} and \cite{Simonett.2001} to the MIWF (\ref{Moebius.flow}) for families of immersions $f_t:\Sigma \longrightarrow \sphere^3$, we have to establish ``Fermi coordinates'' in a sufficiently small open neighborhood $U$ of the Clifford torus 
$\CC = \frac{1}{\sqrt 2}(\sphere^1 \times \sphere^1)$ in $\sphere^3$. At first, we recall from \cite{Lang}, p. 108, that 
the tangent bundle of $\sphere^3$ splits along $\CC$ into 
a direct sum of vector bundles 
$$ 
T\sphere^3\lfloor_{\CC} = 
T\CC \oplus N\CC 
$$
namely into the tangent bundle and into the normal bundle of $\CC$ within $T\sphere^3$. Now, we can construct Fermi coordinates in a canonical way by means of the restriction of the exponential map 
$\exp \equiv \exp^{\sphere^3} : \dom(\exp^{\sphere^3}) \subset T\sphere^3 \longrightarrow \sphere^3$ to the normal bundle $N\CC$ - a smooth subbundle of $T\sphere^3$:
$$ 
\exp\lfloor_{N\CC} : \textnormal{dom}(\exp) \cap N\CC \longrightarrow \sphere^3,
$$ 
because the proof of Theorem 5.1 in \cite{Lang} guarantees, that there is a small open neighborhood $Z$ of the zero-section in the total space of $N\CC$ and an open neighborhood $U$ of the torus $\CC$ in $\sphere^3$, such that 
\begin{equation}  \label{tubular.map}
	\exp\lfloor_{N\CC}: Z  \stackrel{\cong}\longrightarrow  U
\end{equation}  
is a smooth diffeomorphism. Hence, the restriction $\exp\lfloor_{N\CC}$ to a sufficiently small open neighborhood $Z$ of the zero-section of the bundle $N\CC \longrightarrow \CC$ is a ``tubular map'', and the corresponding open neighborhood $U$ of $\CC$ thus turns 
out to be a ``tube'' about $\CC$, 
in the language of differential topology, see e.g. the treatise \cite{Lang}, p. 108. In other words, having chosen a global unit normal field $\nu_{\CC}$ along the orientable surface $\CC$, i.e. a smooth section of the bundle $N\CC$ with constant length $1$, any point $p \in U$ can be written in the form 
``$p = \exp_x(r\,\nu_{\CC}(x))$'' for a unique point $x = x(p) \in \CC$ and a unique number $r=r(p) \in \rel$, 
which are thus globally defined ``Fermi coordinates'' 
$(x,r)\in \CC \times \rel$ for all points $p \in U$.
Hence, having chosen a unit normal field $\nu_{\CC}$ along $\CC$, statement (\ref{tubular.map}) yields a smooth diffeomorphism
\begin{equation}  \label{tubular.map.2}
	X : \CC \times (-a,a) \stackrel{\cong}\longrightarrow \textnormal{image}(X)=:U_a, \qquad  X(x,r):=\exp_x(r\,\nu_{\CC}(x)),
\end{equation}  
onto an open neighborhood $U_a$ of the 
torus $\CC$ in $\sphere^3$, provided $a>0$ is chosen sufficiently small, more precisely smaller than the width of the tube $Z$ about 
the zero-section in the total space of $N\CC$.
See here also Section 4.1 in \cite{Pruess.Simonett} for the construction of the map $X$ in Euclidean space. 
Taking the inverse of the smooth diffeomorphism $X$ in (\ref{tubular.map.2}), we obtain a well-defined and unique 
pair of smooth coordinate functions 
\begin{equation}  \label{coordinate.functions} 
	S: U_a \longrightarrow \CC  \qquad  \textnormal{and} \qquad \Lambda: U_a  \longrightarrow  (-a,a).         
\end{equation} 
Now, suppose there is some smooth manifold 
$\Gamma$ in $U_a$, which has the property that the coordinate 
function $S$ maps $\Gamma$ bijectively onto $\CC$. We thus 
obtain a unique smooth function 
\begin{equation}  \label{rho.Gamma}
	\rho \equiv \rho_{\Gamma} :\CC \longrightarrow (-a,a)   \quad \textnormal{by setting} \quad 
	\rho(x):=\Lambda \circ (S \lfloor_{\Gamma})^{-1}(x), \,\, x \in \CC, 
\end{equation}   
where $\Lambda$ and $S$ are the smooth coordinate functions from line (\ref{coordinate.functions}). Obviously, the function $\rho_{\Gamma}$ measures the pointwise ``geodesic distance'' of 
the set $\Gamma$ from any fixed point $x\in \CC$, and we therefore recover $\Gamma$ as a graph over $\CC$:
$$ 
\Gamma = \textnormal{image}\big{(}[\CC \ni x\mapsto X(x,\rho(x))] \big{)}
$$
just by construction of the diffeomorphism $X$ in (\ref{tubular.map.2}), of $S$ and $\Lambda$ in (\ref{coordinate.functions}) and of $\rho$ in (\ref{rho.Gamma}). Conversely, suppose that there is some function $\rho:\CC \longrightarrow (-a,a)$ of class $h^{2+\alpha}(\CC)$, for some arbitrarily fixed $\alpha \in (0,\beta)$ - where $\beta \in (0,1)$ is already given by the asserted statement of Theorem \ref{Center.manifold} - then the set 
\begin{equation}  \label{Gamma}
	\Gamma(\rho):= \textnormal{image}\big{(}[\CC \ni x\mapsto X(x,\rho(x))] \big{)} \subset U_a 
\end{equation} 
is a $2$-dimensional manifold of class $h^{2+\alpha}$, and $\Gamma(\rho)$ is the level set $\{p \in U_a \,|\,\Phi_{\rho}(p)=0 \,\}$ of the function   
$$
\Phi_{\rho}:U_a \longrightarrow \rel \quad \textnormal{defined by} \quad  \Phi_{\rho}(p):= \Lambda(p) - \rho(S(p)).
$$ 
Obviously, the function $\Phi_{\rho}$ is just as smooth as 
the function $\rho$ is, thus here it is of class $h^{2+\alpha}(U_a)$.
Now suppose, that we have a time-dependent function 
$\rho:\CC \times [0,T) \longrightarrow (-a,a)$ of class 
$C^0([0,T);h^{2+\alpha}(\CC))$ for the above fixed 
$\alpha \in (0,\beta)$. We thus consider the time-dependent function 
$$
\Phi_{\rho_t}(p,t):=\Lambda(p) - \rho(S(p),t),
$$
and we obtain closed and compact $h^{2+\alpha}$-manifolds $\Gamma(\rho_t)$, for $t \in [0,T)$, in the neighborhood $U_a$ of $\CC$ in $\sphere^3$ as level sets: 
\begin{equation}  \label{Gamma.t}
	\Gamma(\rho_t) = \{\,p \in U_a \,|\, \Phi_{\rho_t}(p,t)=0 \,\}. 
\end{equation}   
Now, in combination with equation (\ref{Gamma}) we have the equation 
$$ 
\Phi_{\rho_t}(X(x,\rho(x,t)),t)=0   \qquad \forall (x,t) \in \CC \times [0,T).
$$
The chain rule thus yields: 
\begin{eqnarray*} 
	0= \Big{\langle} \nabla^{\sphere^3}_p \Phi_{\rho_t}(X(x,\rho(x,t)),t),
	\frac{\partial}{\partial t}(X(x,\rho(x,t))) \Big{\rangle}  + 
	\frac{\partial}{\partial t}\Phi_{\rho_t}(X(x,\rho(x,t)),t)                  \\
	=  \pm \Big{|}\nabla^{\sphere^3}_p \Phi_{\rho_t}(X(x,\rho(x,t)),t)\Big{|} \,
	\Big{|}(\partial_t)^{\perp_{X}}(X(x,\rho(x,t)))\Big{|}
	- \frac{\partial \rho}{\partial t}(x,t)
\end{eqnarray*} 
for $(x,t) \in \CC \times [0,T)$, where 
``$\frac{\partial}{\partial t}\Phi_{\rho_t}(X(x,\rho(x,t)),t)$'' 
means 
$\frac{\partial}{\partial t}\Phi_{\rho_t}(p,t)\lfloor_{p=X(x,\rho(x,t))}$.
Here, ``$(\partial_t)^{\perp_{X}}(X(x,\rho(x,t)))$'' denotes
the normal component $V(x,t)$ of the velocity field 
$\partial_t (X(\,\cdot\,,\rho(\,\cdot\,,t)))$ of the family 
$\{\Gamma(\rho_t)\}_{t \in [0,T)}$ of parametrized 
moving manifolds from line (\ref{Gamma.t}), evaluated in their 
points $X(x,\rho(x,t))$. Hence, as on p. 1423 in \cite{Escher.Mayer.Simonett.1998} we obtain the crucial equation   
\begin{equation} \label{normal.speed.1}
	\pm | V(x,t) | = \frac{\frac{\partial \rho}{\partial t}(x,t)}
	{|\nabla^{\sphere^3}_p \Phi_{\rho_t}(X(x,\rho(x,t)),t)|} 
\end{equation} 
for $(x,t) \in \CC \times [0,T)$. Now, if the distance function $\{\rho_t\}$ is additionally of class $C^{\infty}((0,T);C^{\infty}(\CC))$, then the family of functions 
\begin{equation}  \label{theta.exp.rho}
	\theta_{\rho}(x,t):= \exp_x(\rho(x,t)\,\nu_{\CC}(x))	\qquad \textnormal{for} 
	\,\, x\in \CC \quad \textnormal{and} \quad t\in [0,T),  
\end{equation} 
is $C^{\infty}$-smooth on $\CC \times (0,T)$ and parametrizes the family $\{\Gamma(\rho_t)\}_{t \in [0,T)}$ of $C^{\infty}$-smooth surfaces by means of $C^{\infty}$-diffeomorphisms 
$\theta_{\rho}(\,\cdot\,,t):\CC \stackrel{\cong}\longrightarrow \Gamma(\rho_t)$. Moreover, we note that there holds 
\begin{equation}  \label{laplacian}
	\frac{1}{2}\, \triangle_{\Gamma(\rho_t)}^{\perp} \vec H_{\Gamma(\rho_t)} 
	= \triangle_{\Gamma(\rho_t)}^{\perp} (H_{\Gamma(\rho_t)}\, \nu_{\Gamma(\rho_t)})
	= \triangle_{\Gamma(\rho_t)}(H_{\Gamma(\rho_t)}) \, \nu_{\Gamma(\rho_t)}
\end{equation} 
for the ``Beltrami-laplacian'' $\triangle_{\Gamma(\rho_t)}^{\perp}$ in the normal bundle $N\Gamma(\rho_t)$ of the submanifold $\Gamma(\rho_t)$
of $\sphere^3$, by Sections 1 and 2 in \cite{Simons.1968} or also by formula (\ref{laplacian.2}) below, and that  
\begin{equation}  \label{little.computation}
	\frac{1}{2}\, |A^0_{\Gamma(\rho_t)}|^2 \, \vec H_{\Gamma(\rho_t)} =   2 \,(|H_{\Gamma(\rho_t)}|^2 - K_{\Gamma(\rho_t)})\, H_{\Gamma(\rho_t)} \, \nu_{\Gamma(\rho_t)}
\end{equation}
holds on $\Gamma(\rho_t)$ for every $t\in [0,T)$, 
where the symbols ``$2 \,H_{\Gamma(\rho_t)}$'' and 
``$K_{\Gamma(\rho_t)}$'' denote the trace and the 
determinant respectively of the scalar second fundamental form 
$(A_{\Gamma(\rho_t)})_{\sphere^3}$ 
of the submanifold $\Gamma(\rho_t)$ of $\sphere^3$; 
see here \cite{Weiner}, p. 22. 
Hence, as in formulae (1.1)--(2.2) in \cite{Escher.Mayer.Simonett.1998}
we infer from formulae (\ref{Moebius.flow}), (\ref{normal.speed.1}), (\ref{laplacian}) and (\ref{little.computation}), that a family of 
immersions $f_t:\Sigma \longrightarrow \sphere^3$, $t\in [0,T)$,
parametrizing surfaces $\Gamma(\rho_t)$
in the open neighborhood $U_a$ of $\CC$ in $\sphere^3$, which are implicitly given by equation (\ref{Gamma.t}) in terms of a time-dependent distance function $\rho:\CC \times [0,T) \longrightarrow (-a,a)$ of regularity class 
$$ 
\rho \in C^0([0,T);h^{2+\beta}(\CC)) \cap C^{\infty}((0,T);C^{\infty}(\CC)),
$$ 
moves according to the ``relaxed version'' 
\begin{equation}  \label{Moebius.flow.2}
	(\partial_t)^{\perp_{f_t}}(f_t) 
	\stackrel{!}= -\frac{1}{2} \frac{1}{|A^0_{f_t}|^4} \,
	\Big{(} \triangle_{f_t}^{\perp} \vec H_{f_t} + Q(A^{0}_{f_t})(\vec H_{f_t}) \Big{)} \equiv -\frac{1}{|A^0_{f_t}|^4} \,\nabla_{L^2} \Will(f_t)
\end{equation}
of the MIWF (\ref{Moebius.flow}) on $\Sigma \times (0,T)$, 
if and only if the prescribed distance function $\rho=\{\rho_t\}$ satisfies the evolution equation
\begin{eqnarray} \label{normal.speed.3}
	\frac{\partial \rho}{\partial t}(x,t) = 
	\frac{|\nabla^{\sphere^3}_p \Phi_{\rho_t}(\theta_{\rho}(x,t),t)|}
	{|A^0_{\Gamma(\rho_t)}(\theta_{\rho}(x,t))|^4} \,
	\Big{(} 
	\theta_{\rho_t}^*\big{(}\triangle_{\Gamma(\rho_t)}H_{\Gamma(\rho_t)}\big{)}(x,t) +  \nonumber   \\
	+ 2 \,\big{(}|H_{\Gamma(\rho_t)}(\theta_{\rho}(x,t))|^2 - K_{\Gamma(\rho_t)}(\theta_{\rho}(x,t))\big{)} \,H_{\Gamma(\rho_t)}(\theta_{\rho}(x,t)) \, 
	\Big{)}                                       \\  
	=: G(\rho_t)(x)    \qquad     \nonumber   
\end{eqnarray}
for $(x,t) \in \CC \times (0,T)$, whose initial 
value $\rho_0$ is determined by the initial 
$h^{2+\beta}$-manifold 
$\Gamma_0 \subset U_a$ on account of formulae 
(\ref{tubular.map.2}) -- (\ref{Gamma.t}). In equation (\ref{normal.speed.3}) the symbols 
``$2\, H_{\Gamma(\rho_t)}(\theta_{\rho}(x,t))$'' and 
``$K_{\Gamma(\rho_t)}(\theta_{\rho}(x,t))$'' denote 
the trace and the determinant respectively 
of the scalar second fundamental form 
$(A_{\Gamma(\rho_t)})_{\sphere^3}$ 
of the submanifold $\Gamma(\rho_t)$ of $\sphere^3$ 
evaluated in its point $\theta_{\rho}(x,t)$,  
and the symbol ``$\triangle_{\Gamma(\rho_t)}$'' denotes the Beltrami-Laplace operator of the submanifold $\Gamma(\rho_t)$ of $\sphere^3$, operating on functions on $\Gamma(\rho_t)$, which can be obtained from the action of the normal laplacian $\triangle_{\Gamma(\rho_t)}^{\perp}$ in $N\Gamma(\rho_t)$ by means of formula (\ref{laplacian.2}) below. Now, we choose 
some $\beta_0 \in (\alpha,\beta)$, where $\alpha<\beta \in (0,1)$ has been arbitrarily fixed above line (\ref{Gamma}), 
and we define the open subset   
$$
\UU^a_{\beta_0} := 
\{\rho \in h^{2+\beta_0}(\CC)\, |\, 
\parallel \rho \parallel_{L^{\infty}(\CC)} <a \,\} 
$$
of the Banach space $h^{2+\beta_0}(\CC)$, for some sufficiently small $a>0$ as in (\ref{tubular.map.2}). As in Lemma 2.1 in \cite{Escher.Mayer.Simonett.1998} we can prove here:
\begin{lemma}   \label{Lemma.2.1} 
	The differential operator 
	\begin{equation}   \label{G}
		-G(\rho) \equiv - \frac{|\nabla^{\sphere^3}_p \Phi_{\rho}\circ \theta_{\rho}|}{|A^0_{\rho}|^4} \,
		\Big{(} \triangle_{\rho}H_{\rho} +  2 \,H_{\rho} \,(H_{\rho}^2-K_{\rho})\Big{)}   
	\end{equation} 
	from line (\ref{normal.speed.3}), having abbreviated here 
	$A^0_{\rho}:= A^0_{\Gamma(\rho)}\circ \theta_{\rho}$, 
	$H_{\rho}:= H_{\Gamma(\rho)}\circ \theta_{\rho}$ and  
	$K_{\rho}:= K_{\Gamma(\rho)}\circ \theta_{\rho}$, is a  
	uniformly elliptic quasilinear operator. More precisely, $G$ can be decomposed in the following way:
	\begin{equation}  \label{split.G}
		G(\rho) = - P(\rho).\rho + F(\rho)   
	\end{equation} 
	for every $\rho \in \VV^{a}_{\alpha}:=h^{4+\alpha}(\CC) \cap \UU^{a}_{\beta_0}$, where 
	$$
	P :\UU^{a}_{\beta_0} \longrightarrow \Lift(h^{4+\alpha}(\CC),h^{\alpha}(\CC))
	$$
	is a uniformly elliptic quasilinear operator of class 
	$C^{\infty}(\UU^{a}_{\beta_0},\Lift(h^{4+\alpha}
	(\CC),h^{\alpha}(\CC)))$, 
	and $F \in C^{\infty}(\UU^{a}_{\beta_0},h^{\beta_0}(\CC))$ is a non-linear operator of only second order, satisfying $F(0)=0$ on $\CC$. 
	In particular, $-P(\rho)$ generates a strongly continuous analytic semigroup on $h^{\alpha}(\CC)$, i.e. 
	$P(\rho)\in \HH(h^{4+\alpha}(\CC),h^{\alpha}(\CC))$, for every $\rho \in \UU^{a}_{\beta_0}$.  
\end{lemma}	
\proof: We shall closely follow the lines of the proofs of Lemma 2.1 in \cite{Escher.Mayer.Simonett.1998} 
and Lemma 3.1 in \cite{Escher.Simonett.1998} during the proof of this lemma. As in the proof of Lemma 2.1 in \cite{Escher.Mayer.Simonett.1998} we firstly choose an atlas $\{\,\OO_l\, | \,l=1,\ldots,m\,\}$ of open coordinate 
neighborhoods $\OO_l$ on $\CC$, yielding partial derivatives $\partial_j$, $j=1,2$, on $\OO_l$, and we pull back the Euclidean metric $g_{\sphere^3}$ via the diffeomorphism 
$$
X_l:= X\lfloor_{\OO_l \times (-a,a)} : \OO_l \times (-a,a) \stackrel{\cong}\longrightarrow  
\textnormal{image}(X\lfloor_{\OO_l \times (-a,a)}) =: R_l(a) \subset \sphere^3
$$ 
which is a restriction of the diffeomorphism $X$ in line 
(\ref{tubular.map.2}) to $\OO_l \times (-a,a)$. 
We thus obtain a smooth product metric 
\begin{equation}  \label{g.l}
	g_l:= X_l^*(g_{\sphere^3}\lfloor_{R_l(a)})
	= w_l(r) + dr \otimes dr \quad \textnormal{on} \,\,\, T(\OO_l \times (-a,a)),
\end{equation} 
where $w_l(r)$ is the metric on $T(\OO_l \times \{r\})\cong T\OO_l$, whose coefficients can be explicitly given by 
$$ 
(w_l(r))_{jk}(x)=g_{\sphere^3}(\partial_jX_l(x,r), \partial_kX_l(x,r)),
\quad \textnormal{for} \quad (x,r) \in \OO_l \times (-a,a)
$$
and for $l=1,\ldots,m$, where $\partial_j,\partial_k$ are the partial derivatives on $\OO_l$ introduced above. 
Moreover, for $\rho \in \UU^{a}_{\beta_0}$ we use the notion of $w_l(r)$, in order to define the metric $w(\rho)$ on $T\CC$ by setting 
$$ 
(w(\rho))_{jk}(x)=(w_l(\rho(x)))_{jk}(x) \equiv 
g_{\sphere^3}(\partial_jX(x,\rho(x)), \partial_kX(x,\rho(x))),
\quad \textnormal{for} \,\, x \in \OO_l,
$$
and for $l=1,\ldots,m$. See here also \cite{Escher.Mayer.Simonett.1998}, p. 1424, and \cite{Escher.Simonett.1998}, p. 275. In other words, 
$w(\rho)$ is the unique metric $\eta$ on $T\CC$ 
such that  
$$
\eta(x) + dr \otimes dr = (g_l)_{(x,\rho(x))} 
\quad \textnormal{on} \,\,\, T_{(x,\rho(x))}(\OO_l \times (-a,a)) \quad \forall \, x \in \OO_l, 
$$	
and for each $l=1,\ldots,m$. From $w(\rho)$ we also obtain the metric $w^*(\rho)$ on the cotangent bundle $T^*(\OO_l)$.
On the other hand, we obtain another metric $\sigma(\rho)$ 
on $T\CC$ by means of pulling back the 
Euclidean metric $g_{\sphere^3}\lfloor_{\Gamma(\rho)}$ with the $h^{2+\beta_0}$-diffeomorphism 
$\theta_{\rho}:\CC \stackrel{\cong}\longrightarrow  \Gamma(\rho)$, given by $\theta_{\rho}(x):= \exp_x(\rho(x)\,\nu_{\CC}(x))$ for $x\in \CC$.
Now, as in the proof of Lemma 2.1 in \cite{Escher.Mayer.Simonett.1998} 
we can conclude from $\Phi_{\rho}(p) = \Lambda(p) - \rho(S(p))$ and from the definition of the metric $w(\rho)$, that there holds 
\begin{equation}   \label{L.rho}  
	L_{\rho}^2(x):=|\nabla^{\sphere^3}_p \Phi_{\rho}(X_l(x,\rho(x)))|^2
	= 1 + w^*(\rho)_x(d\rho(x),d\rho(x))  \quad \forall \, x \in \OO_l,     
\end{equation}
and for each $l=1,\ldots,m$. Moreover, as in Lemma 3.1 in 
\cite{Escher.Simonett.1998} we pull back the scalar mean curvature 
$$
H_{\rho}(x):= H_{\Gamma(\rho)}(\theta_{\rho}(x)) \equiv 
(\theta_{\rho}^*H_{\Gamma(\rho)})(x)      \qquad \forall\, x\in\CC
$$  
and we write it in terms of $L_{\rho}$, the metric $w(\rho)$, and in terms of the Christoffel-symbols $\Gamma^i_{jk}(\rho)$ of the metric 
$(g_l)_{(s,\rho(s))}$: 
\begin{equation}   \label{H.rho} 
	H_{\rho} = P_1(\rho).\rho + F_1(\rho) \quad \textnormal{on} \,\,\,\CC
\end{equation}
for any $\rho \in \UU^{a}_{\beta_0}$, where we have  
\begin{eqnarray}   \label{P.1} 
	P_1(\rho) = \frac{1}{L_{\rho}^3} \, 
	\Big{(} \big{(} - L_{\rho}^2 \,w^{jk}(\rho) + 
	w^{jl}(\rho) \, w^{km}(\rho) \,\partial_l\rho \partial_m\rho \big{)} \, 
	\partial_j\, \partial_k       +                \qquad                  \nonumber    \\
	+ \big{(} L_{\rho}^2\, w^{jk}(\rho) \, \Gamma^i_{jk}(\rho) + 
	w^{jl}(\rho) \, w^{ki}(\rho) \,\Gamma^3_{jk}(\rho) \, \partial_l\rho        
	+ 2 \, w^{km}(\rho) \Gamma^i_{3k}(\rho) \, \partial_m\rho   -        \quad   \\
	- w^{jl}(\rho) \, w^{km}(\rho)\, \Gamma^i_{jk}(\rho) \, 
	\partial_l\rho \,\partial_m\rho \big{)} \,\partial_i\Big{)},\qquad \nonumber                    
\end{eqnarray}
and 
\begin{equation}  \label{F.1}
	F_1(\rho) = -\frac{1}{L_{\rho}} \, w^{jk}(\rho) \,\Gamma^3_{jk}(\rho),
\end{equation}
where summation only runs from $1$ to $2$ for repeated indices and where we have explicitly 
\begin{equation}  \label{Christoffel}
	\Gamma^3_{jk}(\rho)(x)= 
	g_{\sphere^3}(\partial_j\partial_k X_l, \partial_3X_l)(x,\rho(x)) \qquad
	\textnormal{for} \quad  x\in \OO_l 
\end{equation} 
and for each $l=1,\ldots,m$. As in the proof of Lemma 2.1 in \cite{Escher.Mayer.Simonett.1998}, one can derive from formulae (\ref{L.rho}) and (\ref{P.1}) and Cauchy-Schwarz' inequality, that the symbol $p^{\pi}_1(\rho)$ of the leading 2nd order term of the operator $-P_1(\rho)$ in (\ref{P.1}) satisfies:
\begin{equation}  \label{estimate.p.pi.1} 
	p^{\pi}_1(\rho)(\xi) \geq  \frac{1}{L_{\rho}^3} \,w^*(\rho)(\xi,\xi) \qquad  \forall \,\xi\in T^*(\CC).
\end{equation} 
Furthermore, we recall that the Beltrami-Laplace operator satisfies the transformation law
$$
(\theta_{\rho}^*\triangle_{\Gamma_{\rho}})(\,\cdot\,) 
=\triangle_{\rho}(\theta_{\rho}^*(\,\cdot\,)),  
$$ 
where 
\begin{equation} \label{triangle.rho}
	\triangle_{\rho}(f) := \sigma^{jk}(\rho) \,
	\Big{(} \partial_{jk}(f) - \gamma^i_{jk}(\rho) \,\partial_i(f) \Big{)}
	\qquad \textnormal{for} \,\,\, f\in C^{\infty}(\CC)
\end{equation} 
denotes the Beltrami-Laplace operator on $\CC$ w.r.t. the pullback-metric $\sigma(\rho)= \theta_{\rho}^*(g_{\sphere^3}\lfloor_{\Gamma(\rho)})$ and where $\gamma^i_{jk}(\rho)$, $i,j,k=1,2$, denote the Christoffel-symbols of the metric $\sigma(\rho)$ on $T\CC$. Combining this with formulae (\ref{G}), (\ref{H.rho}) and (\ref{P.1}) one can start proving equation (\ref{split.G}) for any fixed $\rho \in \VV^{a}_{\alpha} \equiv h^{4+\alpha}(\CC)\cap \UU^{a}_{\beta_0}$
as in the proof of Lemma 2.1 in \cite{Escher.Mayer.Simonett.1998}. 
First of all, one can compute by means of formulae (\ref{G}), (\ref{H.rho}) and (\ref{P.1}), that the leading $4$th order term of the operator $P(\rho)$ on the right hand side of equation (\ref{split.G}) is explicitly given by: 
\begin{equation}  \label{P.pi.rho}
	P^{\pi}(\rho) := \frac{1}{L_{\rho}^2\,
		|A^0_{\rho}|^4} \, 
	\sigma^{rs}(\rho)\Big{(} L_{\rho}^2\, w^{jk}(\rho) - w^{jl}(\rho) \,w^{km}(\rho) \, \partial_l \rho \partial_m\rho \Big{)} 
	\partial_r \partial_s \partial_j \partial_k,
\end{equation}
for any fixed $\rho \in \UU^{a}_{\beta_0}$. More precisely, we can write:
\begin{equation}  \label{split.the.P}
	-\frac{1}{|A^0_{\rho}|^4} \, L_{\rho} \, \triangle_{\rho}(P_1(\rho).\rho) 
	= P^{\pi}(\rho).\rho \,+  Q(\rho).\rho \quad \textnormal{for} \,\, 
	\rho \in \UU^{a}_{\beta_0} \cap h^{4+\alpha}(\CC), 
\end{equation}
where $Q(\rho)\in \Lift(h^{3+\alpha}(\CC),h^{\alpha}(\CC))$ is a quasilinear differential operator of third order, which acts on third order partial derivatives linearly, 
for any fixed $\rho \in \UU^{a}_{\beta_0} \cap h^{4+\alpha}(\CC)$. 
Now, comparing the explicit formula (\ref{P.pi.rho}) with formula (\ref{P.1}), and using estimate (\ref{estimate.p.pi.1}) for the symbol $p^{\pi}_1(\rho)$ of the leading $2$nd order term of the operator $-P_1(\rho)$ in (\ref{P.1}),
the symbol $p^{\pi}(\rho)$ of the operator $P^{\pi}(\rho)$ from line (\ref{P.pi.rho}) turns out to satisfy: 
$$ 
p^{\pi}(\rho)(\xi) \geq 
\frac{1}{L_{\rho}^2 \,|A^0_{\rho}|^4} \, 
\sigma^*(\rho)(\xi,\xi) \,w^*(\rho)(\xi,\xi)   \qquad \forall \,\xi \in T^*(\CC),
$$
proving that the operator $P^{\pi}(\rho)$ is uniformly elliptic of fourth order. Moreover, as in the proof of Lemma 2.1 in \cite{Escher.Mayer.Simonett.1998} we 
write the operator $P(\rho)$ on the right hand side of equation (\ref{split.G}) as a sum of the principal quasilinear operator $P^{\pi}(\rho)$ of fourth order and two further quasilinear operators $Q(\rho)$ and $R(\rho)$ of third order, which contain all partial derivatives of third order of the operator $G$ in formula (\ref{G}):  
\begin{equation} \label{P}
	P(\rho):= P^{\pi}(\rho) + Q(\rho) + R(\rho)   \qquad  \textnormal{for} \,\, 
	\rho \in \UU^{a}_{\beta_0} \cap h^{4+\alpha}(\CC),
\end{equation}   
where the quasilinear operator $R(\rho)$ is concretely given by 
$$
R(\rho).\rho:= - \frac{L_{\rho}}{|A^0_{\rho}|^4} \, \Big{(} \triangle_{\rho}\Big{(}\frac{1}{L_{\rho}}\Big{)}\Big{)} \, 
L_{\rho}\, F_1(\rho)\qquad  \textnormal{for} \,\, \rho \in \UU^{a}_{\beta_0} \cap h^{3+\alpha}(\CC).
$$ 
Hence, combining formulae (\ref{G}), (\ref{H.rho}), (\ref{split.the.P}) and (\ref{P}), we see as in the proof of Lemma 2.1 in \cite{Escher.Mayer.Simonett.1998}, that the remaining term in formula (\ref{split.G}) has to be the non-linear operator 
\begin{equation} \label{F.rho}
	F(\rho):= \frac{L_{\rho}}
	{|A^0_{\rho}|^4} \,\Big{(} \triangle_{\rho} F_1(\rho) 
	+ 2\,H_{\rho} \, \big{(}|H_{\rho}|^2-K_{\rho}\big{)} \Big{)}\,   
	+ R(\rho).\rho  \quad  \textnormal{for} \,\, \rho \in 
	\UU^{a}_{\beta_0} \cap h^{3+\alpha}(\CC).
\end{equation}
Moreover, it follows as in the proof of Lemma 2.1 in \cite{Escher.Mayer.Simonett.1998}
and as in the proof of Lemma 2.1 in \cite{Simonett.2001}, that the non-linear operator $F$ is of second order only and smooth, 
more precisely $F$ is of second order and it is of class $C^{\infty}(\UU^{a}_{\beta_0},\Lift(h^{4+\alpha}
(\CC),h^{\alpha}(\CC))$. 
This has completed the proof of formula (\ref{split.G}). 
Furthermore, we can verify that the uniform ellipticity of $P(\rho)$ implies that $-P(\rho)$ is sectorial in $h^{\alpha}(\CC)$ by Theorem 3.3 in \cite{Shao.Simonett.2014}, for any fixed $\rho \in \UU^{a}_{\beta_0}$. 
Since $h^{4+\alpha}(\CC)$ embeds densely into 
$h^{\alpha}(\CC)$, $-P(\rho)$ therefore generates a strongly continuous analytic semigroup in $h^{\alpha}(\CC)$, in classical notation
``$P(\rho)\in \HH(h^{4+\alpha}(\CC), h^{\alpha}(\CC))$'', for every $\rho \in \UU^{a}_{\beta_0}$. Finally, we follow 
the proof of Lemma 3.1 in \cite{Escher.Simonett.1998}, p. 276, 
and compute here by means of formulae
(\ref{tubular.map.2}) and (\ref{Christoffel}): 
\begin{eqnarray}   \label{compute.Chistoffel}
	\Gamma^3_{jk}(\rho)(x)\lfloor_{\rho=0} 
	= g_{\sphere^3}(\partial_j\partial_k X_l, \partial_3X_l)(x,0) 
	= g_{\sphere^3}(\partial_j\partial_k X_l(x,0), \nu_{\CC}(x))    \\
	= g_{\sphere^3}(\nabla^{\sphere^3}_{\partial_j}(\partial_k)(x), \nu_{\CC}(x)) 
	= g_{\sphere^3}(T_{\CC}(\partial_j,\partial_k)(x), \nu_{\CC}(x)) \quad
	\textnormal{for} \quad  x\in \OO_l,    \nonumber
\end{eqnarray} 
and for each $l=1,\ldots,m$, where $\nabla^{\sphere^3}$ and $T_{\CC}$ denote the Riemannian connection on $\sphere^3$ and the second fundamental form tensor of the injection $X(\,\cdot\,,0):\CC \hookrightarrow \sphere^3$,
and where we have used Proposition 2.2.2 on p. 68 in \cite{Simons.1968}. 
Since $\CC$ is a minimal surface in $\sphere^3$, 
i.e. since we have $H_{\rho}\lfloor_{\rho=0} \equiv 0$ on $\CC$, we can conclude from formulae (\ref{F.1}) and (\ref{compute.Chistoffel}) together with the fact that $L_{\rho}\lfloor_{\rho=0} \equiv 1$ on $\CC$ 
- on account of formula (\ref{L.rho}) - that 
\begin{eqnarray*}
	F_1(\rho)\lfloor_{\rho=0} = - \, w^{jk}(0) \,\Gamma^3_{jk}(0) = 
	- w^{jk}(0) \, g_{\sphere^3}(T_{\CC}(\partial_j,\partial_k), \nu_{\CC})   =   \nonumber \\
	= - g_{\sphere^3}(\vec H_{\CC},\nu_{\CC}) \equiv 0 \qquad \textnormal{on} \,\,\, \OO_l,
\end{eqnarray*}
for each $l=1,\ldots,m$, where we have used the definition of 
$\vec H_{\CC}$ on p. 68 in \cite{Simons.1968}. Hence, we can especially infer that $\triangle_{\rho} F_1(\rho)\lfloor_{\rho=0} \equiv 
\triangle_{\CC}(F_1(0)) \equiv 0$ on $\CC$, which finally implies that $F(\rho)\lfloor_{\rho=0} \equiv 0$ on $\CC$ by formula (\ref{F.rho}), just as asserted. 
\qed  \\\\	
\noindent
Relying on the proof of Theorem 2.2 in \cite{Escher.Mayer.Simonett.1998}, 
we infer the following fundamental existence, uniqueness and regularity result for the quasilinear parabolic equation (\ref{normal.speed.3}) from Lemma \ref{Lemma.2.1} and Section 12 
in \cite{Amann.1993}.
See also Theorem 3.1 in \cite{Simonett.1995}. We recall here, that we have chosen $\beta_0 \in (\alpha,\beta)$ above the statement of Lemma \ref{Lemma.2.1}, implying that $h^{2+\beta}(\CC)$ embeds compactly into $h^{2+\beta_0}(\CC)$.  
\begin{theorem}   \label{Short.time.existence}
	For any $\rho_0\in \UU^{a}_{\beta} := 
	\{\rho \in h^{2+\beta}(\CC)\, |\, 
	\parallel \rho \parallel_{L^{\infty}(\CC)} <a \,\}$ 
	there is a unique, non-extendable solution 
	$$
	[t \mapsto \rho(t,\rho_0)] \in C^0([0,t^+),\UU^{a}_{\beta}) \cap 
	C^{\infty}((0,t^+),C^{\infty}(\CC))
	$$  
	of the initial value problem 
	\begin{equation} \label{initial.value.problem}
		\frac{\partial \rho}{\partial t}(x,t) = G(\rho_t)(x),  \,\,\,
		\textnormal{for} \,\,\, (x,t)\in \CC \times (0,t^+),\,\,\,
		\rho(x,0)=\rho_0(x), \,\,\, \textnormal{for}\,\,\,  x \in \CC,
	\end{equation}
	where $t^+=t^+(\rho_0)>0$ denotes the 
	``time of maximal existence'' of the respective solution 
	and where $G$ denotes the quasilinear differential operator from lines \eqref{normal.speed.3} and \eqref{G}. 
	Moreover, the map $[(t, \rho_0) \mapsto \rho(t,\rho_0)]$ defines a smooth local semiflow on $\UU^{a}_{\beta}$ in the sense of Section 12 in \cite{Amann.1993}.
\qed
\end{theorem} 
\noindent	  
Now, using some computations from Weiner's classical article \cite{Weiner} about the Willmore functional, we obtain the following counterpart to Lemma 3.1 in \cite{Escher.Mayer.Simonett.1998}.
\begin{lemma} \label{Lemma.3.1}
	The operator $-G:\VV^{a}_{\alpha}\equiv h^{4+\alpha}(\CC) \cap \UU^{a}_{\beta_0} \longrightarrow h^{\alpha}(\CC)$ from line (\ref{G}) 
	is $C^{\infty}$-smooth, 
	and its Fr\'echet derivative in 
	$\rho =0\in \VV^{a}_{\alpha}$ is precisely the uniformly elliptic linear operator
	$$
	- D_{\rho}G(\rho)\lfloor_{\rho=0} \equiv P(0) - D_{\rho}F(0)
	= \frac{1}{4} \, (\triangle_{\CC}+4) \circ (\triangle_{\CC}+2):
	h^{4+\alpha}(\CC) \longrightarrow  h^{\alpha}(\CC),  
	$$ 
	where $\triangle_{\CC}$ denotes the standard Beltrami-Laplace operator 
	on $\CC$ w.r.t. the Euclidean metric induced by the 
	injection $\CC \hookrightarrow \sphere^3$, i.e.    
	$\triangle_{\CC}=\triangle_{\rho}\lfloor_{\rho=0}$ in the notation of equation (\ref{triangle.rho}).
\end{lemma} 	
\proof 
We can immediately infer from Lemma \ref{Lemma.2.1} that the operator 
$G:\VV^{a}_{\alpha} \longrightarrow h^{\alpha}(\CC)$ is $C^{\infty}$-smooth
and thus continuously Fr\'echet-differentiable in $\VV^{a}_{\alpha}$. 
Moreover, since the manifold $\Gamma(0)$ is simply the Clifford-torus $\CC$ 
in $\sphere^3$ by the above construction, we may use Weiner's computation in \cite{Weiner}, pp. 24--25 and p. 34, and formula 
(\ref{laplacian.2}) below, in order to infer that the 
Fr\'echet derivative of the non-linear mean curvature operator 
$[\rho \mapsto H_{\rho}]$ of the manifold 
$\Gamma_{\rho}$ in $\rho=0$ is concretely given by 
\begin{equation} \label{strike.0}
	D_{\rho}H_{\rho}\lfloor_{\rho=0}= 
	\,-(\triangle_{\CC} + 4):
	h^{2+\alpha}(\CC) \longrightarrow h^{\alpha}(\CC).
\end{equation}  
Now, we compute the Fr\'echet derivative of the first term in 
(\ref{G}). To this end, we can proceed as in the proof of Lemma 3.1 in \cite{Escher.Mayer.Simonett.1998} and obtain from equation (\ref{strike.0}), and on account of $L_0\equiv 1$, $H_{\rho}\lfloor_{\rho=0} \equiv H_{\CC}\equiv 0$ and $|A^0_{\CC}|^2\equiv 2$ by means of the chain rule:
\begin{equation}  \label{first.strike}
	D_{\rho}\Big{(} \frac{L_{\rho}^2}
	{|A^0_{\rho}|^4} \,\triangle_{\rho}H_{\rho} \Big{)}\lfloor_{\rho=0}
	= -\frac{1}{4} \,\triangle_{\CC}\circ (\triangle_{\CC} + 4).    
\end{equation}
Moreover, in order to compute the Fr\'echet derivative of the second term in 
(\ref{G}), we employ again equation (\ref{strike.0}) and the equations $L_0\equiv 1$, $H_{\CC}\equiv 0$ and $K_{\CC}\equiv -1$ and obtain by means of the chain rule:  
\begin{eqnarray}  \label{second.strike}       
	D_{\rho}\Big{(} \frac{L_{\rho}^2}
	{|A^0_{\rho}|^4} \,2\,H_{\rho} \,(H_{\rho}^2-K_{\rho}) \Big{)}\lfloor_{\rho=0}= 
	\frac{1}{2} \, D_{\rho}H_{\rho}\lfloor_{\rho=0} \,(H_{0}^2-K_{0})
	= -\frac{1}{2}	\,(\triangle_{\CC} + 4).    \qquad                 
\end{eqnarray}	
Hence, adding formulae (\ref{first.strike}) and (\ref{second.strike}), we obtain on account of equation 
(\ref{G}):
\begin{eqnarray*} 
	-D_{\rho}G(\rho)\lfloor_{\rho=0} = 
	\frac{1}{4} \,\Big{(} \triangle_{\CC}^2 + 4\,\triangle_{\CC} + 
	2\, \triangle_{\CC} + 8 \Big{)} 
	= \frac{1}{4} \, (\triangle_{\CC}+4) \circ (\triangle_{\CC}+2),
\end{eqnarray*} 
just as asserted in this lemma. 
\qed          \\\\                    
\noindent
In view of the proof of Lemma \ref{Eigenspaces.DG} below,  
we recall here some differential geometric key-insights from 
\cite{Simons.1968} and \cite{Weiner}. 
First of all, since the normal bundle $N\CC$ is only one-dimensional 
and possesses the non-vanishing section $\nu_{\CC}$ of constant length $1$ w.r.t. the Euclidean metric $g_{\sphere^3}$, 
any smooth section $V \in \Gamma(N\CC)$ can be written in the form 
\begin{equation}  \label{section.normal.bundle}
	V = f_V \, \nu_{\CC}  \qquad \textnormal{on} \,\,\, \CC   
\end{equation}
for a uniquely determined smooth function 
$f_V:\CC \longrightarrow \rel$. Hence, we obtain a linear bijection 
\begin{equation}   \label{correspondence}
	\Gamma(N\CC)\ni V \longleftrightarrow f_V  \in C^{\infty}(\CC)  
\end{equation}  
between the set of smooth sections into $N\CC$ and functions of class  
$C^{\infty}(\CC)$. Moreover, the connection $\nabla^{\perp}$ in the normal bundle $N\CC$ - see Section 2.1 in \cite{Simons.1968} - maps sections of $N\CC$ into $N\CC$ again and maps the unit normal $\nu_{\CC}$ to $0$. 
Hence, defining covariant differentiation 
``$\nabla_{\partial_i}^{\CC}$'' of 
smooth functions $f\in C^{\infty}(\CC)$ in the direction of 
some locally defined partial derivative $\partial_i \in \Gamma(T\CC)$ 
by $\nabla_{\partial_i}^{\CC}(f):=\partial_i(f)$ and 
$\nabla_{\partial_i}^{\CC}\nabla_{\partial_j}^{\CC}(f) 
:=\partial_{ij}(f) - (\Gamma_{\CC})^k_{ij} \, \partial_{k}(f)$, 
where $(\Gamma_{\CC})^k_{ij}:=\gamma^k_{ij}(0)$, $i,j,k=1,2$ 
- in the notation of equation (\ref{triangle.rho}) - 
denote the Christoffel-symbols of the Euclidean metric 
induced by the injection $\CC \hookrightarrow \sphere^3$, 
we infer from the Leibniz-rule for linear connections:   
$$
\nabla^{\perp}_{\partial_i}(f \, \nu_{\CC})
= \nabla^{\CC}_{\partial_i}(f) \, \nu_{\CC} \quad 
\textnormal{on} \,\,\, \CC,
$$ 
and thus by definition of the Beltrami-Laplace operator,  
associated to the linear connection $\nabla^{\perp}$ in $N\CC$ 
and to the covariant derivative $\nabla^{\CC}$ on $\CC$:
\begin{equation}  \label{laplacian.2} 
	\triangle_{\CC}^{\perp} (f \, \nu_{\CC})
	= g^{ij}_{\CC} \, \nabla^{\perp}_{\partial_i} 
	\nabla^{\perp}_{\partial_j}(f \, \nu_{\CC})
	=  g^{ij}_{\CC} \, 
	\nabla^{\CC}_{\partial_i}\nabla^{\CC}_{\partial_j}(f) \, \nu_{\CC}
	= \triangle_{\CC}(f) \, \nu_{\CC}
\end{equation} 
for any function $f\in C^{\infty}(\CC)$, where $g_{\CC}$ and 
$\triangle_{\CC}$ denote the Euclidean metric on $T\CC$ and 
the standard Beltrami-Laplace operator on $\CC$, which are both induced by 
the injection $\CC \hookrightarrow \sphere^3$, 
i.e. $\triangle_{\CC}=\triangle_{\rho}\lfloor_{\rho=0}$ in the notation of equation (\ref{triangle.rho}). \\ 
Moreover, in view of the proof of Theorem \ref{Center.manifold} below, we recall and explain 
here some of Weiner's important observations from Sections 3 and 4 in \cite{Weiner}. First of all, we define two different elementary types of smooth sections of the normal bundle $N\CC$ of the Clifford-torus $\CC$ within $T\sphere^3$.
\begin{definition}  \label{parallel.vectorfields} 
	\begin{itemize} 
		\item[1)] We term elements $W$ of the $6$-dimensional Lie-algebra $\Omega$ of the isometry group
		$\textnormal{Iso}(\sphere^3)\equiv \textnormal{O}(4)$ ``Killing fields'' on $\sphere^3$. We denote the orthogonal projections of their restrictions $W\lfloor_{\CC}$ to the Clifford-torus $\CC$ into the normal bundle $N\CC$ by $(W\lfloor_{\CC})^N$, and we term the linear space of all such vector fields by $\Omega^N$. 
		\item[2)] We denote by $\xi$ the $4$-dimensional vector space of tangential projections of smooth parallel vector fields in $\rel^4$ into the tangent bundle of $\sphere^3$. We term the vector space of orthogonal projections $(Z\lfloor_{\CC})^N$ of vector fields $Z \in \xi$, restricted to $\CC$, into the normal bundle $N\CC$ by $\xi^N$.
		\qed  
	\end{itemize}           
\end{definition} 
\noindent
\begin{remark}  \label{conformal.group}
	We should note here, that $\Omega$ and $\xi$ 
	from Definition \ref{parallel.vectorfields} are both 
	vector subspaces of $\Gamma(T\sphere^3)$ with trivial intersection, and that they are both contained in the 
	Lie algebra of the entire conformal group $\textnormal{M\"ob}(\sphere^3)$, which is 
	$10$-dimensional because of the isomorphy 
	$\textnormal{M\"ob}(\sphere^3)\cong \textnormal{SO}^{+}(1,4)$ and since $\textnormal{SO}^{+}(1,n)$ is known to be $\frac{(n+1)\,n}{2}$-dimensional. 
	Hence, counting dimensions the direct sum 
	$\Omega\oplus \xi$ is exactly the vector subspace 
	of $\Gamma(T\sphere^3)$ consisting of all conformal 
	vector fields along $\sphere^3$; see pp. 30--33 in 
	\cite{Weiner}.     
	\qed 
\end{remark}
\noindent
Now, by Lemmata 5.1.3 and 5.1.7 in \cite{Simons.1968} and by
Lemmata 3.3, 3.4 and 3.5 in \cite{Weiner} we have the following basic results:
\begin{lemma}   \label{Simons.Weiner.1}
	\item[1)] Any vector field $W \in \Omega^N$ satisfies the 
	differential equation
	$$ 
	\triangle_{\CC}^{\perp}(W)= - 4 \, W. 
	$$  
	\item[2)] Any vector field $Z \in \xi^N$ satisfies the differential equation
	$$ 
	\triangle_{\CC}^{\perp}(Z)= - 2 \, Z. 
	$$  	  
	Moreover, both $\Omega^N$ and $\xi^N$ are $4$-dimensional 
	$\rel$-vector spaces, and their direct sum $\Gamma_{\CC}$ is the $8$-dimensional vector subspace of $\Gamma(N\CC)$ consisting of all ``normal conformal directions'' along $\CC$ in $T\sphere^3$.  
\qed  
\end{lemma}
\noindent
\begin{remark} \label{Jacobi.fields} 
The operator $-(\triangle_{\CC}^{\perp}+4)$ appearing in 
Lemma \ref{Simons.Weiner.1} is the ``Jacobi operator'' along the Clifford-torus, i.e. corresponds to the second variation of the area functional evaluated in the Clifford-torus 
$\CC$ w.r.t. sections of the normal bundle $N\CC$ 
along the Clifford-torus, and a smooth vector field 
$v\in \Gamma(N\CC)$ is termed a ``Jacobi field'' along $\CC$, 
if it satisfies $(\triangle_{\CC}^{\perp}+4)(v)=0$, i.e. if $v$ is contained in the eigenspace $\textnormal{Eig}_{-4}(\triangle_{\CC}^{\perp})$.
The first part of Lemma \ref{Simons.Weiner.1} therefore 
shows us, that orthogonal projections 
$W \in \Omega^N$ of Killing fields along the 
Clifford-torus $\CC$ into its normal bundle 
are Jacobi fields along $\CC$.  
\qed
\end{remark} 
\noindent
Combining Lemmata \ref{Lemma.3.1} and \ref{Simons.Weiner.1} with the proof of Corollary 1 in \cite{Weiner}, p. 34, we also infer the following important lemma.   
\begin{lemma}  \label{Eigenspaces.DG}
	The spectrum of the Fr\'echet derivative 
	$$
	D_{\rho}G(\rho)\lfloor_{\rho=0} \equiv P(0)-D_{\rho}F(0) 
	=- \frac{1}{4} \, (\triangle_{\CC}+4) \circ (\triangle_{\CC}+2)
	$$ 
	is discrete and non-positive, and its kernel is an $8$-dimensional $\rel$-vector subspace of $C^{\infty}(\CC)$, which corresponds to the vector subspace $\Gamma_{\CC}$ of $\Gamma(N\CC)$ of all normal conformal directions along the Clifford-torus $\CC$ from 
	Lemma \ref{Simons.Weiner.1} via the explicit linear bijection (\ref{correspondence}):
	\begin{equation}  \label{eigenspaces.DG}
	\textnormal{Ker}(D_{\rho}G(\rho)\lfloor_{\rho=0}) 
	\cong \Omega^N \oplus \xi^N = \Gamma_{\CC}.
	\end{equation}  
\end{lemma}
\proof 
Using the particularly simple form of the uniformly 
elliptic operator  
\begin{equation} \label{TC}
T_{\CC}:=- D_{\rho}G(\rho)\lfloor_{\rho=0}=
\frac{1}{4} \,(\triangle_{\CC}+4) \circ (\triangle_{\CC}+2)
\end{equation} 
one can prove as in Section 3 of \cite{Skorzinski.2015}, that $T_{\CC}$ is a compact perturbation of an isomorphism between $W^{4,2}(\CC)$ and $L^{2}(\CC)$ and thus a Fredholm operator of index $0$. Moreover, integration by parts and Cauchy-Schwarz' inequality immediately show, that there is some $c>0$ such that
$$
T_{\CC}+c \,\textnormal{Id}_{W^{4,2}(\CC)}: 
W^{4,2}(\CC) \longrightarrow  L^{2}(\CC)  
$$ 
is injective and thus a topological isomorphism. 
Since the composition 
\begin{equation}  \label{iota.A.C}
	\iota \circ (T_{\CC}+c \,\textnormal{Id}_{W^{4,2}(\CC)})^{-1}: 
	L^{2}(\CC) \longrightarrow L^{2}(\CC) 
\end{equation} 
of the inverse $(T_{\CC}+c \,\textnormal{Id}_{W^{4,2}(\CC)})^{-1}$ with
the compact embedding $\iota: W^{4,2}(\CC) \hookrightarrow  L^{2}(\CC)$ is a compact and selfadjoint operator, classical spectral theory guarantees, that the spectrum of $T_{\CC}$ consists of countably many real and isolated eigenvalues 
$-c < \nu_{1} < \nu_2 < \nu_3 <\ldots \in \rel$.
Moreover, just as in Lemma 3.2 of \cite{Weiner} we remark, that a real number $\lambda$ is an eigenvalue of the laplacian $\triangle_{\CC}$ in $W^{2,2}(\CC)$, 
if and only if the number 
$\frac{1}{4} \, (\lambda+4)\,(\lambda+2)= \frac{1}{4}(\lambda^2+6\lambda+8)$ 
is an eigenvalue of the linear operator 
$T_{\CC}=\frac{1}{4} \, (\triangle_{\CC}+4) \circ (\triangle_{\CC}+2)$ in $W^{4,2}(\CC)$. We note that the polynomial $f(\lambda):=\lambda^2+6\lambda+8$ 
is negative on the open interval $(-4,-2)$, and its roots are the two endpoints $-4$ and $-2$ of this interval.           
It is well-known that the eigenvalues of the Beltrami-laplacian on the Clifford-torus are the numbers $-2(m^2+n^2)$, for $m,n\in \nat_0$, using the fact that the 
Clifford-torus is isometric to the flat torus 
$\com/\Gamma^*$, where the lattice $\Gamma^*$ 
is spanned by the complex numbers $\frac{2\pi}{\sqrt 2}$ 
and $i\frac{2\pi}{\sqrt 2}$.
Since the intersection of the set $\{-2(m^2+n^2)\,|\,m,n\in \nat_0\,\}$ with the interval $(-4,-2)$ is empty, there consequently cannot exist any negative eigenvalues of $T_{\CC}=-D_{\rho}G(\rho)\lfloor_{\rho=0}$. 
Moreover, we see that the set $\{-2(m^2+n^2)\,|\,m,n\in \nat_0\,\}$ contains the two numbers $-2$ and $-4$, and 
we can therefore conclude that exactly  
\begin{equation}  \label{Eigenspaces} 
	\textnormal{Eig}_{0}(T_{\CC}) = \textnormal{Eig}_{-2}(\triangle_{\CC}) \oplus 
	\textnormal{Eig}_{-4}(\triangle_{\CC}). 
\end{equation}  
By formulae (\ref{section.normal.bundle}) -- (\ref{laplacian.2}) 
the eigenvalue problem ``$\triangle_{\CC}(f) \stackrel{!}= \lambda \, f$'', for $f \in C^{\infty}(\CC)$, is equivalent to the eigenvalue problem 
``$\triangle_{\CC}^{\perp}(V) \stackrel{!}= \lambda \, V$'', for $V \in N\CC$. We therefore recall, that by Lemma \ref{Simons.Weiner.1} there holds:
\begin{equation}   \label{Eigenvalue.xi.Omega}
	\textnormal{Eig}_{-2}(\triangle_{\CC}^{\perp}) \supset \xi^N \quad 
	\textnormal{and}  \quad    
	\textnormal{Eig}_{-4}(\triangle_{\CC}^{\perp}) \supset \Omega^N,
\end{equation}  
and that $\xi^N$ and $\Omega^N$ are two $4$-dimensional $\rel$-vector spaces. Moreover, using again the isometry between the Clifford-torus and the flat torus $\com/\Gamma^*$, one can easily show that the eigenspaces $\textnormal{Eig}_{-2}(\triangle_{\CC})$ 
and $\textnormal{Eig}_{-4}(\triangle_{\CC})$
are both $4$-dimensional, as well; compare here also to 
the proof of Corollary 1 on p. 34 in \cite{Weiner}. 
Hence, on account of statement (\ref{Eigenvalue.xi.Omega}) 
the eigenspaces $\textnormal{Eig}_{-2}(\triangle_{\CC}^{\perp})$
and $\textnormal{Eig}_{-4}(\triangle_{\CC}^{\perp})$
exactly coincide with the $4$-dimensional vector spaces 
$\xi^N$ and $\Omega^N$ respectively, i.e. there hold 
exactly equalities in (\ref{Eigenvalue.xi.Omega}).   
Hence, combining this result again with 
formulae (\ref{laplacian.2}) and (\ref{Eigenspaces}), the assertion of the lemma follows. 
\qed  \\\\
\noindent
Now, in the second part of this section we copy the methods of Theorem 1 and of the first part of Theorem 3.1 in \cite{Ruben.MIWF.III}, but using here parabolic H\"older spaces and parabolic Schauder theory instead of optimal $L^p-L^q$-estimates, in order to apply the regularity bootstrap method of Theorem 3 (ii) in \cite{Jakob_Moebius_2016}, via Schauder a-priori estimates.  
To this end, we consider evolution equations (\ref{Moebius.flow}) and (\ref{Moebius.flow.2}) for immersions $f_t:\Sigma \longrightarrow \rel^3$,  
and we recall from the author's article \cite{Jakob_Moebius_2016} as in Section 2 of \cite{Ruben.MIWF.III}, that the differential operator 
\begin{eqnarray}    \label{diff.operator}
	2\, \mid A^0_{f} \mid^{-4} \,\nabla_{L^2} \Will(f)	          
	\equiv \mid A^0_{f} \mid^{-4} 
	\,\Big{(} \triangle^{\perp}_{f} \vec H_{f} + Q(A^{0}_{f})(\vec H_{f}) \Big{)} = \nonumber   \\
	= \mid A^0_{f} \mid^{-4} 
	\Big{(} (\triangle_{f} \vec H_{f})^{\perp_f} + 2\,Q(A_f)(\vec H_f) - \frac{1}{2} \mid \vec H_f \mid^2\,\vec H_f   \Big{)}   
\end{eqnarray} 
arising on the right hand side of evolution equations 
(\ref{Moebius.flow}) and (\ref{Moebius.flow.2}) is not uniformly 
elliptic and that its leading term $(\triangle_{f} \vec H_{f})^{\perp_f}$ can be written as 
\begin{eqnarray}   \label{leading_term_of_invariant}
	(\triangle_{f} \vec H_{f})^{\perp_f}              
	=  g^{ij}_{f} \, g^{kl}_{f} 
	\, \nabla_i^{f} \nabla_j^{f} \nabla_k^{f} \nabla_l^{f}(f) 
	- g^{ij}_{f} \, g^{kl}_{f} \, 
	\langle \nabla_i^{f} \nabla_j^{f} \nabla_k^{f} \nabla_l^{f}(f), 
	\partial_m f \rangle\, g^{mr}_{f} \, \partial_r(f)  
\end{eqnarray}
at least locally, in local coordinates on $\Sigma$, for any $W^{4,2}$-immersion $f:\Sigma \longrightarrow \rel^3$, where $g_f:=f^*(g_{\textnormal{eu}})$
denotes the pullback-metric of the Euclidean metric of $\rel^3$. Applying ``DeTurck's trick'' as in Section 2 of \cite{Ruben.MIWF.III}, we add the globally well-defined differential operator of fourth order
\begin{eqnarray*} 
	\Tan_{F_0}(f):=  g^{ij}_{f} \,g^{kl}_{f}\, g^{mr}_{f} \, 
	\langle \nabla_i^{f} \nabla_j^{f} \nabla_k^{f} \nabla_l^{f}(f), 
	\partial_m f \rangle\, \partial_r f 
	- g^{ij}_{f} g^{kl}_{f}\, \nabla_i^{f} \nabla_j^{f} \,
	\big{(}(\Gamma_{F_0})^m_{kl} - (\Gamma_{f})^m_{kl}\big{)} \, \partial_m(f), 
\end{eqnarray*}
for some fixed smooth immersion $F_0:\Sigma \longrightarrow \rel^3$, to the right hand side of equation (\ref{diff.operator}) and obtain a quasilinear operator of fourth order, whose leading term is $g^{ij}_{f} \, g^{kl}_{f} \nabla_i^{F_0} \nabla_j^{F_0} \nabla_k^{F_0} \nabla_l^{F_0}(f)$, having a uniformly elliptic linearization in any umbilic-free $C^{4,\gamma}$-immersion $f:\Sigma \longrightarrow \rel^3$. We are thus led to consider here the evolution equation    
\begin{eqnarray}  \label{de_Turck_equation}
	\partial_t(f_t) = - \frac{1}{2} \mid A^0_{f_t} \mid^{-4} \,
	\big{(} 2\, \,\nabla_{L^2} \Will(f_t) + \Tan_{F_0}(f_t) \big{)}  =: \Mill_{F_0}(f_t),    \qquad
\end{eqnarray} 
for some arbitrarily fixed $C^{\infty}-$smooth  
immersion $F_0:\Sigma \longrightarrow \rel^3$, 
where the right-hand side 
``$\Mill_{F_0}(f_t)$'' of equation (\ref{de_Turck_equation}) 
can be expressed in local coordinates on $\Sigma$ by: 
\begin{eqnarray}  \label{D_F_0}
	\Mill_{F_0}(f_t)(x) = - \frac{1}{2} \mid A^0_{f_t} \mid^{-4} \,
	g^{ij}_{f_t} \, g^{kl}_{f_t}\, \nabla_i^{f_t} \nabla_j^{f_t} 
	\nabla_k^{F_0} \nabla_l^{F_0}(f_t)(x) 
	+ B(x,D_xf_t,D^2_x f_t, D_x^3 f_t),     
\end{eqnarray} 
for $(x,t) \in \Sigma \times [0,T]$. Here, the symbols 
$D_xf_t, D_x^2f_t, D_x^3 f_t$ abbreviate the matrix-valued functions $(\partial_{1}f_t,\partial_{2}f_t)$, 
$(\nabla^{F_0}_{ij}f_t)_{i,j \in \{1,2\}}$ and 
$(\nabla^{F_0}_{ijk}f_t)_{i,j,k \in \{1,2\}}$, 
and $B: \Sigma \times \rel^6 \times \rel^{12} \times 
\rel^{24} \to \rel^3$ is a globally defined function, 
whose $3$ components are rational functions in 
their $42$ real variables. Following the lines of the 
author's articles \cite{Jakob_Moebius_2016} and \cite{Ruben.MIWF.III}, 
we will collect some basic properties of the linearization of equation (\ref{de_Turck_equation}) or equivalently of equation  
\begin{eqnarray}  \label{de_Turck_equation_2}
	\partial_t(f_t) = - \frac{1}{2} \mid A^0_{f_t} \mid^{-4} \,
	g^{ij}_{f_t} \, g^{kl}_{f_t}
	\, \nabla_i^{f_t} \nabla_j^{f_t} \nabla_k^{F_0} \nabla_l^{F_0}(f_t) + B(\,\cdot\,,D_xf_t,D_x^2f_t,D_x^3f_t), 
\end{eqnarray} 
in any family of $C^{4,\gamma}$-immersions
$f_t:\Sigma \longrightarrow \rel^3$, which is sufficiently close to a given smooth flow line of umbilic-free immersions 
in the parabolic H\"older space 
$C^{4+\gamma,1+\frac{\gamma}{4}}(\Sigma \times [0,T],\rel^3)$,
below in Proposition \ref{Psi.of.class.C_1}.
Now, we fix some umbilic-free immersion 
$F_0 \in C^{\infty}(\Sigma,\rel^3)$, and we denote by 
$T_{\textnormal{max}}(F_0)>0$ the maximal time, such 
that the corresponding unique smooth flow line 
$\{\PP(\,\cdot\,,0,F_0)\}$ of the MIWF (\ref{Moebius.flow})
exists on $\Sigma \times [0,T_{\textnormal{max}}(F_0))$, 
starting to move in $F_0$ at time $t=0$. 
We recall from the proof of Theorem 1 in \cite{Jakob_Moebius_2016}, 
that there is a unique smooth family of smooth diffeomorphisms 
$\phi_t^{F_0}:\Sigma \longrightarrow \Sigma$ with 
$\phi_0=\textnormal{Id}_{\Sigma}$, 
such that the reparametrization $\{\PP(t,0,F_0)\circ \phi_t^{F_0}\}$ satisfies evolution equation (\ref{de_Turck_equation}) respectively (\ref{de_Turck_equation_2}), i.e. the equation
$$
\partial_t(f_t) = \Mill_{F_0}(f_t)  \qquad \textnormal{on} \quad \Sigma \times [0,T],
$$ 
for every final time $0<T<T_{\textnormal{max}}(F_0)$. 
Now, we fix some arbitrary final time  $0<T<T_{\textnormal{max}}(F_0)$, some $\gamma \in (0,1)$ and some open neighborhood $W_{F_0,T,\gamma}$ of the above flow line $\{\PP(t,0,F_0)\circ \phi_t^{F_0}\}$ 
of evolution equation (\ref{de_Turck_equation}) in the 
Banach space
$$
X_T\equiv X_{T,\gamma}:= 
C^{4+\gamma,1+\frac{\gamma}{4}}(\Sigma \times [0,T],\rel^3),
$$ 
and we shall follow the strategy of \cite{Ruben.MIWF.III} and 
\cite{Spener}: using the fact that the restriction of elements 
of $X_{T,\gamma}$ at time $t=0$:
$$
r_0:X_{T,\gamma} \longrightarrow  C^{4,\gamma}(\Sigma,\rel^3)
$$
is a linear and continuous operator, i.e. that 
the ``trace'' of the Banach space $X_T$ at time $t=0$ is exactly 
\begin{equation}  \label{Trace.X.T}
\textnormal{Trace}(X_T)=C^{4,\gamma}(\Sigma,\rel^3),
\end{equation} 
and considering the continuous, non-linear product operator
$$
\Psi^{F_0,T}: W_{F_0,T,\gamma} \subset X_T
\longrightarrow C^{4,\gamma}(\Sigma,\rel^3) \times 
C^{\gamma,\frac{\gamma}{4}}(\Sigma \times [0,T],\rel^3) 
=:Y_{T,\gamma} \equiv Y_{T}
$$
defined by 
\begin{eqnarray}  \label{Psi}
\Psi^{F_0,T}(\{f_t\}_{t\in [0,T]})
:=(f_0,\{\partial_t(f_t) - \Mill_{F_0}(f_t)\}_{t\in [0,T]}). 
\end{eqnarray}
Now using (\ref{Trace.X.T}) and following the 
proofs of Theorem 2.1 in \cite{Ruben.MIWF.III}
and of Theorem 2 in \cite{Jakob_Moebius_2016}, but 
substituting here the role of Proposition 2.1 in \cite{Ruben.MIWF.III} by the powerful combination  
of Propositions 1 and 2 and 
Corollaries 2 and 3 in \cite{Jakob_Moebius_2016}, 
we can easily prove the following counterpart 
of Theorem 2.1 in \cite{Ruben.MIWF.III} 
within the setting of ``parabolic Schauder Theory'', 
aiming at the proof of basic properties of the Fr\'echet derivative of the operator $\Psi^{F_0,T}$ in (\ref{Psi}), evaluated in the smooth solution $\{\PP(t,0,F_0)\circ \phi_t^{F_0}\}_{t\in [0,T]}$ of the modified MIWF-equation (\ref{de_Turck_equation_2}).
\begin{proposition}           \label{Psi.of.class.C_1} 
	Let $\Sigma$ be a smooth compact torus and
	$F_0:\Sigma \longrightarrow \rel^3$ a $C^{\infty}$-smooth 
	and umbilic-free immersion, and let 
	$0<T<T_{\textnormal{max}}(F_0)$ and $\gamma\in (0,1)$ be 
	chosen arbitrarily, where $T_{\textnormal{max}}(F_0)>0$
	denotes the time of maximal existence of 
	the flow line $\{\PP(\,\cdot\,,0,F_0)\}$ of 
	the MIWF (\ref{Moebius.flow}).  
	There is a sufficiently small open neighborhood $W_{F_0,T,\gamma}$ about the smooth solution $\{\PP(t,0,F_0)\circ \phi_t^{F_0}\}_{t\in [0,T]}$ 
	of the modified MIWF-equation (\ref{de_Turck_equation_2}) 
	in the space 
	$X_T = C^{4+\gamma,1+\frac{\gamma}{4}}(\Sigma \times [0,T],\rel^3)$, such that the following three statements 
	hold: 
	\begin{itemize} 
		\item[1)] The map $\Psi^{F_0,T}: W_{F_0,T,\gamma} \longrightarrow Y_T$, defined in line (\ref{Psi}), is of class $C^{1}$ on the open subset $W_{F_0,T,\gamma}$ of the Banach space $X_T$.
		\item[2)] In any fixed element $\{f_t\} \in W_{F_0,T,\gamma}$ the Fr\'echet derivative of the second component of $\Psi^{F_0,T}$ is a linear, uniformly parabolic operator of order $4$ whose leading operator acts diagonally, i.e. on each component of $f=\{f_t\}$ separately:
		\begin{eqnarray} \label{Frechet.of.Psi}
			(D(\Psi^{F_0,T})_2(f)).(\eta) \equiv 
			\partial_t(\eta) - D(\Mill_{F_0})(f).(\eta)= \qquad \\
			= \partial_t(\eta) + \frac{1}{2} \mid A^0_{f_t} \mid^{-4} \,
			g^{ij}_{f_t} \, g^{kl}_{f_t} \, \nabla^{F_0}_{ijkl}(\eta)              
			+  B_3^{ijk} \cdot \nabla^{F_0}_{ijk}(\eta)   
			+ B_2^{ij} \cdot \nabla^{F_0}_{ij}(\eta) 
			+ B_1^{i} \cdot \nabla^{F_0}_{i}(\eta)          \nonumber
		\end{eqnarray} 
		on $\Sigma \times [0,T]$, for any element $\eta=\{\eta_t\}$ of the tangent space $T_{f}W_{F_0,T,\gamma}=X_T$. Here, the coefficients $\mid A^0_{f_t} \mid^{-4} \,
		g^{ij}_{f_t} \, g^{kl}_{f_t}$ of the leading order term are of class $C^{2+\gamma,\frac{2+\gamma}{4}}(\Sigma \times [0,T], \rel^3)$, $B^{ij}_2$ and $B^i_1$ are the coefficients of $\textnormal{Mat}_{3,3}(\rel)-$valued, contravariant tensor fields 
		of degrees $2$ and $1$, which depend on $x$, 
		$D_xf_t,D^2_xf_t,D^3_xf_t$ and on $D^4_xf_t$ and are of class $C^{\gamma,\frac{\gamma}{4}}(\Sigma \times [0,T], \textnormal{Mat}_{3,3}(\rel))$.
		Finally $B_{3}^{ijk}$ are the coefficients of a $\textnormal{Mat}_{3,3}(\rel)-$valued, 
		contravariant tensor field of degree $3$, which depends on $x$, $D_xf_t$ and $D^2_xf_t$ only and is of class 
		$C^{2+\gamma,\frac{2+\gamma}{4}}(\Sigma \times [0,T],\textnormal{Mat}_{3,3}(\rel))$.
		\item[3)] The Fr\'echet derivative of $\Psi^{F_0,T}$ yields a topological isomorphism 
		$$
		D\Psi^{F_0,T}(f): T_{f}W_{F_0,T,\gamma}=X_T \stackrel{\cong} \longrightarrow Y_T
		$$
		in any fixed family of immersions 
		$f \equiv \{f_t\} \in W_{F_0,T,\gamma}$.
	\end{itemize}
	\qed
\end{proposition} 
\noindent 
Combining Proposition \ref{Psi.of.class.C_1} with the proof of Theorem 3 (ii) in \cite{Jakob_Moebius_2016} via Proposition 3 in \cite{Jakob_Moebius_2016}, we obtain the following theorem, similarly to Theorem 3.1 (i) in \cite{Ruben.MIWF.III}.
\begin{theorem}   \label{Frechensbergo}
	Let $\Sigma$ be a smooth compact torus and
	$F_0:\Sigma \longrightarrow \rel^3$ a $C^{\infty}$-smooth 
	and umbilic-free immersion, and let 
	$0<T<T_{\textnormal{max}}(F_0)$ and $\gamma\in (0,1)$ be 
	chosen arbitrarily, where $T_{\textnormal{max}}(F_0)>0$
	denotes the time of maximal existence 
	of the flow line $\{\PP(\,\cdot\,,0,F_0)\}$ of the MIWF 
	(\ref{Moebius.flow}).   
	\begin{itemize}
		\item[1)] There is some small $\rho=\rho(\Sigma,F_0,T,\gamma)>0$ such that for every initial immersion $F \in C^{4,\gamma}(\Sigma,\rel^3)$ with $\parallel F-F_0 \parallel_{C^{4,\gamma}(\Sigma,\rel^3)}<\rho$ 
		there is a unique and classical solution
		$\{\PP^*(t,0,F)\}_{t \in [0,T]}$ of the 
		``DeTurck modification'' (\ref{de_Turck_equation}) respectively (\ref{de_Turck_equation_2}) of the MIWF (\ref{Moebius.flow}) 
		in the Banach space $X_T=C^{4+\gamma,1+\frac{\gamma}{4}}(\Sigma \times [0,T],\rel^3)$, starting to move in the immersion $F$ at time $t=0$, and the resulting evolution operator 
		\begin{equation}  \label{solution.operator.2}
			\PP^*(\,\cdot\,,0,\,\cdot\,):
			B_{\rho}^{4,\gamma}(F_0) \subset C^{4,\gamma}(\Sigma,\rel^3) \longrightarrow X_T, 
		\end{equation}
		mapping any element $F$ of the open ball $B_{\rho}^{4,\gamma}(F_0)$ 
		about $F_0$ in $C^{4,\gamma}(\Sigma,\rel^3)$ to the solution $\{\PP^*(t,0,F)\}_{t \in [0,T]}$, is of class $C^{1}$.
		\item[2)] If the initial immersion 
		$F \in B_{\rho}^{4,\gamma}(F_0)$ from part (1) 
		is additionally of class $C^{\infty}(\Sigma,\rel^3)$, 
		then the resulting classical solution $\{\PP^*(t,0,F)\}_{t \in [0,T]}$ of evolution equation (\ref{de_Turck_equation_2}) from line (\ref{solution.operator.2}) is 
		of class $C^{\infty}(\Sigma \times [0,T],\rel^3)$, and furthermore there is a smooth family of $C^{\infty}$-smooth diffeomorphisms $\psi^F_t:\Sigma \to \Sigma$, with $\psi^F_0=\textnormal{Id}_{\Sigma}$, such that the composition $\PP^*(t,0,F) \circ \psi_t^{F}$ solves evolution equation (\ref{Moebius.flow}) on $\Sigma \times [0,T]$, i.e. such that there holds: 
		$$ 
		\PP^*(t,0,F) \circ \psi_t^{F} = \PP(t,0,F)  \quad  
		\textnormal{on} \,\,\, \Sigma, \quad \forall \,t \in [0,T]. 
		$$   		 	
	\end{itemize} 
\end{theorem}
\proof 
\begin{itemize} 
	\item[1)] Here we can argue exactly as in the proof of Theorem 3.1 (i) in \cite{Ruben.MIWF.III}.
	We assume that $F_0:\Sigma \longrightarrow \rel^3$ is a $C^{\infty}$-smooth and umbilic-free immersion which produces 
	a maximal smooth flow line $\{\PP(\,\cdot\,,0,F_0)\}$
	of the MIWF, starting in $F_0$ at time $t=0$. 
	Hence, the proof of Theorem 1 in \cite{Jakob_Moebius_2016} shows, that there is a unique smooth family of smooth diffeomorphisms 
	$\phi_t^{F_0}:\Sigma \longrightarrow \Sigma$ with 
	$\phi_0^{F_0}=\textnormal{Id}_{\Sigma}$, 
	such that the reparametrization 
	$\{\PP(t,0,F_0)\circ \phi_t^{F_0}\}_{t \geq 0}$ is the unique and maximal smooth solution of evolution equation 
	(\ref{de_Turck_equation}) on $\Sigma$. 
	Moreover, we know from Proposition \ref{Psi.of.class.C_1} above, that there is some open neighborhood $W_{F_0,T,\gamma}$ 
	of the smooth solution $\{\PP(t,0,F_0)\circ \phi_t^{F_0}\}_{t \in [0,T]}$ of equation (\ref{de_Turck_equation}) respectively (\ref{de_Turck_equation_2}) in the space $X_T$, such that the operator $\Psi^{F_0,T}$ from line (\ref{Psi}) is a $C^{1}$-map from $W_{F_0,T,\gamma}$ to $Y_T$, whose Fr\'echet derivative in the particular element 
	$\{\PP(t,0,F_0)\circ \phi_t^{F_0}\}_{t \in [0,T]}\in 
	W_{F_0,T,\gamma}$ is a topological isomorphism between $X_T$ and $Y_T$. Noting also that there holds 
	$$
	\Psi^{F_0,T}(\{\PP(t,0,F_0)\circ \phi_t^{F_0}\})
	=(F_0,0)\in Y_T,
	$$ 
	by definition of the operator $\Psi^{F_0,T}$ in (\ref{Psi}) and since $\{\PP(t,0,F_0)\circ \phi_t^{F_0}\}_{t\in [0,T]}$ solves equation (\ref{de_Turck_equation}),  
	we infer from the inverse mapping theorem for non-linear $C^{1}$-operators, that there is some small open ball 
	$B_{\rho}((F_0,0)) \subset Y_T$, with $\rho=\rho(F_0,T,\gamma)>0$, 
	and an appropriate further open neighborhood $W_{F_0,T,\gamma}^* \subset W_{F_0,T,\gamma}$ about the smooth solution $\{\PP(t,0,F_0)\circ \phi_t^{F_0}\}_{t \in [0,T]}$ of equation (\ref{de_Turck_equation_2}) in $X_T$, such that 
	\begin{equation}  \label{Psi.W.star}
		\Psi^{F_0,T}: W_{F_0,T,\gamma}^* \stackrel{\cong}\longrightarrow  B_{\rho}((F_0,0)) 
	\end{equation}
	is a $C^{1}$-diffeomorphism. Hence, by definition of the map $\Psi^{F_0,T}$ and by Theorem 2.2 (ii) of \cite{Ruben.MIWF.III} the restriction of the inverse mapping $(\Psi^{F_0,T})^{-1}$ from line (\ref{Psi.W.star}) to the product $B_{\rho}^{4,\gamma}(F_0) \times \{0\} \subset B_{\rho}((F_0,0))$ yields exactly the evolution operator of the parabolic evolution equation (\ref{de_Turck_equation}): 
	\begin{eqnarray}  \label{C1.Evolution.operator}
		X_T \supset W_{F_0,T,\gamma}^*\ni 
		\{\PP^*(t,0,F)\}_{t \in [0,T]} = (\Psi^{F_0,T})^{-1}((F,0))                    \\
		\forall \, F \in B_{\rho}^{4,\gamma}(F_0) 
		\subset C^{4,\gamma}(\Sigma,\rel^3),           \nonumber
	\end{eqnarray}   
	and it consequently has to be of class $C^{1}$ as an operator from $C^{4,\gamma}(\Sigma,\rel^3)$ to $X_T$. 
	Here, $\{\PP^*(t,0,F)\}_{t \in [0,T]}$ 
	denotes the restriction of the unique maximal solution 
	$\{\PP^*(t,0,F)\}_{t \in [0,t^+(F)]}$ of equation (\ref{de_Turck_equation}) from Theorem 2.2 (ii) in 
	\cite{Ruben.MIWF.III} to the interval $[0,T]$, noting 
	that $C^{4,\gamma}(\Sigma,\rel^3) \subset  W^{4-\frac{4}{p},p}(\Sigma,\rel^3)$ and 
	that $C^{4+\gamma,1+\frac{\gamma}{4}}(\Sigma \times [0,T],\rel^3)\subset 
	W^{1,p}([0,T];L^p(\Sigma,\rel^3)) \cap L^p([0,T];W^{4,p}(\Sigma,\rel^3))$.
	\item[2)] Here we can partially argue as in the proof of Theorem 3 (ii) in \cite{Jakob_Moebius_2016}. Firstly, we choose some $C^{\infty}$-smooth immersion $F:\Sigma \longrightarrow \rel^3$ in the small ball $B_{\rho}^{4,\gamma}(F_0)$ about the fixed immersion $F_0$ from line (\ref{solution.operator.2}). The fact that the family of immersions 
	$\{f^*_t\}_{t\in [0,T]}:=\{\PP^*(t,0,F)\}_{t\in [0,T]}$ - produced by the solution operator in line (\ref{solution.operator.2}) - 
	solves equation (\ref{de_Turck_equation}), means by equation (\ref{de_Turck_equation_2}) that there holds 
	\begin{eqnarray} \label{fstar_solves_de_Turck}
		\partial_t(f^*_t)(x) 
		= -\frac{1}{2} \mid A^0_{f^*_t}(x) \mid^{-4} \,
		g^{ij}_{f^*_t}(x) \, g^{kl}_{f^*_t}(x) \, \nabla^{F_0}_{ijkl}(f^*_t)(x) \nonumber\\
		- \FF(x,D_xf^*_t(x),D_x^2f^*_t(x),D_x^3f^*_t(x))  \qquad
	\end{eqnarray}
	for $(x,t) \in \Sigma \times [0,T]$, for some globally defined function $\FF:\Sigma \times \rel^{6} \times \rel^{12} \times \rel^{24} \longrightarrow \rel^3$, see equations (24) and (45) 
	in \cite{Jakob_Moebius_2016}. We note, that the corresponding linear differential operator 
	\begin{eqnarray*}
		L_F:=\partial_t + \frac{1}{2} \mid A^0_{f^*_t} \mid^{-4} \,
		g^{ij}_{f^*_t} \, g^{kl}_{f^*_t} \nabla^{F_0}_{ijkl}:
		C^{4+\gamma,1+\frac{\gamma}{4}}(\Sigma \times [0,T],\rel^3)
		\longrightarrow C^{\gamma, \frac{\gamma}{4}}(\Sigma \times [0,T],\rel^3)
	\end{eqnarray*}
	is uniformly parabolic and of diagonal form, 
	and that its coefficients are of class 
	$C^{2+\gamma,\frac{2+\gamma}{4}}(\Sigma \times [0,T])$, on 
	account of $\{f^*_t\} \in W_{F_0,T,\gamma}^*$, thus meets all conditions of Proposition 1 in \cite{Jakob_Moebius_2016}, 
	for some appropriate constant $\Lambda \geq 1$, 
	depending on $\Sigma$, $F_0$, $T$ and on the size of the narrow open tube $W_{F_0,T,\gamma}^*$ from line (\ref{C1.Evolution.operator}) about the flow line $\{\PP^*(t,0,F_0)\}_{t\in [0,T]}$. Moreover, 
	the known $C^{4+\gamma,1+\frac{\gamma}{4}}(\Sigma \times [0,T])$-regularity of $\{f^*_t\}$ - on account of the first part of this theorem - implies, that the composition $\FF(\,\cdot\,,D_xf^*,D_x^2f^*,D_x^3f^*)$ 
	is of class $C^{1+\gamma,\frac{1+\gamma}{4}}(\Sigma \times [0,T],\rel^3)$. Since the family $\{f^*_t\}$ solves the reformulation  
	\begin{eqnarray} \label{booty} 
		L_F(f^*)(x,t) = - \FF(x,D_xf^*_t(x),D_x^2f^*_t(x),D_x^3f^*_t(x))  \quad \textnormal{on} \,\,\, \Sigma \times [0,T]\quad \\
		\textnormal{with} \,\,\, f^*_0=F \,\,\,\textnormal{on} \,\,\, \Sigma, \quad \nonumber 
	\end{eqnarray} 
	of equation (\ref{de_Turck_equation_2}) classically on $\Sigma \times [0,T]$, since the operator $L_F$ and the right hand side of equation (\ref{booty}) meet all requirements of 
	Proposition 3 in \cite{Jakob_Moebius_2016} for $k=1$,
	and since we also assume that the initial immersion $F$ in problem (\ref{booty}) is of class $C^{\infty}(\Sigma,\rel^3)$, 
	we may apply the Schauder Regularity Theorem, 
	Proposition 3 in \cite{Jakob_Moebius_2016}, and infer that $\{f^*_t\}$ is of class $C^{5+\gamma,\frac{5+\gamma}{4}}(\Sigma \times [0,T],\rel^3)$. Thus, $L_F$, the right hand side of equation (\ref{booty}) and the initial surface $F$ turn out to satisfy all requirements of Proposition 3 in \cite{Jakob_Moebius_2016} for $k=2$. We can therefore repeat the above argument and obtain by induction that $\{f^*_t\}$ is actually of class 
	$C^{4+k+\gamma,1+\frac{k+\gamma}{4}}(\Sigma \times [0,T],\rel^3)$
	for any $k \in \nat_0$, i.e. that $\{f^*_t\} \in C^{\infty}(\Sigma \times [0,T],\rel^3)$,
	just as asserted in this part of the theorem.
	The last assertion of the theorem now follows immediately from the proven $C^{\infty}$-regularity of $\{f^*_t\}$, combined with a simple comparison of the evolution equations (\ref{Moebius.flow}) and (\ref{de_Turck_equation}), just as elaborated in the proof of Theorem 1 in \cite{Jakob_Moebius_2016}.
	\qed 
\end{itemize} 

\section{Proofs of Theorems \ref{Center.manifold}--\ref{main.result.2}}	
\underline{Proof of Theorem \ref{Center.manifold}}\\\\
\noindent
i) Without loss of generality we can assume 
that $F^*:\Sigma \stackrel{\cong}\longrightarrow 
\frac{1}{\sqrt 2}(\sphere^1 \times \sphere^1) \equiv \CC$ is a diffeomorphic parametrization of exactly the Clifford-torus in $\sphere^3$. Following the proof of Theorem 1.2 in \cite{Escher.Mayer.Simonett.1998} respectively the lines of Section 6 in \cite{Escher.Simonett.1998}, 
we use Lemmata \ref{Lemma.2.1}, \ref{Lemma.3.1} and \ref{Eigenspaces.DG} and Theorem \ref{Short.time.existence},
and we copy the procedure in Section 4 of \cite{Simonett.1995}, 
in order to construct an invariant ``center manifold'' $\Mill^c$ for the flow equation (\ref{normal.speed.3}) 
as a graph over the space $\textnormal{Ker}(T_{\CC})$.     
To this end, we use again the fact that the composition 
$\iota \circ (T_{\CC}+c \,\textnormal{Id}_{W^{4,2}(\CC)})^{-1}$ 
from line (\ref{iota.A.C}) is compact and selfadjoint, 
and thus classical spectral theory yields, that the 
Hilbert space $L^{2}(\CC)$ 
can be decomposed orthogonally
w.r.t. $\langle \,\cdot \,, \,\cdot \,\rangle_{L^2(\CC)}$   
into the finite dimensional eigenspaces 
$\textnormal{Ker}(T_{\CC})\subset C^{\infty}(\CC)$
and $\textnormal{Eig}_{\mu_j}(T_{\CC})\subset C^{\infty}(\CC)$ of $T_{\CC}$, for the positive eigenvalues $\mu_j$ of the linear operator $T_{\CC}$ in $L^2(\CC)$, which means precisely:
\begin{eqnarray}    \label{decomposition.a}
	L^{2}(\CC) = \textnormal{Ker}(T_{\CC}) 
	\oplus \,\overline{\bigoplus_{\mu >0} \textnormal{Eig}_{\mu}(T_{\CC})}^{L^{2}},
\end{eqnarray} 
where we used the first statement of Lemma \ref{Eigenspaces.DG}.
Now, by the second statement of Lemma \ref{Eigenspaces.DG} we 
can choose an $L^2(\CC)$-orthonormal system of $8$ eigenfunctions $\{Y_k\}_{k=1,\ldots,8}$ of $T_{\CC}$ 
in $\textnormal{Ker}(T_{\CC})$, 
and we define the continuous linear projection 
\begin{equation}  \label{Projection}
\pi^c:= \sum_{k=1}^8 \langle \,\cdot \,,Y_k \rangle_{L^2(\CC)} \, Y_k: 
L^2(\CC) \longrightarrow \textnormal{Ker}(T_{\CC})
\end{equation}
of $L^2(\CC)$ onto the ``center subspace'' 
$\textnormal{Ker}(T_{\CC})$ 
of the linear operator $T_{\CC}$, which is 
orthogonal w.r.t. $\langle \,\cdot \,, \,\cdot \,\rangle_{L^2(\CC)}$ on account of (\ref{decomposition.a}). The restrictions of $\pi^c$ 
in (\ref{Projection}) to the Banach spaces $h^{r}(\CC)$, for $r>0$, are still continuous, linear projections onto the finite dimensional subspace $\textnormal{Ker}(T_{\CC})$, 
and we therefore obtain as in Section 4 of 
\cite{Simonett.1995} a unique decomposition of $h^{r}(\CC)$ into two closed linear subspaces: 
\begin{equation}  \label{decomposition.b} 
	h^{r}(\CC) = 
	\textnormal{range}(\pi^c\lfloor_{h^{r}(\CC)}) 
	\oplus \textnormal{ker}(\pi^c\lfloor_{h^{r}(\CC)}) 
	\equiv \textnormal{Ker}(T_{\CC}) \oplus  h_s^{r}(\CC),  
\end{equation} 
for any fixed $r>0$, especially for $r =\alpha$ or 
$r=4+\alpha$. Obviously, $T_{\CC}$ is a symmetric
operator in $L^2(\CC)$ on account of the concrete formula (\ref{TC}), which implies together with (\ref{Projection}):
$$
\pi^c(T_{\CC}(f))= 
\sum_{k=1}^8 \langle T_{\CC}(f), Y_k \rangle_{L^2(\CC)} \, Y_k
=\sum_{k=1}^8 \langle\, f, T_{\CC}(Y_k) \rangle_{L^2(\CC)} \, Y_k = 0
$$ 
for every function $f \in W^{4,2}(\CC)$.  
In particular, there holds therefore:  
$$
T_{\CC} \circ \pi^c\lfloor_{h^{4+\alpha}(\CC)} 
=0= \pi^c\lfloor_{h^{\alpha}(\CC)} \circ T_{\CC} \quad 
\textnormal{on} \quad h^{4+\alpha}(\CC),
$$ 
i.e. that $T_{\CC}$ descends to a direct sum $T_{\CC}^c \oplus T_{\CC}^s$ of linear operators, which respects the direct sum decomposition $\textnormal{Ker}(T_{\CC}) \oplus h_s^{4+\alpha}(\CC)$ of $h^{4+\alpha}(\CC)$ in (\ref{decomposition.b}). 
Obviously, by (\ref{decomposition.a}) and (\ref{decomposition.b}) the restriction $T_{\CC}^s$ of 
$T_{\CC}$ to $h_s^{4+\alpha}(\CC)$ has only positive eigenvalues $0<\mu_1<\mu_2<\ldots$. As in \cite{Escher.Mayer.Simonett.1998}, \cite{Escher.Simonett.1998} and \cite{Simonett.1995}, we shall consider besides 
the projection $\pi^c$ in (\ref{Projection}) the 
$L^2$-orthogonal projection
$\pi^s:=\textnormal{id}_{L^2(\CC)} - \pi^c$
of $L^2(\CC)$ onto $\textnormal{ker}(\pi^c)$, 
whose restriction to $h^{r}(\CC)$ 
maps $h^{r}(\CC)$ onto the ``stable subspace'' 
$h_s^{r}(\CC)\subset h^{r}(\CC)$ 
in (\ref{decomposition.b}) w.r.t. $T_{\CC}$. 
Now, following Section $2$ in \cite{Simonett.1994}, 
Section $4$ of \cite{Simonett.1995} or Section 9.2.1 in \cite{Lunardi}, we write equation (\ref{normal.speed.3}) 
in the equivalent form: 
\begin{equation}  \label{reform}
	\partial_t \rho_t + T_{\CC}(\rho_t) \equiv
	\partial_t \rho_t + (P(0)-D_{\rho}F(0)).(\rho_t)= g(\rho_t),
\end{equation}  
with $g(\rho):= (P(0)-P(\rho)).(\rho)+ F(\rho)-D_{\rho}F(0).(\rho)$, satisfying $g(0)=0$ and $Dg(0)=0$, because of $F(0)=0$ by Lemma \ref{Lemma.2.1}. Moreover, as in formula (4.21) in \cite{Simonett.1995} or as in Section 9.2.1 in \cite{Lunardi}, decomposition (\ref{decomposition.b}) yields the equivalent formulation  
\begin{eqnarray} \label{decomposition.2}
\partial_t x_t + T_{\CC}^c(x_t) = \pi^cg(x_t,y_t)                \\
\partial_t y_t + T_{\CC}^s(y_t) = \pi^sg(x_t,y_t)              \nonumber
\end{eqnarray}
of equation (\ref{reform}) respectively of equation (\ref{normal.speed.3}) as a coupled system, for two seperate functions $x:[0,T]\longrightarrow \textnormal{Ker}(T_{\CC})$ and $y:[0,T]\longrightarrow h_s^{4+\alpha}(\CC)$.
Furthermore, on account of Lemmata \ref{Lemma.2.1}, 
\ref{Lemma.3.1} and \ref{Eigenspaces.DG} and Theorem \ref{Short.time.existence}, 
we may apply Theorem 4.1 in \cite{Simonett.1995}. This theorem guarantees us, that for some fixed $m\in \nat$ there exists a neighborhood $U=U(m)$ of $0$ in $\textnormal{Ker}(T_{\CC})$ and a function 
\begin{equation}  \label{gamma} 
\gamma \in C^m(U,h_s^{4+\alpha}(\CC)),\,\,\textnormal{with} \,\, \gamma(0)=0 \,\, \textnormal{and}\,\, D\gamma(0)=0
\end{equation} 
such that $\Mill^c := \textnormal{graph}(\gamma)$ is a ``locally invariant center manifold'' for the semiflow generated by the unique maximal solutions of equation (\ref{normal.speed.3}) respectively of the coupled system (\ref{decomposition.2}), provided by Theorem \ref{Short.time.existence} above. 
Obviously, by construction and statement (\ref{gamma}) $\Mill^c$ is a submanifold 
of $h^{4+\alpha}(\CC)$ with tangent space $T_0(\Mill^c) = \textnormal{Ker}(T_{\CC})$, which is $8$-dimensional on account of Lemma \ref{Eigenspaces.DG}. In addition, the invariant manifold $\Mill^c$ is ``exponentially attractive'' by Theorem 5.8 in \cite{Simonett.1995}. 
This means here precisely the following: 
Due to Lemmata \ref{Lemma.2.1}, \ref{Lemma.3.1} and \ref{Eigenspaces.DG} and Theorem \ref{Short.time.existence}, 
we may apply Theorem 5.8 in \cite{Simonett.1995}, and this 
theorem guarantees us, that there is some appropriate 
$\omega \in (0, \mu_1)$ - where $\mu_1$ is the smallest positive eigenvalue of $T_{\CC}$ respectively $T^s_{\CC}$ by 
decomposition (\ref{decomposition.b}) 
- a positive constant $c=c(\omega,\beta,\alpha)$ and a neighborhood $W$ of $0$ in $h^{2+\beta}(\CC)$, such that
\begin{equation}  \label{3.6}
	\parallel \pi^s(\rho(t, \rho_0)) - \gamma(\pi^c(\rho(t, \rho_0))) \parallel_{h^{4+\alpha}(\CC)} 
	\leq c \,\,\frac{e^{-\omega\, t}}{t^{1-\theta}} \,
	\parallel \pi^s(\rho_0) - \gamma(\pi^c(\rho_0)) \parallel_{h^{2+\beta}(\CC)} 
\end{equation}
for each $\rho_0 \in W$ and for $t \in (0, t^+(\rho_0))$, 
as long as there holds $\pi^c\rho(t, \rho_0) \in U$, and where we set $\theta:= \frac{(2 + \beta - \alpha)}{4}$. 
Here, $[(t,\rho_0) \mapsto \rho(t,\rho_0)]$, 
$t\in [0,t^+(\rho_0))$, denotes the unique classical and maximal solution of initial value problem (\ref{initial.value.problem}) from Theorem \ref{Short.time.existence} above.                      \\
ii) According to estimate (\ref{3.6}) we know that the invariant manifold $\Mill^c$ contains all smooth equilibria $\rho$ of equation (\ref{normal.speed.3}), which are contained in a sufficiently small neighborhood of $0$ in $h^{2+\beta}(\CC)$.
Now, again following closely the proof of Theorem 1.2 in 
\cite{Escher.Mayer.Simonett.1998} respectively of Proposition 6.4 in \cite{Escher.Simonett.1998}, we show that - at least locally about $0$ - $\Mill^c$ consists only of equilibria of equation (\ref{normal.speed.3}), 
more precisely of smooth distance functions $\rho$, whose induced maps $\theta_{\rho}(x):= \exp_x(\rho(x) \,\nu_{\CC}(x))$, $x\in \CC$, from line (\ref{theta.exp.rho}) yield $C^{\infty}$-diffeomorphisms between $\CC$ and compact tori, which are congruent to the Clifford-torus $\CC$ in $\sphere^3$ modulo the action of $\textnormal{M\"ob}(\sphere^3)$. 
Following the notation in the proof of Theorem 1.2 in \cite{Escher.Mayer.Simonett.1998}, we call 
the set of these special equilibria of equation \eqref{normal.speed.3} ``$\Mill$''.  
To this end we firstly recall from 
Remark \ref{conformal.group}, that $\textnormal{M\"ob}(\sphere^3)\cong \textnormal{SO}^{+}(1,4)$ is a Lie-group with $10$-dimensional Lie-algebra, which is the direct sum of the vector spaces $\xi$ and $\Omega$ from Definition \ref{parallel.vectorfields}. Now 
we choose a system of $10$ globally defined and linearly independent conformal vector fields 
$\{\vec v_k\}_{k=1,\ldots,10} \subset \xi \oplus \Omega 
\subset \Gamma(T\sphere^3)$. 
For any tuple $z=(z_1,\ldots,z_{10}) \in B_1^{10}(0)$ the linear combination $V_z:=\sum_{k=1}^{10} z_k \,\vec v_k \in \Gamma(T\sphere^3)$ is smooth and generates a smooth 
family of conformal transformations $T_z(t) \in \textnormal{M\"ob}(\sphere^3)$, 
$t \in \rel$, namely in terms of the flow $\Psi_z:\sphere^3 \times \rel \longrightarrow \sphere^3$, which is 
generated by the flow lines of the initial value problem:
\begin{equation}  \label{generating.Moeb.group}
	\partial_ty(t)= V_z(y(t)), \qquad y(0)=y_0 \in \sphere^3,
\end{equation} 
for an unknown smooth function 
$y:\rel \longrightarrow \sphere^3$, and then 
setting: $T_z(t):=\Psi_z(\,\cdot\,,t)$, for every fixed 
$z \in B_1^{10}(0)$. 
Since the vector fields $V_z=\sum_{k=1}^{10} z_k \,\vec v_k$ are smooth sections of $T\sphere^3$ and also  
depend smoothly on their $10$ parameters $(z_1,\ldots,z_{10}) \in B_1^{10}(0)$, we infer from the smooth dependence of solutions to ordinary initial value problems on additional 
real parameters, which the right hand side of the underlying differential equation (\ref{generating.Moeb.group}) 
depends on, that the map 
\begin{equation}  \label{T.t}
	T_{(\,\cdot\,)}(1): B_1^{10}(0) \longrightarrow \textnormal{M\"ob}(\sphere^3)
\end{equation}    
is $C^{\infty}$-smooth. Now, the images $T_z(t)(\CC)$ are compact tori in $\sphere^3$ being conformally equivalent to $\CC$, for any $z\in B_1^{10}(0)$ and for any $t\in \rel$. 
In particular, the tori $T_z(t)(\CC)$ are Willmore tori, i.e. any immersion $f_{z,t}:\Sigma \longrightarrow \sphere^3$ parametrizing $T_z(t)(\CC)$ is a critical point of $\Will$. 
Now, we choose some small $\varepsilon >0$, consider the tori 
$\CC_z:=T_z(1)(\CC)$ for any 
$z\in B_{\varepsilon}^{10}(0)$ and obtain via  
Fermi coordinates (\ref{coordinate.functions}) 
respectively via formula (\ref{rho.Gamma}) a unique 
smooth function $\rho_{z} \equiv \rho_{\CC_z}$, 
which measures the pointwise geodesic distance between 
general points $x\in \CC$ and the torus $\CC_z$. Hence, the function $[x\mapsto X(x,\rho_{z}(x))]$ parametrizes $\CC_z$ diffeomorphically:
\begin{equation}  \label{parametrizing.C.z}
	X(\,\cdot\,,\rho_z(\,\cdot\,)): \CC \stackrel{\cong}\longrightarrow\CC_z, 
\end{equation} 
as a graph over $\CC$ via the exponential map, for any fixed $z\in B_{\varepsilon}^{10}(0)$, provided $\varepsilon>0$ is sufficiently small. On account of the smoothness of the map $T_{(\,\cdot\,)}(1)$ in (\ref{T.t}), and on account of formula (\ref{rho.Gamma}) combined with the smoothness of the Fermi coordinate functions $S$ and $\Lambda$ in \eqref{coordinate.functions}, we easily infer 
also the smoothness of the non-linear operator 
\begin{equation} \label{rho.z}
	\rho_{(\,\cdot\,)}: B_{\varepsilon}^{10}(0) \longrightarrow  h^{4+\alpha}(\CC),
\end{equation}
assigning to $z\in B_{\varepsilon}^{10}(0)$ the unique 
smooth distance function $\rho_{z}$, 
that we have just obtained via formula (\ref{rho.Gamma}).   
Since the tori $\CC_{z}$ are Willmore tori, every distance function $\rho_z$ is an equilibrium of the corresponding evolution equation (\ref{normal.speed.3}), for $z\in B_{\varepsilon}^{10}(0)$, which implies that 
\begin{equation}  \label{rho.in.Mill} 
	\rho_z \in \Mill^c  \qquad \forall \,z\in  B_{\varepsilon}^{10}(0)
\end{equation} 
according to estimate (\ref{3.6}), provided $\varepsilon>0$ 
is sufficiently small, where we have used that 
$T_{0}(1)=\textnormal{id}_{\sphere^3}$ implies 
$\rho_{0}=0$ in $h^{4+\alpha}(\CC)$.
As in the proof of Theorem 1.2 in \cite{Escher.Mayer.Simonett.1998},
we consider now the composition 
\begin{equation}  \label{F}
	F:= \pi^c \circ \rho_{(\,\cdot\,)}: B_{\varepsilon}^{10}(0) \longrightarrow \textnormal{ker}(B_{\CC}),
\end{equation}    
for $\varepsilon>0$ as small as in (\ref{rho.in.Mill}).  
We note that statement (\ref{rho.z}) implies, that the map $F$ is a smooth map between finite dimensional flat manifolds, and that $\rho_{0}=0$ implies that $F(0)=0$ in $\textnormal{ker}(B_{\CC})$. 
Now, the vector space $\Omega \oplus \xi \subset T(\sphere^3)$ of all conformal vector fields on $\sphere^3$ is $10$-dimensional, whereas the vector space $\Gamma_{\CC} \equiv \xi^N \oplus \Omega^N \subset \Gamma(N\CC)$ of ``normal conformal directions along $\CC$'' is only $8$-dimensional, and the proofs of Lemmata 3.4 and 3.5 in \cite{Weiner} show, that the kernel of the linear projection $\Omega \longrightarrow \Omega^N$, sending $V \mapsto V^N$, is two-dimensional, 
whereas the projection of $\xi$ onto $\xi^N$ is isomorphic. 
Hence, we may assume without loss of generality that 
the $10$ basis vectors $\vec v_k$ of $\Omega \oplus \xi$
are chosen in such a way that $(\vec v_k)^N \equiv 0$, for $k=9,10$, i.e. such that $\{(\vec v_k)^N\}_{k=1,\ldots,8}$
is a basis of the vector space $\Gamma_{\CC}=\Omega^N \oplus \xi^N$ of ``normal conformal directions'' along $\CC$.  
Moreover, we infer from the definition of the parametrizations 
$X(\,\cdot\,,\rho_z(\,\cdot\,))$ of the tori $\CC_z=T_z(1)(\CC)$ in (\ref{parametrizing.C.z}) via the construction of the map $T_{(\,\cdot\,)}(1)$ in \eqref{generating.Moeb.group}--\eqref{T.t} and 
diffeomorphism \eqref{tubular.map.2}, that there holds
for each of the first $8$ conformal directions $\vec v_k$:
\begin{eqnarray}   \label{v.k}
(\vec v_k)^N(x) = D_{z_k}X(x,\rho_z(x))\lfloor_{z=0} 
\equiv D_{z_k}\big{(}\exp_x(\rho_z(x)\,\nu_{\CC}(x))\big{)}
\lfloor_{z=0}  =
\nonumber \\
= D_v(\exp_{x})(0) \,
\Big{(} D_{z_k} \big{(} \rho_z(x)\,\nu_{\CC}(x) \big{)}\lfloor_{z=0} \Big{)}    
= D_{z_k} \big{(} \rho_z(x) \,\nu_{\CC}(x) \big{)}\lfloor_{z=0}   \\       
= D_{z_k} \rho_z(x)\lfloor_{z=0}\, \nu_{\CC}(x)    \qquad \forall \,x \in \CC. \nonumber
\end{eqnarray} 
Using the isomorphism (\ref{eigenspaces.DG}) between the vector space $\textnormal{ker}(T_{\CC})$ and the 
$8$-dimensional vector space 
$\Gamma_{\CC} \equiv \xi^N \oplus \Omega^N \subset \Gamma(N\CC)$, we obtain from the chosen basis vectors $\{(\vec v_k)^N\}_{k=1,\ldots,8}$ of $\Gamma_{\CC}$ 
unique coordinate functions $\{v_k\}_{k=1,\ldots,8}$, 
which form a basis of the vector space $\textnormal{ker}(T_{\CC})$ and which satisfy by equation (\ref{v.k}):
$$
v_k(x)\, \nu_{\CC}(x) = (\vec v_k)^N(x) =  
D_{z_k} \rho_z(x)\lfloor_{z=0}\, \nu_{\CC}(x), \qquad \forall \,x \in \CC,
$$ 
and therefore:
$$ 
D_{z_k} \rho_z\lfloor_{z=0} = v_k  \in \textnormal{ker}(T_{\CC}), \quad
\textnormal{for} \,\, k=1,\ldots,8.
$$ 
On account of the definition of $F$ in (\ref{F}) and the chain rule, the partial derivative of $F$ in $z=0$ in direction of the coordinate $z_k$ turns out to be:
$$ 
D_{z_k}F(0) = \pi^c\big{(}D_{z_k} \rho_z\lfloor_{z=0}\big{)} = \pi^c(v_k) = v_k, \quad \textnormal{for} \,\, k=1,\ldots,8,
$$
showing that the entire differential 
$$
DF(0) :\rel^{10} \longrightarrow \textnormal{ker}(T_{\CC})
$$ 
is an epimorphism. Since we also know that there holds $F(0)=0$, we can now infer from the classical ``open mapping theorem'' for $C^1$-maps between finite dimensional
vector spaces, that there is some small open ball $B_{\delta}(0)\subset U$ about $0$ in $\textnormal{ker}(T_{\CC})$, depending 
on the size of $\varepsilon$ in 
(\ref{rho.in.Mill}) and (\ref{F}), 
which satisfies $B_{\delta}(0) \subset F(B_{\varepsilon}^{10}(0))$. 
By definition of the map $F$ in (\ref{F}) this means, that the projection $\pi^c$ of $h^{4+\alpha}(\CC)$ onto $\textnormal{ker}(T_{\CC})$ restricted 
to the set of distance functions $\{\rho_z \,|\,z\in B_{\varepsilon}^{10}(0)\}$ 
covers the open ball $B_{\delta}(0)$ in $\textnormal{ker}(T_{\CC})$.    
Now, on account of statement (\ref{rho.in.Mill}) we know that 
$\{\rho_z \,|\,z\in B_{\varepsilon}^{10}(0)\}$ is contained in the $8$-dimensional manifold $\Mill^c$, which is the graph of the function $\gamma \in C^m(U,h_s^{4+\alpha}(\CC))$ in (\ref{gamma}) over the neighborhood $U$ of $0$ in $\textnormal{Ker}(T_{\CC})$, provided $\varepsilon >0$ has been chosen sufficiently small. Since we also know that the manifold 
$\Mill^c=\textnormal{graph}(\gamma)$ intersects each fiber 
of the projection $\pi^c$ over the $8$-dimensional tangent space $T_0(\textnormal{graph}(\gamma))
=\textnormal{ker}(T_{\CC})$ only once and that the set $\{\rho_z \,|\,z\in B_{\varepsilon}^{10}(0)\}$ lies over the ball $B_{\delta}(0)$ in $\textnormal{ker}(T_{\CC})$ w.r.t. 
$\pi^c$, we conclude that $\textnormal{graph}(\gamma\lfloor_{B_{\delta}(0)})$ 
is contained in the set of smooth equilibria of equation 
(\ref{normal.speed.3}) of type $\Mill$. 
Hence, at least locally about the zero-function the center manifold $\Mill^c=\textnormal{graph}(\gamma)$ only consists of smooth equilibria of equation (\ref{normal.speed.3}) of type $\Mill$.   \\
(iii) As in the proof of Theorem 1.2 in \cite{Escher.Mayer.Simonett.1998},
we can infer from the result of step (ii) that the locally ``reduced flow'' of equation (\ref{normal.speed.3}) on $\Mill^c=\textnormal{graph}(\gamma)$
- which is determined by flow lines $\{z_t\}$ of class 
$C^0([0,t^+),\textnormal{ker}(T_{\CC})) \cap 
C^{\infty}((0,t^+),\textnormal{ker}(T_{\CC}))$ of the ``reduced equation'' 
\begin{equation}  \label{reduced}  
	\partial_t z_t + T_{\CC}^c(z_t) = \pi^cg(z_t,\gamma(z_t)), \qquad z_0 \in 
	B_{\delta}(0)\subset \textnormal{ker}(T_{\CC})     
\end{equation} 
according to the decomposition in (\ref{decomposition.2}) - 
consists of equilibria only, i.e. the locally ``reduced flow''
of equation (\ref{normal.speed.3}) does not move at all 
in a sufficiently small neighborhood about $0$ in $\textnormal{ker}(T_{\CC})$.  
In particular, the zero-function is a stable equilibrium for the ``reduced flow''. Hence, by Proposition 3.2 respectively Theorem 3.3 in \cite{Simonett.1994} also the point $(0,\gamma(0))=(0,0)$ is a stable equilibrium for the original evolution equation (\ref{normal.speed.3}) in $h^{2+\beta}(\CC)$. 
This means precisely, that there exists for every neighborhood $W_1$ of $0$ in $h^{2+\beta}(\CC)$ another neighborhood $W_2$ of $0$ in $h^{2+\beta}(\CC)$, such that a solution of evolution equation (\ref{normal.speed.3}) exists globally and stays within $W_1$, provided its initial value $\rho_0$ is contained in $W_2$. Combining this with statement (\ref{3.6}), we obtain even more precise information: There is a neighborhood $W$ of $0$ in $h^{2+\beta}(\CC)$,
depending on the size of the neighborhood $U$ of $0$ in 
$\textnormal{ker}(T_{\CC})$ from line (\ref{gamma}), 
such that any flow line of evolution equation (\ref{normal.speed.3}), which starts moving in some initial function $\rho_0 \in W$, exists globally and approaches the center manifold $\Mill^c=\textnormal{graph}(\gamma)$ asymptotically in the $h^{4+\alpha}(\CC)$-norm for all $t > 0$, according to estimate (\ref{3.6}). See here also some technical explanations in the proof of Theorem 6.5 in \cite{Escher.Simonett.1998}. \\
(iv) Again following the last step of the proof of Theorem 1.2 in \cite{Escher.Mayer.Simonett.1998}, i.e. using the reasoning of the proof of Theorem 6.5 in \cite{Escher.Simonett.1998}, which means adjusting the proof of Proposition 9.2.4 in \cite{Lunardi} and using Proposition 5.4(a) in \cite{Simonett.1995}, 
in combination with the bootstrap-technique of Proposition 6.6 in \cite{Escher.Simonett.1998}, one can draw the following conclusions from the results of steps (i)--(iii):
For any fixed $k \in \nat$ and for some appropriately chosen $\omega \in (0, \mu_1)$ - where $\mu_1$ denotes the smallest positive eigenvalue of $A_{\CC}$ respectively of $A^s_{\CC}$ by equation (\ref{decomposition.b}) - 
there exists a neighborhood $W = W(k,\omega)$ of $0$ in $h^{2+\beta}(\CC)$ with the following properties: 
Given an initial function $\rho_0 \in W$, the unique maximal classical solution $\{\rho(t, \rho_0)\}_{t\in [0,t^+(\rho_0))}$ of initial value problem \eqref{initial.value.problem} exists globally, and there exist a constant $c=c(k,\omega)>0$ and a unique function 
$z_0 = z_0(\rho_0) \in B_{\delta}(0)\subset \textnormal{ker}(T_{\CC})$, such that
\begin{eqnarray} 
	\parallel \big{(}\pi^c(\rho(t, \rho_0)), \pi^s(\rho(t, \rho_0))\big{)} 
	- (z_0, \gamma(z_0)) \parallel_{C^k(\CC)} 
	\leq   \label{exponential.convergence}   \\ 
	\leq  c\, e^{-\omega\, t} \, 
	\parallel \pi^s(\rho_0) - \gamma(\pi^c(\rho_0)) \parallel_{h^{2+\beta}(\CC)}          \nonumber
\end{eqnarray} 
holds for all $t \geq 1$. According to step (ii), we know that for $z_0 \in B_{\delta}(0)$ the pair $(z_0, \gamma(z_0)) \in \Mill^c$ has to be contained in the set of equilibria of equation (\ref{normal.speed.3})
of type $\Mill$ as well. Hence, $(z_0, \gamma(z_0))$ is a smooth distance function on $\CC$ whose induced map 
$\theta_{(z_0, \gamma(z_0))}(x):= 
\exp_x((z_0, \gamma(z_0))(x) \,\nu_{\CC}(x))$, 
$x\in \CC$, from line (\ref{theta.exp.rho}) yields a  
$C^{\infty}$-smooth diffeomorphism between $\CC$ and some compact torus, which is conformally equivalent to the Clifford-torus $\CC$ in $\sphere^3$. 
Thus statement (\ref{exponential.convergence}) guarantees that: having fixed some $k \in \nat$ and some appropriate $\omega \in (0, \mu_1)$, for any initial distance function 
$\rho_0$ taken from a sufficiently small neighborhood $W= W(k,\omega)$ of $0$ in $h^{2+\beta}(\CC)$ the unique maximal solution $\{\rho(\,\cdot\,,\rho_0)\}$ of equation (\ref{normal.speed.3}) exists globally and converges fully to a smooth distance function $(z_0, \gamma(z_0))$, which yields - via formula (\ref{theta.exp.rho}) -
a smooth diffeomorphic parametrization $\theta_{(z_0, \gamma(z_0))}:\CC \stackrel{\cong}\longrightarrow \sphere^3$ of a compact torus in $\sphere^3$, which is conformally equivalent to the Clifford-torus, and this convergence is at an exponential rate w.r.t. the $C^k(\CC)$-norm as $t \to \infty$. 
Finally, we remark that we require the initial immersion $F_1:\Sigma \longrightarrow \sphere^3$ 
to be $C^{\infty}$-smooth. On account of Theorem 1  
in \cite{Jakob_Moebius_2016} this implies, that there is a
unique and maximal flow line 
$\{\PP(t,0,F_1)\}_{t\in [0,T_{\textnormal{max}})}$ 
of the MIWF, starting in $F_1$ at time $t=0$, which is additionally $C^{\infty}$-smooth on $\Sigma \times [0,T_{\textnormal{max}})$, 
and moreover the proof of Theorem 1 in \cite{Jakob_Moebius_2016} shows, that for any given $C^{\infty}$-smooth solution $\{f_t\}_{t\in [0,T)}$ of the ``relaxed MIWF-equation'' (\ref{Moebius.flow.2}) on $\Sigma \times [0,T)$, with $f_0=F_1$ and with $T>0$ arbitrarily fixed, yields the unique smooth flow line $\{\PP(\,\cdot\,,0,F_1)\}_{t\in [0,T)}$ of the 
original MIWF, starting in $F_1$ at time $t=0$, 
only by means of reparametrization 
with a $C^{\infty}$-smooth family of $C^{\infty}$-smooth diffeomorphisms from $\Sigma$ onto itself. See here 
also the second part of Theorem \ref{Frechensbergo} above. 
Hence, on account of the ``correspondence'' 
between $C^{\infty}$-smooth flow lines of evolution equation (\ref{Moebius.flow.2}) and $C^{\infty}$-smooth flow lines of 
evolution equation (\ref{normal.speed.3}) - as 
explained in formulae (\ref{Gamma})--(\ref{normal.speed.3}) of 
Section \ref{Prep.proofs} - the above results prove 
the assertion of this theorem.                    
\qed   \\\\
\noindent
\underline{Proof of Theorem \ref{main.result.1}}\\\\	
\noindent 
On account of the assumptions of the theorem  
there is some smooth family of smooth diffeomorphisms $\Psi_t:\Sigma \stackrel{\cong}
\longrightarrow \Sigma$, $t>0$, such that the reparametrized flow line $\{\PP(t,0,F_0)\circ \Psi_t\}$ of the MIWF in $\sphere^3$ converges smoothly and fully to a diffeomorphism $F^*:\Sigma \stackrel{\cong}\longrightarrow M(\CC)$, where $M\in \textnormal{M\"ob}(\sphere^3)$ denotes some appropriate 
M\"obius-transformation and $\CC$ the Clifford-torus in $\sphere^3$. 
Now we choose some $k \in \nat$ and $\beta \in (0,1)$, and 
we obtain from Theorem \ref{Center.manifold} of this article 
some small neighborhood $W=W(\Sigma,F^*,k)$ about the limit immersion $F^*$ in $h^{2+\beta}(\Sigma,\rel^4)$, such that for every $C^{\infty}$-smooth initial immersion 
$F_1:\Sigma \longrightarrow \sphere^3$ being contained 
in $W$ the unique smooth flow line 
$\{\PP(t,0,F_1)\}_{t \geq 0}$ of the MIWF exists globally and converges - up to smooth reparametrization - fully and exponentially fast in $C^k(\Sigma,\rel^4)$ to a parametrization $E^*_{F_1}:\Sigma \stackrel{\cong}\longrightarrow M_1(\CC)$ of a compact torus, which is congruent to the Clifford-torus $\CC$ via some 
conformal transformation $M_1\in \textnormal{M\"ob}(\sphere^3)$. 
Now using the fact that the flow line 
$\{\PP(t,0,F_0)\circ \Psi_t\}$ of the MIWF in 
$\sphere^3$ converges smoothly and fully to a parametrization $F^*:\Sigma \stackrel{\cong}
\longrightarrow M(\CC)$ of the ``Clifford-torus'' $M(\CC)$, we can choose some large but finite time $T=T(F_0,k,F^*,\beta)>>1$, such that the immersion $\PP(T,0,F_0)\circ \Psi_T$ is contained in the specified neighborhood $W$ of the limit immersion $F^*$ in $h^{2+\beta}(\Sigma,\rel^4)$. 
Now, we recall that the MIWF (\ref{Moebius.flow}) 
is conformally invariant and that the stereographic $\stereo$ projection from $\sphere^3\setminus \{(0,0,0,1)\}$ into $\rel^3$ is a conformal diffeomorphism. 
Moreover, on account of the 
compactness of $\Sigma$ and on account of the conformal invariance of the MIWF we may assume, that the image of the initial immersion $F_0:\Sigma \longrightarrow \sphere^3$ does not contain the north pole $(0,0,0,1)$ of $\sphere^3$.  
Therefore, the result of Theorem \ref{Frechensbergo}
can be transported - at least to some satisfactory 
extend - in our situation from $\rel^3$ to 
$\sphere^3$ by means of stereographic projection $\stereo:\sphere^3 \setminus \{(0,0,0,1)\}
\longrightarrow \rel^3$ and its inverse conformal diffeomorphism $\stereo^{-1}$. To this end, we firstly 
see that the requirements of Theorem \ref{Frechensbergo}
are trivially satisfied here for the initial immersion
$\tilde F_0:=\stereo \circ F_0:
\Sigma \longrightarrow \rel^3$ and for any final 
time $T>0$. Hence, the first part of Theorem \ref{Frechensbergo} guarantees, that there is for 
$\tilde F_0$ and for any fixed $\gamma \in (0,1)$ and $T>0$ 
some small $\rho=\rho(\Sigma,\tilde F_0,T,\gamma)>0$, 
such that for every immersion 
$\tilde F \in C^{4,\gamma}(\Sigma,\rel^3)$ with 
$\parallel \tilde F - \tilde F_0 \parallel_{C^{4,\gamma}(\Sigma,\rel^3)}<\rho$ there is a unique, classical solution 
$\{\PP^*(t,0,\tilde F)\}_{t \in [0,T]}$ of the 
``DeTurck modification'' (\ref{de_Turck_equation_2}) of 
the MIWF-equation (\ref{Moebius.flow}) in the parabolic H\"older space 
$X_T\equiv X_{T,\gamma}=C^{4+\gamma,1+\frac{\gamma}{4}}
(\Sigma \times [0,T],\rel^3)$, starting to move in the immersion $\tilde F$ at time $t=0$, and such
that this unique solution
$\{\PP^*(t,0,\tilde F)\}_{t \in [0,T]}$ of 
equation (\ref{de_Turck_equation_2}) in 
$X_{T,\gamma}$ depends in a $C^1$-fashion on its initial immersion $\tilde F$, in the sense of statement \eqref{solution.operator.2} in Theorem \ref{Frechensbergo}.    
Now, combining this information also with the second part of Theorem \ref{Frechensbergo} and applying again inverse stereographic projection $\stereo^{-1}$, we can therefore 
at least infer, that in our situation there is for any $\varepsilon>0$ some sufficiently small $r=r(\Sigma,F_0,T,\varepsilon,\gamma)>0$, 
such that for every $C^{\infty}$-smooth immersion 
$F:\Sigma \longrightarrow \sphere^3$ 
with $\parallel F-F_0 \parallel_{C^{4,\gamma}(\Sigma,\rel^4)}<r$
the unique, maximal smooth flow line 
$\{\PP(\,\cdot \,,0,F)\}$ of the MIWF in $\sphere^3$
exists at least on $\Sigma \times [0,T]$, and such that the smooth flow lines $\{\PP(t,0,F_0)\}_{t\in [0,T]}$ and 
$\{\PP(t,0,F)\}_{t\in [0,T]}$ of the MIWF in 
$\sphere^3$ can be reparametrized by smooth families of smooth diffeomorphisms 
$\phi^{F_0}_t:\Sigma \longrightarrow \Sigma$ and
$\phi^F_t:\Sigma \longrightarrow \Sigma$ in such a way, 
that the reparametrized flow lines  
$\{\PP(t,0,F_0) \circ \phi^{F_0}_t\}_{t \in [0,T]}$ and
$\{\PP(t,0,F) \circ \phi^F_t\}_{t \in [0,T]}$ 
satisfy:
\begin{equation}  \label{close.flow.lines}
	\parallel \PP(t,0,F) \circ \phi^F_t - 
	\PP(t,0,F_0)\circ \phi^{F_0}_t
	\parallel_{C^{4,\gamma}(\Sigma,\rel^4)} < \varepsilon,
\end{equation} 
for every $t\in [0,T]$. Now, we had chosen $T=T(F_0,k,F^*,\beta)$ that large, such that the immersion $\PP(T,0,F_0)\circ \Psi_T$ is contained in the neighborhood $W$ of the limit immersion $F^*$ in $h^{2+\beta}(\Sigma,\rel^4)$. Hence, by estimate (\ref{close.flow.lines}) the diffeomorphism
$\Theta^{F_0}_{T}:=(\phi^{F_0}_T)^{-1}\circ \Psi_T:
\Sigma \stackrel{\cong}\longrightarrow \Sigma$ has the 
property that both immersions 
$\PP(T,0,F_0)\circ \phi^{F_0}_T \circ \Theta^{F_0}_{T} = \PP(T,0,F_0)\circ \Psi_T$
and $\PP(T,0,F) \circ \phi^F_T \circ \Theta^{F_0}_{T}$ 
are contained in the neighborhood $W$ of the limit immersion 
$F^*$ in $h^{2+\beta}(\Sigma,\rel^4)$, at time $t=T$, provided $\varepsilon>0$ has been chosen sufficiently small in estimate (\ref{close.flow.lines}), the initial smooth immersion $F:\Sigma \longrightarrow \sphere^3$ is contained in the open ball $B_r(F_0)$ about $F_0$ in $C^{4,\gamma}(\Sigma,\rel^4)$, and provided $r=r(\Sigma,F_0,T,\varepsilon,\gamma)
=r(\Sigma,F_0,F^*,k,\gamma,\beta,W)>0$
had also been chosen sufficiently small in estimate (\ref{close.flow.lines}). Here, we have also used the obvious embedding $C^{4,\gamma}(\Sigma,\rel^4) 
\hookrightarrow h^{2+\beta}(\Sigma,\rel^4)$.
Recalling that the neighbourhood $W$ of $F^*$ 
only depends on $\Sigma$, $F^*$ and $k$, we 
finally see that the radius $r=r(\Sigma,F_0,F^*,k,\gamma,\beta,W)$ 
only depends on $\Sigma,F_0$ and $F^*$ and on the 
parameters $k,\gamma,\beta$, i.e. $r=r(\Sigma,F_0,F^*,k,\gamma,\beta)$,  
and we conclude: if this number $r$ is chosen  
sufficiently small, then for any smooth immersion 
$F:\Sigma \longrightarrow \sphere^3$ being contained 
in the open ball $B_r(F_0)$, the reparametrized immersion 
$\PP(T,0,F) \circ \phi^F_T \circ \Theta^{F_0}_{T}$ 
is an element of the prescribed neighborhood $W$ about the limit immersion $F^*$ in $h^{2+\beta}(\Sigma,\rel^4)$.
Now since we know already from above, that the entire 
reparametrized flow line 
$\{\PP(t,0,F) \circ \phi^F_t\}_{t \in [0,T]}$ of 
the MIWF in $\sphere^3$ is of class 
$C^{\infty}(\Sigma \times [0,T],\rel^4)$, we  
especially conclude that the immersion 
$\PP(T,0,F) \circ \phi^F_T \circ \Theta^{F_0}_{T}$ 
is $C^{\infty}$-smooth on $\Sigma$. 
We can therefore choose the above initial immersion $F_1$ 
from the statement of Theorem \ref{Center.manifold}
of this article to be $F_1:= \PP(T,0,F) \circ \phi^F_T \circ \Theta^{F_0}_{T}$ and infer from Theorem \ref{Center.manifold}, that the unique flow line 
$\{\PP(t,0,\PP(T,0,F) \circ \phi^F_T \circ \Theta^{F_0}_{T})\}_{t \geq 0}$ of the MIWF in $\sphere^3$, starting to move in the smooth immersion 
$\PP(T,0,F) \circ \phi^F_T \circ \Theta^{F_0}_{T}$ 
at time $t=0$, converges - up to smooth reparametrization - fully and exponentially fast in $C^k(\Sigma,\rel^4)$ to a diffeomorphic parametrization 
$E^*_{F}:\Sigma \stackrel{\cong}\longrightarrow M^F(\CC)$ 
of a compact torus, which is conformally equivalent to the Clifford-torus $\CC$ in $\sphere^3$, provided $F:\Sigma \longrightarrow \sphere^3$ is a smooth immersion being contained in the ball $B_r(F_0)$ about $F_0$ in $C^{4,\gamma}(\Sigma,\rel^4)$ and $r=r(\Sigma,F_0,F^*,k,\gamma,\beta)>0$ is sufficiently small. 
On account of the invariance of the MIWF w.r.t. 
time-independent smooth reparametrizations, this means that any flow line $\{\PP(t,0,F)\}_{t \geq 0}$ of the MIWF 
in $\sphere^3$, which starts moving in some arbitrarily chosen smooth immersion 
$F:\Sigma \longrightarrow \sphere^3$ belonging to the $C^{4,\gamma}$-ball $B_r(F_0)$ about $F_0$, converges - 
up to smooth reparametrization - fully and exponentially 
fast in the $C^k(\Sigma,\rel^4)$-norm to a diffeomorphic parametrization $E^*_{F}:\Sigma \stackrel{\cong}\longrightarrow M^F(\CC)$ of a conformally transformed Clifford-torus $M^F(\CC)$ in $\sphere^3$, provided $r=r(\Sigma,F_0,F^*,k,\gamma,\beta)>0$ was chosen sufficiently small.
\qed     \\\\
\noindent
\underline{Proof of Theorem \ref{Convergence.to.local.minimizer}} \\\\
We shall combine here the methods of the proof of Lemma 4.1 in \cite{Chill.Schatz.2009} with the methods of 
the proof of Theorem 1.2 in \cite{Dall.Acqua.Spener.2016},
until we will have to involve Rivi\`ere's and Bernard's
technique from \cite{Riviere.2008} and \cite{Bernard.2016}, 
in order to accomplish the entire argument.
First of all, we require by assumption that the initial immersion $f_0:\Sigma \longrightarrow \rel^3$ is smooth and satisfies:   
\begin{equation}   \label{smallness} 
\parallel f_0 - F^* \parallel_{C^{k,\alpha}(\Sigma)} <\varepsilon 
\end{equation} 
for some sufficiently small $\varepsilon \in (0,\varepsilon_0)$, 
to be determined later. As in the proof of Lemma 4.1 in \cite{Chill.Schatz.2009}, we conclude from condition (\ref{smallness}) that we can represent the immersion $f_0$ as a graph over the Willmore immersion $F^*$, which means precisely on account of Theorem 5.1 
in \cite{Skorzinski.2015}: There is a smooth section 
$N_0$ of the normal bundle of $F^*$ and a smooth diffeomorphism 
$\Phi_0:\Sigma \stackrel{\cong}\longrightarrow \Sigma$, 
such that there holds: 
\begin{equation}  \label{normal.representation} 
f_0 \circ \Phi_0 = F^* + N_0 \quad \textnormal{on} \,\, \Sigma.
\end{equation} 
Furthermore, we infer from condition (\ref{smallness}) and 
again from Theorem 5.1 in \cite{Skorzinski.2015} combined  
with Lemma 3.1 in \cite{Farkas.Garay.2000}
\footnote{See here also Section 2 in \cite{Ebin.Marsden.1970}.}, 
that there is some continuous and monotonically increasing 
function $C^o:[0,\varepsilon_0]\longrightarrow \rel_+$,
with $C^{o}(0)=0$, depending also on the immersion $F^*$ 
and on $k$, such that the estimate
\begin{equation}   \label{smallness.2} 
	\parallel N_0 \parallel_{C^{k,\alpha}(\Sigma)} 
	= \parallel f_0 \circ \Phi_0 - F^* \parallel_{C^{k,\alpha}(\Sigma)}<C^o(\varepsilon) 
\end{equation} 
holds for the $\varepsilon$ from line \eqref{smallness}. 
Now we consider the modified Cauchy-problem
\begin{equation}  \label{modified.Cauchy} 
	\partial_t^{\perp_{\tilde f_t}}(\tilde f_t) = 
	- \frac{1}{|A^0_{\tilde f_t}|^4} \,
	\nabla_{L^2}\Will(\tilde f_t), \quad  
	\tilde f_0 \stackrel{!}=f_0 \circ \Phi_0 \quad \textnormal{on} \,\,\Sigma
\end{equation} 
of the MIWF, which is solved by any smooth reparametrization 
$\{\PP(t,0,f_0) \circ \Phi_t\}_{t\geq 0}$ 
of the smooth flow line $\{\PP(\,\cdot\,,0,f_0)\}$
of the MIWF, starting in $f_0$.
On account of condition (\ref{smallness}), equation  (\ref{normal.representation}) and again
Theorem 5.1 in \cite{Skorzinski.2015}, we may assume 
the existence of a smooth short-time solution $\{\tilde f_t\}$
of equation (\ref{modified.Cauchy}) being of the special form: 
$\tilde f_t= F^*+ N_t$\, on $\Sigma$,  
for a family of smooth normal sections $N_t$ along $F^*$ 
and starting in the immersion $F^*+ N_0 = f_0 \circ \Phi_0$ 
at time $t=0$. As in the proof of Lemma 4.1 in \cite{Chill.Schatz.2009}, we can introduce the function 
$\phi_t:= \langle N_t,\nu_{F^*} \rangle_{\rel^3}$, 
i.e. the signed length of the normal section $N_t$,  
and compute that the family of immersions 
$\tilde f_t = F^* + N_t$ solves 
equation (\ref{modified.Cauchy}), if and only if
the family of functions $\{\phi_t\}$ solves a uniformly parabolic quasi-linear differential equation of fourth order, namely: 
\begin{eqnarray}  \label{parabolic.equation}
	\partial_t(\phi_t) + 
	\frac{1}{2} \, \frac{1}{|A^0_{F^*+ \phi_t \nu_{F^*}}|^4} \,\,g^{ij}_{F^*+ \phi_t \nu_{F^*}} \, 
	g^{kl}_{F^*+ \phi_t \nu_{F^*}} \, \nabla^{F^*}_{ijkl}(\phi_t) 
	= B(\,\cdot\,,\phi_t,\ldots,D^3_x\phi_t)        
\end{eqnarray}  	 
on $\Sigma$, for some globally defined function 
$B:\Sigma \times \rel^{1+2+4+8} \longrightarrow \rel$, which is rational in its $15$ real arguments and has smooth
coefficients, depending on the fixed immersion $F^*$ only, 
at least as long as 
\begin{equation}  \label{delta.bounded} 
	\parallel N_t \parallel_{C^2(\Sigma)} \equiv  
	\parallel \tilde f_t - F^* \parallel_{C^2(\Sigma)} 
	< \tilde \delta(F^*)
\end{equation} 
holds, for some sufficiently small chosen positive number 
$\tilde \delta(F^*)>0$, such that
\begin{equation}  \label{good.projection} 
	| P^{\perp_{\tilde f_t}}(\hat N_t) | >  
	\frac{1}{2} \, |\hat N_t| \qquad \textnormal{on} \,\, \Sigma  
\end{equation} 
holds for the projection $P^{\perp_{\tilde f_t}}(\hat N_t)$ of any smooth normal field $\hat N_t$ along $F^*$ into the 
normal bundle of $\tilde f_t$, and such that also
\begin{equation}  \label{no.umbilics}  
	\min_{\Sigma} |A^0_{F^*+ \phi_t \nu_{F^*}}|^2 >  
	\frac{1}{2}\, \min_{\Sigma} |A^0_{F^*}|^2 >0
\end{equation} 
holds at sufficiently small times $t\geq 0$. 
Moreover, we infer from condition \eqref{smallness.2}:
\begin{equation}    \label{smallness.3}
\parallel \phi_0 \parallel_{C^{k,\alpha}(\Sigma)} 
<C^o(\varepsilon)  
\end{equation}
for the smooth initial function of the solution  
$\phi_t = \langle N_t, \nu_{F^*} \rangle_{\rel^3}$
of equation \eqref{parabolic.equation},  
where $C^{o}(\varepsilon)$ denotes here the same function 
as in \eqref{smallness.2}. Now, we shall introduce the 
parabolic H\"older space
$$ 
Z_{T,\beta} := C^{4+\beta,1+\frac{\beta}{4}}
(\Sigma \times [0,T],\rel),
$$ 
for any fixed $\beta \in (\alpha,1)$ and $T>0$, 
and its open subsets 
\begin{equation} \label{U.phi.varrho.p} 
	U_{F^*,\beta,\varrho,T} :=
	\Big{\{} \, \{\varphi_t\} \in Z_{T,\beta} \,| \,
	\parallel \varphi_t \parallel_{C^2(\Sigma)}  
	<\varrho  \quad \forall\, t \in [0,T] \,\Big{\}}
\end{equation} 
with $0<\varrho < \tilde \delta(F^*)$ that small, 
such that the normal field 
$N_t := \varphi_t \,\nu_{F^*}$ satisfies  
inequality (\ref{delta.bounded}) $\forall \,t \in [0,T]$, 
for any fixed function $\varphi \in U_{F^*,\beta,\varrho,T}$. 
Now, by statements (\ref{delta.bounded})--(\ref{no.umbilics}) and (\ref{U.phi.varrho.p}) the non-linear differential operator 
\begin{eqnarray}   \label{parabolic.operator}  
	\partial_t + \frac{1}{2} \, 
	\frac{1}{|A^0_{F^*+\, (\,\cdot\,)\, \nu_{F^*}}|^4} \,\,
	g^{ij}_{F^*+ \, (\,\cdot\,)\, \nu_{F^*}} \, 
	g^{kl}_{F^*+ \, (\,\cdot\,)\, \nu_{F^*}} \, \nabla^{F^*}_{ijkl}:\\ \nonumber 
	U_{F^*,\beta,\varrho,T} \subset Z_{T,\beta}  \longrightarrow  C^{\beta,\frac{\beta}{4}}(\Sigma \times [0,T],\rel)
\end{eqnarray}
is well-defined for any fixed 
$0<\varrho < \tilde \delta(F^*)$ as in (\ref{U.phi.varrho.p}), 
and we can infer exactly as in Theorem 2 in \cite{Jakob_Moebius_2016}, that it is a $C^1$-map, 
that the highest order term of its Fr\'echet derivative in any chosen $\varphi \in U_{F^*,\beta,\varrho,T}$ is the uniformly 
parabolic linear operator 
\begin{eqnarray}  \label{leading.operator}  
	\partial_t  + \Lift_{F^*,\varphi} 
	:= \partial_t + 
	\frac{1}{2} \, \frac{1}{|A^0_{F^*+ \varphi \nu_{F^*}}|^4} \,
	g^{ij}_{F^*+ \varphi \nu_{F^*}} \, 
	g^{kl}_{F^*+ \varphi \nu_{F^*}} \, \nabla^{F^*}_{ijkl}: \nonumber\\
	Z_{T,\beta} \longrightarrow 
	C^{\beta,\frac{\beta}{4}}(\Sigma \times [0,T],\rel)
\end{eqnarray}
for any fixed $\beta \in (\alpha,1)$ and $T>0$, 
and that this linear operator satisfies all 
requirements of Proposition 2.1 in \cite{Ruben.MIWF.III}.
We can therefore argue as in the proof of that 
proposition, that for any fixed 
$\varphi \in U_{F^*,\beta,\varrho,T}$ 
the linear differential operators 
\begin{eqnarray} \label{leading.operator.t}
	\Lift_{F^*,\varphi_t} 
	:= \frac{1}{2} \, \frac{1}{|A^0_{F^*+ \varphi_t \nu_{F^*}}|^4} \,g^{ij}_{F^*+ \varphi_t \nu_{F^*}} \, 
	g^{kl}_{F^*+ \varphi_t \nu_{F^*}} \, \nabla^{F^*}_{ijkl}: \nonumber\\ 
	C^{4,\beta'}(\Sigma,\rel)  \longrightarrow  C^{0,\beta'}(\Sigma,\rel)     \qquad 
\end{eqnarray}
are $(\Wil,\frac{3}{4}\pi,4)$-elliptic for 
some appropriate constant 
$\Wil=\Wil(F^*,\varrho)>1$ in 
the terminology of \cite{Shao.Simonett.2014}, p. 228, for 
any fixed $\beta' \in (0,\beta]$ and 
for any fixed $t\in [0,T]$. We can therefore derive
from Theorem 3.3 in \cite{Shao.Simonett.2014}, 
that the linear operators $\Lift_{F^*,\varphi_t}$ 
in (\ref{leading.operator.t}) are sectorial
in the H\"older space $C^{0,\beta'}(\Sigma,\rel)$ 
with constants $\omega>0$ and $\NN>1$ depending only on   
$F^*,\parallel \varphi 
\parallel_{C^{4+\beta,1+\frac{\beta}{4}}(\Sigma \times [0,T])}$ 
and $\beta'$, in the terminology of Theorem 3.3 in \cite{Shao.Simonett.2014}, 
for any fixed $\beta' \in (0,\beta)$ and uniformly for every fixed 
$t\in [0,T]$, having chosen $T>0$ and $\beta \in (\alpha,1)$ 
already in line (\ref{U.phi.varrho.p}). See here also 
Definition 2.0.1 in \cite{Lunardi}. Choosing now 
$\beta'=\alpha$ and recalling also the smoothness of the 
initial function 
$\phi_0:= \langle N_0,\nu_{F^*} \rangle_{\rel^3}$,
we may therefore apply Theorems 8.1.1 and 8.1.3 and  
Proposition 8.2.1 in \cite{Lunardi} with Banach space 
pair $C^{4,\alpha}(\Sigma,\rel)\subset C^{0,\alpha}(\Sigma,\rel)$ and obtain the existence 
of a unique and maximal strict solution 
$\{\phi_t\}_{t\in [0,T_{\textnormal{max}})}$ of 
equation (\ref{parabolic.equation}) with values in the 
Banach space $C^{4,\alpha}(\Sigma,\rel)$, 
meeting the condition: 
\begin{equation} \label{smaller.than.varrho} 
\parallel \phi_t \parallel_{C^2(\Sigma)}  
<\varrho  \quad \forall \,t \in [0,T_{\textnormal{max}}) 
\end{equation} 
from line (\ref{U.phi.varrho.p}) and 
starting in the initial function 
$\phi_0$ at time $t=0$, provided there holds 
$C^o \,\varepsilon <\varrho$. Moreover, this solution is of class $C^{\gamma}([0,T],C^{4,\alpha}(\Sigma,\rel)) \cap 
C^{1,\gamma}([0,T],C^{0,\alpha}(\Sigma,\rel))$, 
$\forall \, \gamma \in (0,1)$ and for every 
$T \in (0,T_{\textnormal{max}})$, where $\alpha \in (0,1)$ 
had been chosen already in the initial 
assumption (\ref{smallness}). We may therefore apply 
Proposition 3 of \cite{Jakob_Moebius_2016} - 
successively for every $k\in \nat_0$ -
in order to conclude, that the maximal solution $\{\phi_t\}$ 
of equation (\ref{parabolic.equation}) also satisfies the 
following Schauder a-priori estimates:
\begin{eqnarray}    \label{Schauder.estimates}
	\parallel \phi_t  \parallel_{C^{4+l+\mu,1+\frac{l+\mu}{4}}
		(\Sigma \times [0,T])}   \leq                            \nonumber \\
	\leq C \,\Big{(}  \parallel
	B(\,\cdot\,, \phi(\,\cdot\,,t),\ldots,
	\partial_{xxx} \phi(\,\cdot\,,t))   \parallel_{C^{l+\mu,\frac{l+\mu}{4}}(\Sigma \times [0,T])}  + \nonumber  \\
	+ \parallel \phi_t  \parallel_{L^{\infty}(\Sigma \times [0,T])}
	+ \parallel \phi_0 \parallel_{C^{4+l+\mu}(\Sigma)} \Big{)},       \qquad
\end{eqnarray}
for every fixed $T\in (0,T_{\textnormal{max}})$, for 
every $l \in \nat_0$ and every $\mu \in (0,\alpha]$, 
and for some large constant $C=C(\Sigma,F^*,T,\mu,l)$.
From this result we immediately infer the 
$C^{\infty}$-smoothness of the maximal solution 
$\{\phi_t\}$ of equation (\ref{parabolic.equation}) 
and also the existence of a corresponding solution of equation 
(\ref{modified.Cauchy}) of the special form 
$\tilde f_t = F^* + N_t$, with 
$N_t=\phi_t \,\nu_{F^*}$, of class 
$C^{\infty}(\Sigma \times [0,T],\rel^3)$ for 
every $T \in [0,T_{\textnormal{max}})$, which starts 
moving in the immersion $\tilde f_0 =f_0 \circ \Phi_0$ 
at time $t=0$. As pointed out in the proofs of Lemma 4.1 in \cite{Chill.Schatz.2009} and of Theorem 1.2 in \cite{Dall.Acqua.Spener.2016}, this smooth solution 
$\tilde f_t = F^* + N_t$ of equation 
(\ref{parabolic.equation}) can be reparametrized 
by a smooth family of smooth diffeomorphisms 
$\Psi_t:\Sigma \stackrel{\cong}\longrightarrow \Sigma$, 
with \,$\Psi_0 = \textnormal{Id}_{\Sigma}$, such that: 
\begin{equation}  \label{zusammengemogili}
	\tilde f_t \circ  \Psi_t = 
	\PP(t,0,f_0 \circ \Phi_0)
	\equiv \PP(t,0,f_0)\circ \Phi_0 \qquad 
	\textnormal{for} \,\, t\in [0,T_{\textnormal{max}}), 
\end{equation} 
which will be later of great importance. 
Now, we fix some positive $\sigma <\min\{\delta,\varrho\}$ - 
where $\delta$ was determined in assumption (\ref{local.minimum})
and $\varrho$ in line (\ref{U.phi.varrho.p}) - 
such that the Lojasiewicz-Simon-gradient-inequality 
for the Willmore functional, Theorem 3.1 in \cite{Chill.Schatz.2009}, holds for every $C^4$-immersion 
$f:\Sigma \longrightarrow \rel^3$ with 
$\parallel f - F^* \parallel_{C^4(\Sigma)} \leq \sigma$, 
and we choose the $\varepsilon>0$ in (\ref{smallness}) and (\ref{smallness.3}) that small, such that we have: 
$C^{o}(\varepsilon)<\sigma$. 
As in \cite{Chill.Schatz.2009}, p. 359, or as in \cite{Dall.Acqua.Spener.2016}, p. 2190, we shall now choose a possibly smaller ``maximal'' time $T(\sigma)\in (0,T_{\textnormal{max}}]$ - depending on $\sigma$, but not on $\varepsilon$ - by means of imposing the following additional, quantitative smallness condition on the normal sections $N_t =\phi_t \, \nu_{F^*}$:
\begin{eqnarray}  \label{enclosure.T}
	\parallel \tilde f_t - F^* \parallel_{C^{k}(\Sigma,\rel^3)} 
	\equiv \parallel N_t \parallel_{C^{k}(\Sigma,\rel^3)}     
	\stackrel{!}\leq \sigma \qquad \forall \, t\in [0,T(\sigma)).\quad
\end{eqnarray}
We should note here, that condition (\ref{enclosure.T}) 
implies the inequality:
\begin{equation}  \label{enclosure.T.2}
	\parallel \phi_t \parallel_{C^{k}(\Sigma,\rel)} 
	\leq  \sigma  \qquad \forall \, t\in [0,T(\sigma)).
\end{equation}
Comparing condition (\ref{smaller.than.varrho}) with 
inequality (\ref{enclosure.T.2}) and 
recalling initial condition (\ref{smallness.2})  
and our requirements that $C^{o}(\varepsilon)<\sigma$
and $\sigma<\varrho$, we can conclude that 
$$
0<T(\sigma) \leq T_{\textnormal{max}}\leq \infty.
$$ 
In order to prove, that statement (\ref{enclosure.T}) actually 
holds for $T(\sigma)=\infty$, we shall follow the strategy 
of the proof of Theorem 1.2 in \cite{Dall.Acqua.Spener.2016}, 
pp. 2190--2191: We assume firstly, that the 
time $T(\sigma)$ was finite and that there would hold 
$T(\sigma) < T_{\textnormal{max}}$.
Now, we note that we can rewrite 
equation (\ref{parabolic.equation}) in divergence 
form, i.e. in the form:  
\begin{eqnarray}  \label{Divergence.form.equation}
	\partial_t(\phi_t) + 
	\frac{1}{2} \,\nabla^{F^*}_{ij} 
	\Big{(} \frac{1}{|A^0_{F^*+\, \phi_t\, \nu_{F^*}}|^4} \,
	g^{ij}_{F^*+ \, \phi_t\, \nu_{F^*}} \, 
	g^{kl}_{F^*+ \phi_t\, \nu_{F^*}}  \,
	\nabla^{F^*}_{kl}(\phi_t) \Big{)} = \nonumber \\
	= \tilde B(\,\cdot\,,\phi_t,\ldots,D^3_x\phi_t,D^4_x \phi_t) \quad 
\end{eqnarray}
on $\Sigma \times [0,T(\sigma)]$ for some globally defined 
function $\tilde B:\Sigma \times \rel^{1+2+4+8+16} \longrightarrow \rel$, which is rational in its $31$ real arguments and 
contains the right-hand side  
$B(\,\cdot\,,\phi_t,\ldots,D^3_x\phi_t)$ of 
equation (\ref{parabolic.equation}) as a summand. 
On account of inequality (\ref{enclosure.T.2}), on account of  
the assumption that the immersion $F^*$ is umbilic-free 
and since we assume here that $T(\sigma)$ is finite and 
smaller than $T_{\textnormal{max}}$, we know that 
there holds:
\begin{equation}  \label{tilde.B.bounded}
	\parallel 
	\tilde B(\,\cdot\,,\phi_t,\ldots,D^3_x\phi_t,D^4_x \phi_t)
	\parallel_{L^{\infty}(\Sigma)} 
	\leq C(\Sigma,F^*,\sigma)  \quad 
	\forall \,t\in [0,T(\sigma)],
\end{equation}
for some appropriate, large constant $C=C(\Sigma,F^*,\sigma)$, 
and also that the coefficients on the left hand side 
of (\ref{Divergence.form.equation}) are of class 
$C^{0,\frac{1}{2}}(\Sigma)$ and bounded in this norm 
by another large constant $C=C(\Sigma,F^*,\sigma)$, 
uniformly for every $t\in [0,T(\sigma)]$.  
Hence, we may apply Theorem 2.1 in \cite{Dong.Zhang.2015} 
to the unique smooth solution $\{\phi_t\}$ of equation (\ref{parabolic.equation}) respectively of equation  
(\ref{Divergence.form.equation}) - 
here with $m=2$, $a=\frac{1}{2}$, $g_t:= \phi_t$ 
for $t \in [0,T(\sigma)]$ and by means of a finite smooth 
atlas of the torus $\Sigma$ - 
and we infer in combination with inequalities 
\eqref{enclosure.T.2}, \eqref{tilde.B.bounded},
and again with differential equation \eqref{parabolic.equation}
the a-priori-estimate:
\begin{eqnarray}  \label{2+half.parabolic.estimate}   
	\sup_{\Sigma \times [0,T(\sigma)]} |\phi_t|  
	+ [\phi_t]_{\frac{1}{2},\frac{1}{2}+\frac{1}{8},
		\Sigma \times [0,T(\sigma)]} 
	+[D_x \phi_t]_{\frac{1}{2},\frac{1}{4}+\frac{1}{8},
		\Sigma \times [0,T(\sigma)]}         
	+[D^2_x \phi_t]_{\frac{1}{2},\frac{1}{8},
		\Sigma \times [0,T(\sigma)]}     \leq  \qquad          \\      
	\leq C(\Sigma,F^*,\sigma) \, 
	\Big{(} \parallel \phi_t \parallel_{L^2(\Sigma \times [0,T(\sigma)])} 
	+ \parallel \tilde B(\,\cdot\,,\phi_t,\ldots,D^3_x\phi_t,
	D^4_x\phi_t)
	\parallel_{L^{\infty}(\Sigma \times [0,T(\sigma)])}  \nonumber     \\                  
	+ \parallel \phi_t \parallel_{C^{2+\frac{1}{2},0}(\Sigma\times [0,T(\sigma)])}                                           
	+ \parallel \partial_t(\phi_t) 
	\parallel_{L^{\infty}(\Sigma \times [0,T(\sigma)])}  \Big{)}   \leq  C(\Sigma,F^*,\sigma,T(\sigma)),       \nonumber 
\end{eqnarray}
for some appropriate large constants $C=C(\Sigma,F^*,\sigma)$ 
and $C=C(\Sigma,F^*,\sigma,T(\sigma))$,
where we adopted Dong's and Zhang's notation in \cite{Dong.Zhang.2015}. This shows in particular, that the coefficients of the uniformly parabolic linear operator $\Lift_{F^*,\phi}$ in (\ref{leading.operator}) - here 
evaluated in $\varphi = \phi$ - and thus also of the highest order elliptic operator in (\ref{parabolic.equation}) 
are bounded in the parabolic H\"older space  
$C^{\frac{1}{2},\frac{1}{8}}(\Sigma \times [0,T(\sigma)],\rel)$ by some appropriate constant $C=C(\Sigma,F^*,\sigma,T(\sigma))$. 
We can therefore apply Proposition 2.1 in \cite{Ruben.MIWF.III} to the linear operator $\Lift_{F^*,\phi}$ in (\ref{leading.operator}) 
- here on $\Sigma \times [0,T(\sigma)]$ and with any fixed
$p\in (1,\infty)$ - and we infer from that proposition and estimate (\ref{2+half.parabolic.estimate}) 
together with conditions (\ref{smallness.3}) and (\ref{enclosure.T.2}), that the smooth solution 
$\{\phi_t\}$ of equation (\ref{parabolic.equation}) 
is bounded in the parabolic $L^p$-space 
\begin{equation}   \label{X.T}
	X_{T,p} := W^{1,p}([0,T];L^p(\Sigma,\rel)) \cap 
	L^p([0,T];W^{4,p}(\Sigma,\rel)),
\end{equation} 
for any fixed $p\in (1,\infty)$: 
\begin{eqnarray}  \label{Lp.Lp.estimate} 
	\parallel \{\phi_t\} \parallel_{X_{T(\sigma),p}} 
	\leq C(\Sigma,F^*,\sigma,T(\sigma),p) \, 
	\Big{(} \parallel B(\,\cdot\,,\phi_t,\ldots,D^3_x\phi_t) \parallel_{L^{p}([0,T(\sigma)];L^p(\Sigma,\rel))} + \nonumber        \\ 
	+ \parallel \phi_0 \parallel_{W^{4,p}(\Sigma,\rel)} \Big{)}
	\leq C^*(\Sigma,F^*,\sigma,T(\sigma),p), \qquad
\end{eqnarray}
for some approriate large constant $C^*(\Sigma,F^*,\sigma,T(\sigma),p)$. 
Moreover, as explained in Theorem B.5 in \cite{Spener} we 
can use interpolation results from \cite{Amann.2000},
in order to obtain for any $p \in (1,\infty)$ and 
$\theta \in (0,1)$ with $4(1-\theta) \not \in \nat$ the 
continuous embedding: 
\begin{eqnarray}  \label{embedding.interpol}
	X_{T,p}=W^{1,p}([0,T];L^p(\Sigma,\rel)) 
	\cap L^{p}([0,T];W^{4,p}(\Sigma,\rel))   \\
	\hookrightarrow  
	\big{(}W^{1,p}([0,T];L^p(\Sigma,\rel)), 
	L^{p}([0,T];W^{4,p}(\Sigma,\rel))\big{)}_{\theta,p} = \nonumber \\
	=W^{\theta,p}([0,T];W^{4(1-\theta),p}(\Sigma,\rel)),  \nonumber
\end{eqnarray} 	
for any finite $T>0$. Moreover, the proof of Lemma 3.3 in \cite{Spener} can be slightly adapted, in order to see that for 
$p \in (6,\infty)$ and for $\theta := \frac{1+\epsilon}{p} \in (0,\frac{1}{4})$ the general Sobolev-embedding from Theorem B.4 
in \cite{Spener} yields:
\begin{equation}  \label{Sobolev.embedding}		
	W^{\theta,p}([0,T];W^{4(1-\theta),p}(\Sigma,\rel))
	\hookrightarrow 
	C^{q_1}([0,T];C^{3,q_2}(\Sigma,\rel)),
\end{equation}
for any finite $T>0$ and for sufficiently small 
exponents $q_1,q_2\in (0,\frac{1}{8})$.
Hence, combining embeddings (\ref{embedding.interpol}) 
and (\ref{Sobolev.embedding}) with estimate (\ref{Lp.Lp.estimate}) - here with any fixed $p\in (6,\infty)$ - we obtain the existence of sufficiently small $q_1,q_2 \in (0,\frac{1}{8})$, such that the smooth solution 
$\phi_t = \langle N_t,\nu_{F^*} \rangle_{\rel^3}$ of 
equation (\ref{parabolic.equation}) satisfies:    
\begin{eqnarray}  \label{Super.estimate}
	\parallel \{\phi_t\} 
	\parallel_{C^{q_1}([0,T(\sigma)],C^{3,q_2}(\Sigma,\rel))} \leq C(\Sigma,F^*,q_1,q_2,p,T(\sigma),\sigma), \qquad
\end{eqnarray} 
for some appropriate constant 
$C=C(\Sigma,F^*,q_1,q_2,p,T(\sigma),\sigma)>0$. 
On account of the mean value theorem, estimate (\ref{Super.estimate}) particularly implies the estimate:
\begin{eqnarray}  \label{Super.estimate.2} 
	\parallel  B(\,\cdot\,,\phi_t,\ldots,D^3_x\phi_t)
	\parallel_{C^{\bar \mu,\frac{\bar \mu}{4}}(\Sigma \times [0,T(\sigma)],\rel)}                 \leq C(\Sigma,F^*,T(\sigma),\sigma,\bar \mu),                   
\end{eqnarray} 
for any small $\bar \mu \in (0,\frac{1}{8})$ with 
$\bar \mu <\min \{4\,q_1,q_2\}$ and for another appropriate 
constant $C=C(\Sigma,F^*,T(\sigma),\sigma,\bar \mu)$, which 
does not depend on any more data of the solution $\{\phi_t\}$,
especially not on the size of $\varepsilon$ from lines 
(\ref{smallness}) and (\ref{smallness.2}). 
Hence, we obtain from the parabolic 
Schauder a-priori estimates (\ref{Schauder.estimates}) - 
here with $l=0$ and $\mu =\bar \mu$ - combined with statement 
(\ref{leading.operator}), with conditions 
(\ref{smallness.3}) and (\ref{enclosure.T.2}) and with 
estimate (\ref{Super.estimate.2}), that the above 
smooth solution $\{\phi_t\}$ of equation 
(\ref{parabolic.equation}) satisfies:
\begin{equation}  \label{Schauder.estimates.normal.a}
	\parallel  \{\phi_t\} \parallel_{C^{4+\bar \mu,1+\frac{\bar \mu}{4}}(\Sigma \times [0,T(\sigma)],\rel)}
	\leq C_0(\Sigma,F^*,T(\sigma),\sigma,\bar \mu),
\end{equation} 
for some sufficiently small H\"older-exponent 
$\bar \mu \in (0,\frac{1}{8})$, where the above constant 
$C_0=C_0(\Sigma,F^*,T(\sigma),\sigma,\bar \mu)$ 
does not depend on the size of $\varepsilon>0$ from lines 
(\ref{smallness}) and (\ref{smallness.2}) neither. \\
Now, we recall that ``$k$'' in conditions  
(\ref{smallness}) and (\ref{smallness.2}) was a
fixed integer $\geq 4$. If $k=4$, then estimate 
(\ref{Schauder.estimates.normal.a}) does not have to be 
improved any more. But if $k>4$, then estimate (\ref{Schauder.estimates.normal.a}) should be used, in order to improve estimate (\ref{Super.estimate.2}) 
by means of another application of the mean value theorem:
\begin{eqnarray}  \label{Super.estimate.3} 
\parallel  B(\,\cdot\,,\phi_t,\ldots,D^3_x\phi_t)
\parallel_{C^{1+\bar \mu,\frac{1+\bar \mu}{4}}(\Sigma \times [0,T(\sigma)],\rel)}    
\leq \tilde C_1(\Sigma,F^*,T(\sigma),\sigma,\bar\mu),    
\end{eqnarray} 
for the same exponent $\bar \mu$ as in estimate 
(\ref{Schauder.estimates.normal.a})  
and for another appropriate constant 
$\tilde C_1=\tilde C_1(\Sigma,F^*,T(\sigma),\sigma,\bar \mu)$.
Since we have proved already, that the solution 
$\{\phi_t\}$ of equation (\ref{parabolic.equation}) is 
$C^{\infty}$-smooth, estimates (\ref{Schauder.estimates.normal.a}) and (\ref{Super.estimate.3}) can be combined again with conditions (\ref{smallness.3}) and (\ref{enclosure.T.2}), in order to infer from another application of Schauder estimates (\ref{Schauder.estimates}) 
- but now with $l=1$ and $\mu = \bar \mu$:
\begin{equation}  \label{Schauder.estimates.normal.b}
	\parallel \phi_t \parallel_{C^{5+\bar \mu,\frac{5+\bar \mu}{4}}
	(\Sigma \times [0,T(\sigma)],\rel)}   
	\leq  C_1(\Sigma,F^*,T(\sigma),\sigma,\bar \mu),
\end{equation}
for the same exponent $\bar \mu$ as in estimate 
(\ref{Schauder.estimates.normal.a})  
and for another appropriate constant 
$C_1= C_1(\Sigma,F^*,T(\sigma),\sigma,\bar \mu)$.
Hence, by finite induction - stopping after exactly $k-4$ steps 
on account of condition (\ref{smallness.3}) - we arrive 
in this way at the optimal Schauder estimate:  
\begin{equation}  \label{Schauder.estimates.normal.k-4}
	\parallel \phi_t \parallel_{C^{k+\bar \mu,\frac{k+\bar \mu}{4}}
	(\Sigma \times [0,T(\sigma)],\rel)}   
	\leq  C_{k-4}(\Sigma,F^*,T(\sigma),\sigma,\bar \mu),
\end{equation}
for the same exponent $\bar \mu$ as in estimate 
(\ref{Schauder.estimates.normal.a})  
and for another appropriate constant 
$C_{k-4}=C_{k-4}(\Sigma,F^*,T(\sigma),\sigma,\bar \mu)$,
which does not depend on the size of $\varepsilon>0$ from condititons (\ref{smallness.2}) and (\ref{smallness.3}) neither. Estimate (\ref{Schauder.estimates.normal.k-4}) 
immediately implies the estimate 
\begin{equation}  \label{Schauder.estimates.normal.k-4.b} 
	\parallel \tilde f_t - F^* 
	\parallel_{C^{k,\bar \mu}(\Sigma,\rel^3)}
	\equiv 	\parallel N_t  
	\parallel_{C^{k,\bar \mu}(\Sigma,\rel^3)}
	\leq \tilde C_{k-4}(\Sigma,F^*,T(\sigma),\sigma,\bar \mu)
	\qquad \forall \, t \in  [0,T(\sigma)],
\end{equation}
for the corresponding smooth solution 
$\tilde f_t = F^* + N_t \equiv F^* + \phi_t \, \nu_{F^*}$ 
of equation (\ref{parabolic.equation}), 
for some small exponent $\bar \mu \in (0,\frac{1}{8})$, 
which is the analogue here of formula (5.13) in 
\cite{Dall.Acqua.Spener.2016}. 
Now, on account of condition (\ref{enclosure.T}),  
due to the choice $\sigma <\min\{\delta,\varrho\}$
and since $F^*$ was supposed to be a $C^k$-local 
minimizer of the Willmore functional, we know that 
\begin{equation}  \label{bounded.below.W.F} 
	\Will(\tilde f_t)\geq \Will(F^*)
\end{equation}
for $t\in [0,T(\sigma)]$. Moreover, using the fact that 
the smooth family $\{\tilde f_t\}=\{F^* + \phi_t \,\nu_{F^*}\}$
solves equation (\ref{modified.Cauchy}), we can infer that 
\begin{eqnarray}  \label{monotonic.decrease} 
	\frac{d}{dt} \Will(\tilde f_t) 
	= \int_{\Sigma}  \, 
	\langle \partial_t^{\perp_{\tilde f_t}}(\tilde f_t),
	\nabla_{L^2}(\tilde f_t) \rangle_{\rel^3} \, d\mu_{\tilde f_t}
	= - \int_{\Sigma} \frac{1}{|A^0_{\tilde f_t}|^4} \, 
	|\nabla_{L^2}(\tilde f_t)|^2 \, d\mu_{\tilde f_t} \leq 0  
\end{eqnarray}
for $t\in [0,T(\sigma)]$, i.e. that 
$\Will(\tilde f_t)$ does not increase for $t\in [0,T(\sigma)]$. 
Moreover, due to $T(\sigma)<T_{\textnormal{max}}$ 
there holds equation (\ref{zusammengemogili}) on $[0,T(\sigma)]$, 
implying that we have
\begin{equation} \label{equal.energy} 
	\Will(\tilde f_t) = \Will(\PP(t,0,f_0)) \qquad \textnormal{for}\quad t\in [0,T(\sigma)].
\end{equation} 
Now, we can combine equations 
(\ref{bounded.below.W.F})--(\ref{equal.energy}) with 
Theorem 2.3 (ii) in \cite{Ruben.MIWF.III}, 
in order to apply the same argument as on page 360 in 
\cite{Chill.Schatz.2009}, ruling out the special 
case, in which there might 
hold $\Will(\tilde f_s) = \Will(F^*)$ for some 
$s\in [0,T(\sigma))$. Hence, in the sequel we can assume without loss of generality, that there holds 
$\Will(\tilde f_t) > \Will(F^*)$ 
for every $t\in [0,T(\sigma))$. Finally, we have to observe that condition (\ref{enclosure.T}) implies inequalities (\ref{delta.bounded}) and (\ref{no.umbilics}) for $t\in [0,T(\sigma)]$, since we chose 
$\sigma <\varrho$ and $\varrho<\tilde \delta$. Hence, 
there is some small constant $c=c(F^*,\sigma)>0$, such that 
$$ 
c(F^*,\sigma) \leq |A^0_{\tilde f_t}|^2 \leq  \frac{1}{c(F^*,\sigma)} 
\qquad \textnormal{for} \,\, t\in [0,T(\sigma)].
$$ 
We can therefore introduce the 
smooth, non-increasing and positive function
$[t \mapsto (\Will(\tilde f_t) - 
\Will(F^*))^{\theta}]$, for $t \in [0,T(\sigma))$, where $\theta=\theta(F^*) \in (0,1/2]$ denotes the exponent appearing in the Lojasiewicz-Simon-gradient-inequality for the Willmore-functional, Theorem 3.1 in \cite{Chill.Schatz.2009}, in order to compute by means of H\"older's inequality, again equation (\ref{modified.Cauchy}) 
and by the usual chain rule:
\begin{eqnarray}
	- \frac{d}{dt} (\Will(\tilde f_t) - \Will(F^*))^{\theta}
	= - \theta\, (\Will(\tilde f_t) - \Will(F^*))^{\theta-1}
	\, \int_{\Sigma} \langle \nabla_{L^2} \Will(\tilde f_t),
	\partial_t \tilde f_t \rangle \,d\mu_{\tilde f_t}                     \nonumber   \\
	= - \theta\, (\Will(\tilde f_t) - \Will(F^*))^{\theta-1}
	\, \int_{\Sigma} \langle \nabla_{L^2} \Will(\tilde f_t),
	\partial^{\perp_{\tilde f_t}}(\tilde f_t) \rangle \,d\mu_{\tilde f_t} \nonumber   \\   
	=\theta \,(\Will(\tilde f_t) - \Will(F^*))^{\theta-1} \,
	\int_{\Sigma}   \,\frac{1}{|A^0_{\tilde f_t}|^4} \, 
	|\nabla_{L^2}\Will(\tilde f_t)|^2 \, 
	d\mu_{\tilde f_t}          \nonumber    \\
	\geq \, c(F^*,\sigma)^2 \, \theta \,
	(\Will(\tilde f_t) - \Will(F^*))^{\theta-1} \,
	\, \int_{\Sigma} \,
	|\nabla_{L^2} \Will(\tilde f_t)|^2 \, d\mu_{\tilde f_t} \nonumber\\
	\geq \, c(F^*,\sigma)^4 \,\, \theta \,
	(\Will(\tilde f_t) - \Will(F^*))^{\theta-1} \,
	\Big{(} \int_{\Sigma}
	|\nabla_{L^2} \Will(\tilde f_t)|^2   \, d\mu_{\tilde f_t} \Big{)}^{1/2} \cdot \nonumber \\
	\cdot \Big{(} \int_{\Sigma} \,\frac{1}{|A^0_{\tilde f_t}|^8} \, 
	|\nabla_{L^2} \Will(\tilde f_t)|^2 \, d\mu_{\tilde f_t} \Big{)}^{1/2}               \nonumber\\
	\geq \frac{c(F^*,\sigma)^4\,\theta}{C_1^*(F^*)}  \, 
	\parallel \,\partial_t^{\perp_{\tilde f_t}}(\tilde f_t) 
	\parallel_{L^2(\mu_{\tilde f_t})}  \quad \textnormal{for} \,\,\,
	t \in [0,T(\sigma)), \qquad  \label{Lojasiewicz.Simon}
\end{eqnarray}
where we have been able to apply the Lojasiewicz-Simon-gradient inequality in line (\ref{Lojasiewicz.Simon}) in a $C^4$-ball of 
radius $\sigma$ about the Willmore immersion $F^*$ with 
appropriate constants $C^*_1=C^*_1(F^*)>0$ and 
$\theta=\theta(F^*) \in (0,1/2]$, taking estimate (\ref{enclosure.T}) for $t \in [0,T(\sigma)]$ into account. 
Now, estimate (\ref{enclosure.T}) also implies inequality (\ref{good.projection}) for $t\in [0,T(\sigma)]$, 
on account of $\sigma <\varrho<\tilde \delta$, and the time derivative $\partial_t \tilde f_t \equiv \partial_t N_t$ is actually a smooth section of the normal bundle of $F^*$, just as 
$N_t$ is, for every $t \in [0,T(\sigma)]$.   
Hence, we infer from an integration of inequality (\ref{Lojasiewicz.Simon})
w.r.t. time and again from estimate (\ref{enclosure.T}) -
implying inequality (\ref{good.projection}) to hold for 
$t\in [0,T(\sigma)]$ - the estimate:
\begin{eqnarray}  \label{Just.do.it}
	\int_0^s \parallel \partial_t \tilde f_t  \,
	\parallel_{L^2(\mu_{F^*})} \, dt
	\leq C(\sigma) \, 
	\int_0^s  \parallel \partial_t \tilde f_t  \,
	\parallel_{L^2(\mu_{\tilde f_t})} \, dt                     
	\leq  \nonumber \\
    \leq 2 \, C(\sigma) \, \int_0^s  \parallel
	\Big{(} \partial_t \tilde f_t \Big{)}^{\perp_{\tilde f_t}}
	\parallel_{L^2(\mu_{\tilde f_t})} \, dt  
	\leq  -\frac{2 \, C(\sigma) \,C^*_1}{c(F^*,\sigma)^4\,\theta}
	\int_0^s \frac{d}{dt} 
	\Big{(} (\Will(\tilde f_t) - 
	\Will(F^*))^{\theta} \Big{)} \, dt = \nonumber \\
	= \frac{2 \, C(\sigma) \,C^*_1}{c(F^*,\sigma)^4\,\theta} \,
	\Big{(} (\Will(\tilde f_0) - \Will(F^*))^{\theta}
	- (\Will(\tilde f_s) - \Will(F^*))^{\theta} \Big{)} 
	\leq \nonumber \\
	\leq \frac{2 \, C(\sigma) \,C^*_1}{c(F^*,\sigma)^4\,\theta} \,
	\Big{(} \Will(F^* + N_0) - 
	\Will(F^*) \Big{)}^{\theta} <\infty \quad  \forall \,s \in [0,T(\sigma)).  \qquad
\end{eqnarray}
Now, we can derive from a combination of condition (\ref{enclosure.T}) and estimate (\ref{Just.do.it}), together 
with the triangle inequality for the $L^2(\mu_{F^*})$-norm:
\begin{eqnarray} \label{estim.L2}
	\parallel \tilde f_s - F^*
	\parallel_{L^2(\mu_{F^*})}
	\leq \parallel N_0 \parallel_{L^2(\mu_{F^*})}
	+ \frac{2\, C(\sigma)\,C^*_1}{c(F^*,\sigma)^4\,\theta}
	\, (\Will(F^*+ N_0) - \Will(F^*))^{\theta} 
	\leq   \nonumber \\
	\leq C \, \parallel N_0 \parallel_{C^2(\Sigma,\rel^3)}^{\theta}
	\qquad  \forall \,s \in [0,T(\sigma)],          \qquad
\end{eqnarray}
for some appropriate constant $C=C(\Sigma,F^*,\theta,\sigma)>0$.
By Theorem 6.4.5 (iii) in \cite{Bergh.Loefstroem.1976} 
we can interpolate the Besov space   
$B_{p^*,p^*}^{\beta\,(k+\bar \mu)}(\Sigma,\rel^3)$, 
for $p^*=\frac{2}{1-\beta}>>1$ and $\beta \in (0,1)$ 
close to $1$, between the spaces 
$C^{k,\bar\mu}(\Sigma,\rel^3)
=B_{\infty,\infty}^{k+\bar \mu}(\Sigma,\rel^3)$ 
and $L^2(\Sigma,\rel^3)=B^0_{2,2}(\Sigma,\rel^3)$, 
and we can then use the fact that 
$B_{p^*,p^*}^{\beta\,(k+\bar \mu)}(\Sigma,\rel^3)$ 
embeds into $C^{k}(\Sigma,\rel^3)$ by 
the fractional Sobolev embedding theorem, provided 
there holds $\beta\,(k+\bar \mu) - \frac{2}{p^*} \equiv  
\beta\,(k+\bar \mu)+\beta-1>k$. 
Consequently, we infer from estimates (\ref{smallness.2}),
(\ref{Schauder.estimates.normal.k-4.b}) and (\ref{estim.L2}):
\begin{eqnarray} \label{interpolari}
	\parallel \tilde f_s - F^* \parallel_{C^{k}(\Sigma,\rel^3)}	
	\leq C\, \parallel \tilde f_s - F^* 
	\parallel_{C^{k,\bar \mu}(\Sigma,\rel^3)}^{\beta}
	\, \parallel \tilde f_s - F^* \parallel_{L^2(\mu_{F^*})}^{1-\beta} \leq \quad  \\
	\leq C\, \parallel N_0 \parallel_{C^2(\Sigma,\rel^3)}^{(1-\beta) \,\theta}
	\leq C^* \,(C^{o}(\varepsilon))^{(1-\beta) \,\theta} \quad  \forall \,s \in [0,T(\sigma)],    \nonumber
\end{eqnarray}
for some appropriately large constant
$C^*=C^*(\Sigma,F^*,\bar \mu,k,\beta,\theta,T(\sigma),\sigma)$, 
which is independent of $\varepsilon$.
It therefore turns out now, that we should choose $\varepsilon>0$ 
above in estimate \eqref{smallness} that small, such that
\begin{eqnarray}\label{final.contradiction}
C^* \, (C^{o}(\varepsilon))^{(1-\beta) \,\theta} 
\leq \frac{\sigma}{2}
\end{eqnarray}
holds, implying by estimate (\ref{interpolari}) that 
we thus had:  
$$ 
\parallel \tilde f_t - F^* \parallel_{C^{k}(\Sigma,\rel^3)}\lfloor_{t=T(\sigma)}
\leq \frac{\sigma}{2}.
$$ 
But this contradicts the fact that  
$\parallel \tilde f_t - F^* \parallel_{C^{k}(\Sigma,\rel^3)}\lfloor_{t=T(\sigma)}
=\sigma$ would have to hold at time $t=T(\sigma)$ 
on account of condition (\ref{enclosure.T}), 
if the ``maximal time'' $T(\sigma)$ in (\ref{enclosure.T}) 
would have actually been finite 
and also smaller than $T_{\textnormal{max}}$.
In the remaining special case 
``$T(\sigma)=T_{\textnormal{max}}<\infty$'' we could infer, 
that estimates (\ref{enclosure.T}), (\ref{enclosure.T.2}) and (\ref{Schauder.estimates.normal.k-4.b}) would hold on every 
compact interval $[0,T]$ with $T<T_{\textnormal{max}}<\infty$.
Since we also know that $\sigma <\varrho$, a comparison 
of conditions (\ref{smaller.than.varrho}) and
(\ref{enclosure.T.2}) shows us, that in this situation 
the smooth solution $\{\phi_t\}$ of 
the quasilinear parabolic equation (\ref{parabolic.equation}) 
could be extended - using e.g. the methods of Theorems 2 and 3 in \cite{Jakob_Moebius_2016} - from 
$\Sigma \times [0,T_{\textnormal{max}})$ to 
$\Sigma \times [0,T']$, for some $T'>T_{\textnormal{max}}$, 
of class $C^{4+\mu,1+\frac{\mu}{4}}
(\Sigma \times [0,T'],\rel)$, for any fixed 
$\mu \in (0,\bar \mu)$, and thus also of class 
$C^{\infty}(\Sigma \times [0,T'],\rel)$ on account 
of the above bootstrap argument employing estimates 
(\ref{Schauder.estimates}), but without violating condition (\ref{smaller.than.varrho}) for $t\in [0,T']$. 
This contradicts the definition  
of the maximal time $T_{\textnormal{max}}$.  
Hence, we have proved, that there actually has 
to hold ``$T(\sigma)=\infty$'' in estimates 
(\ref{enclosure.T}) and (\ref{enclosure.T.2}), 
i.e. that the particular smooth solution 
$\tilde f_t = F^* + N_t$ of equation (\ref{modified.Cauchy}) 
exists globally and satisfies the smallness condition 
(\ref{enclosure.T}) at arbitrarily large times $t$:
\begin{equation}  \label{enclosure.T.infty}
	\parallel \tilde f_{t} - F^* \parallel_{C^{k}(\Sigma,\rel^3)}
	\leq \sigma  \qquad \forall \,t\in [0,\infty),
\end{equation}
provided the initial immersion $\tilde f_{0}$ satisfies
condition (\ref{smallness.2}) with $\varepsilon>0$ chosen that small in condition (\ref{smallness}),
such that inequality (\ref{final.contradiction}) finally holds. Since we also know now that $T_{\textnormal{max}}=\infty$, we obtain equation (\ref{zusammengemogili}) for every $t\geq 0$, 
i.e.: there is a smooth family of smooth diffeomorphisms 
$\Psi_t:\Sigma \stackrel{\cong}\longrightarrow \Sigma$, 
$\Psi_0 = \textnormal{Id}_{\Sigma}$, such that: 
\begin{equation}  \label{zusammengemogili.infty}
	\tilde f_t = \PP(t,0,f_0)\circ \Phi_0 \circ \Psi_t^{-1} 
	\qquad \textnormal{for} \,\, t\in [0,\infty). 
\end{equation} 
Now, having chosen $\varepsilon>0$ in 
(\ref{smallness}) sufficiently small, we can now let tend 
$s \to \infty$ in estimate (\ref{Just.do.it}) and obtain:
\begin{equation}   \label{finite.improper.int} 
\int_0^{\infty} \parallel \partial_t \tilde f_t
\parallel_{L^2(\mu_{F^*})} \, dt
\leq \frac{2 \, C(\sigma) \,C^*_1}{c(F^*,\sigma)^4\,\theta}\,
\Big{(}\Will(F^*+N_0)-\Will(F^*) \Big{)}^{\theta} <\infty,
\end{equation}
implying the existence of a unique function
$F_{\infty} \in L^2((\Sigma,\mu_{F^*}),\rel^3)$, 
such that
\begin{equation}  \label{full.L2.convergence}
	\tilde f_t \equiv F^* + N_t  
	\longrightarrow F_{\infty} \qquad
	\textnormal{in} \quad L^2(\Sigma,\mu_{F^*})
\end{equation}
as $t \to \infty$. Combining now the full convergence (\ref{full.L2.convergence}) and equation (\ref{zusammengemogili.infty}) with 
estimate (\ref{enclosure.T.infty}), we can easily infer: 
\begin{equation}  \label{full.Ck-1.convergence}
	\PP(t,0,f_0)\circ \Phi_0 \circ \Psi_t^{-1} 
	= \tilde f_t \longrightarrow F_{\infty}  \qquad
	\textnormal{in} \quad C^{k-1,\alpha}(\Sigma,\rel^3)
\end{equation}
as $t \to \infty$, where the limit function 
$F_{\infty}$ turns out to be a umbilic-free 
$C^k$-immersion, because it additionally satisfies estimate  (\ref{enclosure.T.infty}) in the limit, i.e.: 
\begin{equation}  \label{enclosure.infty}
	\parallel F_{\infty} - F^* \parallel_{C^{k}(\Sigma,\rel^3)}
	\leq \sigma,
\end{equation}
and since estimate (\ref{enclosure.T.infty}) implies inequalities (\ref{delta.bounded}) and (\ref{no.umbilics}) to hold for every 
$t\geq 0$. It remains to prove, that the limit immersion $F_{\infty}$ of convergence (\ref{full.Ck-1.convergence}) is ``Willmore'' 
and a $C^k$-local minimizer of the Willmore functional $\Will$. To this end, we have to distinguish two different cases: \\  
\underline{Case 1}: If $k>4$, then we infer from 
convergences (\ref{finite.improper.int})
and (\ref{full.Ck-1.convergence}) immediately:
\begin{eqnarray*} 
	0 \longleftarrow 
	\parallel \,\partial_t^{\perp_{\tilde f_t}}(\tilde f_t) 
	\parallel_{L^2(\mu_{\tilde f_t})}^2\lfloor_{t=t_i}
	=  \int_{\Sigma}   \,\frac{1}{|A^0_{\tilde f_t}|^8} \, 
	|\nabla_{L^2}\Will(\tilde f_t)|^2 \, d\mu_{\tilde f_t}\lfloor_{t=t_i}\\
	\longrightarrow 
	\int_{\Sigma}  \,\frac{1}{|A^0_{F_{\infty}}|^8} \, 
	|\nabla_{L^2}\Will(F_{\infty})|^2 \, d\mu_{F_{\infty}},
\end{eqnarray*}
for some appropriate sequence $t_i\nearrow \infty$, 
showing that $F_{\infty}$ is a umbilic-free Willmore immersion, satisfying statement (\ref{enclosure.infty}) for some $k>4$. \\
\underline{Case 2}: If $k=4$, then we cannot conclude directly, 
since the convergence in (\ref{full.Ck-1.convergence}) 
is not strong enough. However, since the immersions $\tilde f_t$ are known to be $C^{\infty}$-smooth on $\Sigma$, we can apply here Theorem 1.2 of Bernard's paper \cite{Bernard.2016} and obtain for any test function $\xi \in C^{\infty}(\Sigma,\rel^3)$: 
\begin{eqnarray}   \label{Bernard.Riviere}
	\int_{\Sigma}  \,\Big{\langle} \frac{1}{|A^0_{\tilde f_t}|^4} \, 
	\nabla_{L^2}\Will(\tilde f_t),\xi \Big{\rangle}_{\rel^3}  \, d\mu_{\tilde f_t}  
	= \sum_{j=1}^2 \int_{\Sigma} \,\Big{\langle} \nabla^{\tilde f_t}_{j} 
	\Big{(}\frac{1}{|A^0_{\tilde f_t}|^4} \, \xi \Big{)}, 
	\Tan^j(\tilde f_t)  \Big{\rangle}_{\rel^3} \,d\mu_{\tilde f_t} 
\end{eqnarray} 
for two particular, globally defined differential operators 
$\Tan^j:C^{\infty}_{\textnormal{Imm}}(\Sigma,\rel^3) 
\longrightarrow  C^{\infty}(\Sigma,\rel^3)$ of third order, 
which are explicitly given in Theorem 1.2 of \cite{Bernard.2016}. 
Combining equation (\ref{Bernard.Riviere}) with equation 
(\ref{modified.Cauchy}) and convergences (\ref{finite.improper.int})
and (\ref{full.Ck-1.convergence}) for $k=4$, we obtain here 
along the same sequence $t_i\nearrow \infty$ as in Case 1:
\begin{eqnarray} 
	0 \longleftarrow 
	\int_{\Sigma}  \,\Big{\langle} \partial_t^{\perp_{\tilde f_t}}(\tilde f_t),\xi \Big{\rangle}_{\rel^3}  \, d\mu_{\tilde f_t}\lfloor_{t=t_i} 
	=- \sum_{j=1}^2 \int_{\Sigma}  \,\Big{\langle} \nabla^{\tilde f_t}_{j} 
	\Big{(}\frac{1}{|A^0_{\tilde f_t}|^4} \, \xi \Big{)}, 
	\Tan^j(\tilde f_t) \Big{\rangle}_{\rel^3}  \, d\mu_{\tilde f_t}
	\lfloor_{t=t_i}  
    \nonumber    \\ \label{Bernard.action} 
	\longrightarrow  
	-\sum_{j=1}^2 \int_{\Sigma}  \,\Big{\langle} \nabla^{F_{\infty}}_{j} 
	\Big{(}\frac{1}{|A^0_{F_{\infty}}|^4} \, \xi \Big{)}, 
	\Tan^j(F_{\infty}) \Big{\rangle}_{\rel^3}  \, d\mu_{F_{\infty}} \quad 
\end{eqnarray}  
as $i\to \infty$. Since statement \eqref{enclosure.infty} holds here for $k=4$, we can again integrate by parts in \eqref{Bernard.action} 
and infer again from Theorem 1.2 in \cite{Bernard.2016}, that 
$F_{\infty}$ is a umbilic-free Willmore immersion
of class $C^k$, satisfying statement \eqref{enclosure.infty} for $k=4$. Moreover, combining statement \eqref{enclosure.infty} again with the Lojasiewicz-Simon-gradient-inequality, 
Theorem 3.1 in \cite{Chill.Schatz.2009}, 
$F_{\infty}$ turns out to satisfy $\Will(F_{\infty})=\Will(F^*)$, 
proving that $F_{\infty}$ is actually a $C^k$-local 
minimizer of $\Will$ as well, 
for any fixed $k\geq 4$.
\qed    \\
\noindent
\begin{remark}\,
\begin{itemize} 
\item[1)]
	It should be noted here, that the proof of Lemma 4.1 
	in \cite{Chill.Schatz.2009} was partially inspired by Simon's proof of Theorem 2 of his pioneering work 
	\cite{Simon.1983} on parabolic $L^2$-gradient flows 
	corresponding to real analytic functionals; see 
	the respective quotation on p. 359 in \cite{Chill.Schatz.2009}.  
	However, in the proof of the above Theorem 
	\ref{Convergence.to.local.minimizer}
	we actually avoid most of Simon's arguments, besides 
	the decisive usage of the ``Lojasiewicz-Simon-gradient-inequality''. 
	We especially do not use the assertion, that the Fr\'echet derivative of the quasilinear differential operator 
	$$
	\eta \mapsto \frac{1}{2} \, \frac{1}{|A^0_{F^*+ \eta \nu_{F^*}}|^4} \,g^{ij}_{F^*+ \eta \nu_{F^*}} \, 
	g^{kl}_{F^*+ \eta \nu_{F^*}} \, \nabla^{F^*}_{ijkl}(\eta)
	- B(\,\cdot\,,\eta,\ldots,D^3_x\eta)
	$$
	appearing in equation (\ref{parabolic.equation}) - 
	mapping $C^{4,\alpha}(\Sigma,\rel)$ into $C^{0,\alpha}(\Sigma,\rel)$ - 
	in the function $\eta \equiv 0$ \underline{was symmetric} w.r.t. an appropriately chosen  
	$L^2(\Sigma)$-scalar product. Such a symmetry-property of this Fr\'echet derivative - see here formulae (0.1), (1.5) and  
	(1.6) in Simon's paper \cite{Simon.1983} - 
	is a crucial point in the proof of certain a-priori-estimates in Section 4 of \cite{Simon.1983}, where key-consequences of the classical ``Spectral Theorem'' for selfadjoint 
	compact linear operators between Banach spaces 
	are utilized. See here also formulae (1.20) 
	and (4.10)--(4.16) in \cite{Simon.1983}.   
	\item[2)]  In contrast to the final steps of the 
	proof of Lemma 4.1 in \cite{Chill.Schatz.2009}, we cannot combine neither estimate \eqref{enclosure.T.infty} nor the full $C^{k-1,\alpha}$-convergence in \eqref{full.Ck-1.convergence} 
	with ``localized $L^{\infty}$-estimates''
	for covariant derivatives $\nabla^m A_{\tilde f_t}$ 
	of the second fundamental forms
	of the converging immersions $\tilde f_t$ in (\ref{full.L2.convergence}) and (\ref{full.Ck-1.convergence}), 
	in order to improve the quality of convergence in (\ref{full.Ck-1.convergence}) furthermore, because such strong estimates have not been proved so far for flow lines of the MIWF, only for flow lines of the classical Willmore flow; 
	see Section 5 in \cite{Kuwert.Schaetzle.2001} and 
	Section 4 in \cite{Kuwert.Schaetzle.2002}.
	\qed 
\end{itemize}        
\end{remark} 
\noindent
\underline{Proof of Theorem \ref{main.result.2}}: \\\\ 
The proof of Theorem \ref{main.result.2} essentially follows the lines of the proof of Theorem \ref{main.result.1}, i.e. it relies on a combination of both parts of Theorem \ref{Frechensbergo} and on our second full convergence theorem, Theorem \ref{Convergence.to.local.minimizer}. 
The only striking point here is the fact, that the evolution operator 
$$ 
\PP^*(\,\cdot\,,0,\,\cdot\,):
B_{\rho}^{4,\gamma}(F_0) \subset C^{4,\gamma}(\Sigma,\rel^3) \longrightarrow C^{4+\gamma,1+\frac{\gamma}{4}}(\Sigma \times [0,T],\rel^3) 
$$
from line (\ref{solution.operator.2}) for the 
modified MIWF-equation (\ref{de_Turck_equation_2}) 
is a non-linear operator of class $C^1$, for any fixed 
$\gamma \in (0,1)$, which is exactly ``compatible'' with the requirement in Theorem \ref{Convergence.to.local.minimizer} 
- to be applied here with $k=4$ - to let the flow line $\PP(\,\cdot\,,0,f_0)$ 
of the MIWF start moving in any smooth immersion 
$f_0:\Sigma \longrightarrow \rel^3$ which  
is sufficiently close to the given umbilic-free 
$C^4$-local minimizer $F^*$ of the Willmore energy $\Will$ 
in the $C^{4,\alpha}(\Sigma,\rel^3)$-norm, for any 
fixed $\alpha \in (0,1)$, e.g. for $\alpha=\gamma$. Therefore, both parts of Theorem \ref{Frechensbergo} can be combined here with the statement of Theorem \ref{Convergence.to.local.minimizer}, and the proof of Theorem \ref{main.result.2} follows the lines of the 
proof of Theorem \ref{main.result.1}. 
\qed  \\\\
\noindent 
\underline{Acknowledgements}\\\\
	The M\"obius-invariant Willmore flow was originally motivated by Professor Ben Andrews and Professor Dr. Reiner M. Sch\"atzle.
	The author would like to express his deep gratitude to both of them, in particular for many invaluable comments and constant support. Moreover, the author would like to thank Professor Gieri Simonett for having directed the author's attention to some of his papers about the surface diffusion flow and center manifold techniques. Finally, the author would like to thank Professor Itai Shafrir and Professor Yehuda Pinchover for their hospitality and strong support at the Mathematics department of the ``Israel Institute of Technology''.

\end{document}